\documentclass[11pt]{article}

\usepackage{mathrsfs}
\usepackage{amsbsy,amsmath,amsthm,amssymb,graphicx,subcaption}
\usepackage{algorithm}
\usepackage{algpseudocode}
\usepackage{verbatim}
\usepackage{natbib}
\usepackage{hyperref}
\usepackage{makecell,booktabs}
\usepackage{xcolor}
\usepackage[normalem]{ulem}
\usepackage{dsfont}
\usepackage{enumitem}
\usepackage{geometry}

\newtheorem{theorem}{Theorem}
\newtheorem{lemma}[theorem]{Lemma}
\newtheorem{proposition}[theorem]{Proposition}
\newtheorem{remark}[theorem]{Remark}
\newtheorem{corollary}[theorem]{Corollary}

\title{\bf On spectral gap decomposition for Markov chains} 
\author{Qian Qin \\ School of Statistics \\ University of Minnestoa}
\date{}

\makeatletter
\def\namedlabel#1#2{\begingroup
	#2%
	\def\@currentlabel{#2}%
	\phantomsection\label{#1}\endgroup
}
\makeatother

\newcommand{\df}{\mathrm{d}}
\newcommand{\X}{\mathsf{X}}
\newcommand{\Y}{\mathsf{Y}}
\newcommand{\Z}{\mathsf{Z}}
\newcommand{\C}{\mathcal{C}}
\newcommand{\B}{\mathcal{B}}
\newcommand{\A}{\mathcal{A}}

\newcommand{\ind}{\mathds{1}}
\newcommand{\Proj}{E}
\newcommand{\Projj}{F}

\DeclareMathOperator*{\essinf}{ess\,inf}
\DeclareMathOperator*{\esssup}{ess\,sup}

\begin{document}
	\maketitle
	
	\begin{abstract}
		Multiple works regarding convergence analysis of Markov chains have led to spectral gap decomposition formulas of the form  
		\[
		\mathrm{Gap}(S) \geq c_0 \left[\essinf_z \mathrm{Gap}(Q_z)\right] \mathrm{Gap}(\bar{S}),
		\]
		where $c_0$ is a constant, $\mathrm{Gap}$ denotes the right spectral gap of a reversible Markov operator, $S$ is the Markov transition kernel (Mtk) of interest, $\bar{S}$ is an idealized or simplified version of $S$, and $\{Q_z\}$ is a collection of Mtks characterizing the differences between $S$ and $\bar{S}$.  
		
		This type of relationship has been established in various contexts, including:  
		1. decomposition of Markov chains based on a finite cover of the state space,  
		2. hybrid Gibbs samplers, and  
		3. spectral independence and localization schemes.  
		
		We show that multiple key decomposition results across these domains can be connected within a unified framework, rooted in a simple sandwich structure of~$S$.  
		Within the general framework, we establish new instances of spectral gap decomposition for hybrid hit-and-run samplers and hybrid data augmentation algorithms with two intractable conditional distributions.  
		Additionally, we explore several other properties of the sandwich structure, and derive extensions of the spectral gap decomposition formula.
	\end{abstract}

	\section{Introduction}
	
	Convergence analysis for Markov chains is a fundamental research area in probability theory, statistics, and computer science.
	In particular, convergence analysis plays an important role in the study of Markov chain Monte Carlo (MCMC) algorithms \citep{jones2001honest}.
	
	A common strategy is to compare the behavior of different Markov chains. 
	{ This literature includes Markov-operator and Dirichlet-form comparisons \citep{peskun1973optimum,diaconis1993comparison,roberts1997geometric,hobert2008theoretical,pozza2024fundamental}, perturbation bounds \citep{rudolf2018perturbation,caprio2023calculus}, comparisons of scan strategies in component-wise MCMC algorithms \citep{he2016scan,qin2020convergence,chlebicka2023solidarity,gaitonde2024comparison}, comparison between exact and approximate algorithms \citep{roberts1997geometric,jones2014convergence,pillai2014ergodicity,wang2022theoretical,power2024weak}, among others.}
	This article aims to unify and extend a series of comparison results {  involving the decomposition of Markovian dynamics} that are scattered throughout the literature. 
	Our focus is on results exhibiting a property we call ``spectral gap decomposition," expressed as follows:
	\begin{equation} \label{ine:generic}
		\text{Gap}(S) \geq c_0 \left[\essinf_z \text{Gap}(Q_z)\right] \text{Gap}(\bar{S}).
	\end{equation}
	Here, $c_0$ is a positive constant, and $\text{Gap}(\cdot)$ represents the right spectral gap of the Markov transition kernel (Mtk) of a reversible Markov chain, with a larger gap loosely translating to faster mixing.
	In this framework:
	\begin{itemize}
		\item $S$ is the Mtk of interest.
		\item $\bar{S}$ is an idealized or simplified version of $S$, serving as a benchmark  for comparison.
		\item $\{Q_z\}$ is a collection of Mtks characterizing the differences between $S$ and $\bar{S}$.
	\end{itemize}
	{ In \eqref{ine:generic} and the analogous displays below, the essential infimum is taken with respect to a mixing measure on the index~$z$, denoted later by~$\varpi$.}
	We will also consider something more delicate, i.e., decomposition formulas of the form $\mathcal{E}_S(f) \geq c_0 \, [\essinf_z \mathrm{Gap}(Q_z)] \, \mathcal{E}_{\bar{S}}(f)$, where $\mathcal{E}_S(f)$ is the Dirichlet form of a function~$f$ associated with~$S$ (see Section~\ref{sec:preliminaries}), and $\mathcal{E}_{\bar{S}}(f)$ is defined similarly.

	
	Some of the first instances of spectral gap decomposition were derived in the works of \cite{caracciolo1992two} and \cite{madras2002markov}.
	These authors investigated the decomposition of Markov chains based on a finite covering of the state space.
	To be more specific, assume that $S$ is the Mtk of a reversible chain defined on a state space $\X$, where $\X$ is partitioned into (or covered by) a finite number of subsets, say, $\X = \bigcup_z \X_z$.
	Then, under mild conditions, \eqref{ine:generic} holds, with $\bar{S}$ characterizing the chain's transition between the subsets, and $Q_z$ depicting the chain's movement within a subset $\X_z$.
	This type of relation has been applied to Metropolis-Hastings algorithms \citep{guan2007small}, tempering algorithms \citep{woodard2009conditions} and the reversible jump algorithm \citep{qin2023geometric}.
	See \cite{madras1996factoring}, \cite{jerrum2004elementary}, and \cite{ge2018simulated} for alternative versions of spectral gap decomposition with forms different from~\eqref{ine:generic}.
	See also \cite{atchade2021approximate} for a decomposition result involving approximate spectral gaps.

	Almost independently of the aforementioned works, another line of research discovered spectral gap decomposition for a number of hybrid Gibbs algorithms.
	In \cite{andrieu2018uniform}, it was found that \eqref{ine:generic} holds when~$S$ is a hybrid data augmentation (two-component Gibbs) algorithm, with one conditional distribution replaced by a Markovian (e.g., Metropolis-Hastings) approximation.
	In this context, $\bar{S}$ is the Mtk of the ideal data augmentation algorithm, and $\{Q_z\}$ correspond to the Markovian approximations.
	Another spectral gap decomposition result was derived for random-scan hybrid Gibbs samplers in \cite{qin2024spectral}; see also \cite{ascolani2024scalability}.
	For studies of similar flavors involving slice samplers, which are a subclass of data augmentation algorithms, see \cite{latuszynski2014convergence,power2024weak}.
	
	
	A third relevant line of research stems from theoretical computer science.
	In recent years, a technique called ``spectral independence" has found success in the analysis of some Markov chains, particularly Glauber dynamics \citep{anari2021spectral,chen2021rapid,feng2022rapid,chen2023rapid,qin2024spectrala,chen2024rapid}.
	See also \cite{carlen2003determination} for an earlier work with similar flavors.
	In \cite{chen2022localization}, the spectral independence technique is further developed into a class of methods called ``localization schemes."
	As we will see, in the framework of \cite{chen2022localization}, a key technical result involving variance conservation can be viewed as spectral gap decomposition.
	 
	While the three lines of works seem to share conceptual similarities, as noted in, e.g., \cite{chen2022localization} and \cite{liu2023spectral}, no comprehensive effort has been made to connect them into a unified framework.
	The goal of this paper is to fill this gap.
	It is shown that spectral gap decomposition can be derived for all Mtks possessing a simple sandwich structure, which is rigorously described in in Section~\ref{sec:unified}.
	Roughly speaking, $S$ has the sandwich structure when it can be written into the form
	\[
	S(x, \mathsf{A}) = \int P^*(x, \df (y,z)) \, Q_z(y, \df y') \, P((y',z), \mathsf{A}),
	\]
	where $P$ defines a contraction, and $P^*$ is in some sense the adjoint of~$P$.
	(The actual structure is slightly more general than this.)
	Then \eqref{ine:generic} holds with 
	\[
	\bar{S}(x, \mathsf{A}) = \int P^*(x, \df (y,z)) \, \varpi_z(\df y') \, P((y',z), \mathsf{A}),
	\]
	where $\varpi_z$ is the stationary distribution of $Q_z$.
	We demonstrate that multiple spectral gap decomposition results from the three lines of works can find their roots in this structure.
	This includes key results from \cite{caracciolo1992two}, \cite{madras2002markov}, \cite{andrieu2018uniform}, \cite{chen2022localization}, and \cite{qin2024spectral}.
	
	{ The sandwich representation may be viewed as a reversible
		disintegration of~$S$ on a refined space after a congruence
		transformation. It arises most naturally for Markov dynamics built
		from local or low-dimensional updates. The present framework does not,
		in general, provide an algorithm-specific disintegration of
		gradient-based methods such as MALA, HMC, or underdamped Langevin
		algorithms. Nevertheless, when the relevant assumptions hold, these
		chains may still be treated using the partition- or cover-based
		decompositions in \cite{madras2002markov}; they may also serve as the inner kernels~$Q_z$.
		In either case, quantitative bounds for gradient-based methods typically require separate contraction, functional-inequality, or
		hypocoercivity arguments.}

	We discover two new instances of spectral gap decomposition that fall into the unified framework:
	\begin{enumerate}
		\item Spectral gap decomposition holds when $S$ is a hybrid hit-and-run sampler, and $\bar{S}$ is its idealized counterpart.
		\item Spectral gap decomposition holds when $S$ is a hybrid data augmentation algorithm with two intractable conditional distributions, as opposed to one in the study of \cite{andrieu2018uniform}.
	\end{enumerate}

	Additionally, a unified framework allows one to extend spectral gap decomposition in a generic setting.
	For instance, while most existing spectral gap decomposition results concern reversible Markov chains, with the exception of \cite{qin2023geometric}, it is almost trivial to derive a similar generic result for non-reversible chains based on the sandwich structure.
	The new spectral gap decomposition result for hybrid data augmentation algorithm with two intractable conditional distributions mentioned above is derived in this way.
	Another extension involves weak Poincar\'{e} inequalities, which is a weaker notion than admitting a right spectral gap, and closely related to subgeometric convergence \citep{andrieu2022comparison}.
	It is shown that, under the sandwich structure, weak Poincar\'{e} inequalities involving $\bar{S}$ and the $Q_z$'s can lead to a weak Poincar\'{e} inequality for~$S$.
	This generalizes a key result from \cite{power2024weak} regarding hybrid slice samplers.

	To obtain the sandwich structure, we take inspiration from the pioneering work of \cite{caracciolo1992two}.
	The structure is closely related to the operator factorization used by \cite{rudolf2018comparison}, who establish Peskun orderings among particular hit-and-run, slice-sampling, and random-walk Metropolis kernels through an algorithm-specific common lifting.  Our formulation instead uses a measurable family of inner kernels $Q_z$ to obtain quantitative spectral-gap decomposition and its norm, non-reversible, and weak-Poincar\'{e} extensions.
	See also \cite{andersen2007hit}, which unified hit-and-run samplers with certain Gibbs and data augmentation algorithms.

	While the sandwich framework under in Section \ref{sec:unified}
	unifies most existing results of the form \eqref{ine:generic} that
	we know of, it does not encompass all variants. Notable qualifications
	concern the bounds in \cite{jerrum2004elementary} and
	\cite{ge2018simulated}.
	{  The result in \cite{jerrum2004elementary} does not strictly adhere to the structure \eqref{ine:generic}; it involves, in the context of Markov chain decomposition over a partition of the state space, an additional parameter quantifying the point-wise maximum amount of flow out of a subset $\X_z$, not captured by $\bar{S}$.}
	{  Ge et al.~\citep{ge2018simulated} treat soft-mixture decompositions that need not admit an exact sandwich factorization. We show at the end of Section~\ref{sec:unified} and in Appendix~\ref{app:beyond} that our framework can be extended by replacing the exact factorization with the weaker requirement that the Dirichlet form of~$S$ dominates that of a sandwich-structured operator. This extension recovers a discrete-time analogue of Theorem~6.1 in that work. It does not, however, encompass their more specialized decomposition bounds for tempering algorithms.

	We obtain two new quantitative applications of the framework. First, combining our hybrid hit-and-run decomposition result with quantitative bounds for
	one-dimensional random-walk Metropolis kernels and exact hit-and-run,
	we obtain, to our knowledge, the first explicit spectral-gap and $L^2$-mixing-time bounds for a multiscale
	Metropolis-within-hit-and-run algorithm targeting a strongly
	log-concave and log-smooth distribution. 
	Section~\ref{sec:example-hybrid}
	develops this application and illustrates it using Bayesian logistic
	regression.
	
	Second, we apply our new non-reversible comparison result to a hybrid
	data augmentation algorithm in which both conditional distributions
	are intractable and are approximated by Metropolis--Hastings kernels.
	For strongly log-concave and log-smooth targets, we derive explicit
	operator-norm contraction and $L^2$ mixing-time bounds. This application
	is developed in Section~\ref{sec:da-example}.
	}

	The rest of this article is organized as follows.
	After introducing some preliminary facts in Section~\ref{sec:preliminaries}, we list existing and new examples of spectral gap decomposition in Section~\ref{sec:examples}.
	We build a unified framework around the aforementioned sandwich structure in Section~\ref{sec:unified}.
	In this framework, we establish decomposition formulas involving Dirichlet forms, spectral gaps, and norms.
	In Section~\ref{sec:special} and Appendix \ref{app:derive}, we demonstrate that the decomposition results from Section~\ref{sec:examples} can be derived within the unified framework.
	Some other simple consequences of the sandwich structure are derived in Section~\ref{sec:simple}.
	In Sections~\ref{sec:example-hybrid} and \ref{sec:da-example}, we demonstrate how one can obtain concrete mixing time bounds using the new decomposition results.
	Some technical results and proofs are given in Appendices~\ref{app:lemma} and~\ref{app:proofs}.
	Appendix \ref{app:beyond} demonstrates spectral gap decomposition beyond the sandwich structure.

	\section{Preliminaries} \label{sec:preliminaries}
	
	\subsection{$L^2$ space with respect to a probability measure}
	
	Let $(\Omega, \mathcal{F}, \rho)$ be a probability space.
	Define a Hilbert space $L^2(\rho)$ of real functions on $\Omega$ that are square integrable with respect to $\rho$, modulo the equivalence relation of $\rho$-a.e. equality.
	For $f, g \in L^2(\rho)$, their inner product is $\langle f, g \rangle_{\rho} = \int_{\Omega} f(w) \, g(w) \, \rho(\df w)$, and the norm of $f$ is $\|f\|_{\rho} = \sqrt{\langle f, f \rangle_{\rho}}$.
	It will be convenient to consider the subspace $L_0^2(\rho)$, which consists of functions $f \in L^2(\rho)$ such that $\rho f := \int_{\Omega} f(w) \, \rho(\df w) = 0$.
	
	Let $(\Omega', \mathcal{F}', \rho')$ be another probability space.
	Let $J: L_0^2(\rho) \to L_0^2(\rho')$ be a bounded linear transformation.
	The adjoint of $J$ is the linear transformation $J^*: L_0^2(\rho') \to L_0^2(\rho)$ satisfying $\langle Jf,  g \rangle_{\rho'} = \langle f, J^* g \rangle_{\rho}$ for $f \in L_0^2(\rho)$, $g \in L_0^2(\rho')$.
	
	A linear operator $J: L_0^2(\rho) \to L_0^2(\rho)$ is self-adjoint if $J = J^*$.
	It is positive semi-definite if it is self adjoint and $\langle f, Jf \rangle_{\rho} \geq 0$.
	
	Assume that $J: L_0^2(\rho) \to L_0^2(\rho)$ is self-adjoint, and its operator norm $\|J\|_{\rho} \leq 1$.
	For $f \in L_0^2(\rho)$, define the Dirichlet form $\mathcal{E}_J(f) = \|f\|_{\rho}^2 - \langle f, Jf \rangle_{\rho}$.
	Define the right spectral gap of $J$ to be 
	\[
	\mathrm{Gap}(J) = 1 - \sup_{f \in L_0^2(\rho) \setminus \{0\}} \frac{\langle f, Jf \rangle_{\rho}}{\|f\|_{\rho}^2} = \inf_{f \in L_0^2(\rho) \setminus \{0\}} \frac{\mathcal{E}_J(f)}{\|f\|_{\rho}^2}.
	\]
	Note that $\mathrm{Gap}(J) \in [0,2]$, and, if $J$ is positive semi-definite, $\mathrm{Gap}(J) = 1 - \|J\|_{\rho} \leq 1$.
	We also define the left spectral gap to be 
	\[
	\mathrm{Gap}_-(J) = \inf_{f \in L_0^2(\rho) \setminus \{0\}} \frac{\langle f, Jf \rangle_{\rho}}{\|f\|_{\rho}^2} + 1 = 2 - \sup_{f \in L_0^2(\rho) \setminus \{0\}} \frac{\mathcal{E}_J(f)}{\|f\|_{\rho}^2}.
	\]
	Note that $\mathrm{Gap}_-(J) \in [0,2]$, and $\mathrm{Gap}_-(J) \in [1,2]$ if $J$ is positive semi-definite.
	Moreover, $1 - \|J\|_{\rho} = \min\{\mathrm{Gap}(J), \mathrm{Gap}_-(J)\}$.
	{ We keep track of the left spectral gap because a large right spectral gap alone does not rule out (approximate) eigenvalues near $-1$; such eigenvalues may affect the operator norm.
	This quantity will be tracked in Propositions~\ref{pro:left-gap} and~\ref{pro:iteration-1}.
	}

	Finally, we state a useful fact: for a bounded linear transformation $J: L_0^2(\rho) \to L_0^2(\rho')$ whose norm is no greater than~1, it holds that $\mathrm{Gap}(J^*J) = 1 - \|J^*J\|_{\rho} = 1 - \|J J^*\|_{\rho'} = \mathrm{Gap}(J J^*)$.

	\subsection{Linear operator associated with a Markov transition kernel}
	
	A function $K: \Omega \times \mathcal{F}' \to [0,1]$ is an Mtk if $w \mapsto K(w, \mathsf{F})$ is measurable for $\mathsf{F} \in \mathcal{F}'$, and $\mathsf{F} \mapsto K(w, \mathsf{F})$ is a probability measure for $w \in \Omega$.
	For the rest of this subsection, assume that $(\Omega, \mathcal{F}) = (\Omega',\mathcal{F}')$, so $K$ defines the transition law of a Markov chain whose state space is~$\Omega$.
	
	For a probability measure $\mu: \mathcal{F} \to [0,1]$, write $\mu K(\cdot) = \int_{\mathsf{\Omega}} K(w, \cdot) \, \mu(\df w)$.
	For $t \in \mathbb{N}_+ := \{1,2,\dots\}$, $K^t(w, \mathsf{F})$ is $K(w, \mathsf{F})$ if $t = 1$, and $\int_{\Omega} K^{t-1}(w,\df w') \, K(w', \mathsf{F})$ if $t \geq 2$.

	For the rest of this subsection, assume that $\rho K = \rho$.
	Then $K$ defines a linear operator on $L_0^2(\rho)$ via the formula $Kf(w) = \int_{\Omega} K(w, \df w') \, f(w')$.
	It always holds that $\|K\|_{\rho} \leq 1$.
	The Mtk $K$ is reversible with respect to $\rho$ if and only if the operator $K$ is self-adjoint.
	
	The norm $\|K\|_{\rho}$ is closely related to the mixing rate of the corresponding Markov chain.
	It is well-known that, for a large class of initial distributions $\mu$, the $L^2$ distance between $\mu K^t$ and $\rho$ goes to~0 at a rate of $\|K\|^t$ or faster \citep[see, e.g.,][Section 22.2]{douc2018markov}. 
	When $K$ is reversible, this rate is essentially exact \citep{roberts1997geometric}.
	{ We can quantify this notion using the mixing time, which is defined to be
	\[
	t_*(K,\mu,\epsilon) = \min\{t \geq 0: \, \mathrm{TV}(\mu K^t, \rho) \leq \epsilon \}
	\]
	for a given $\epsilon > 0$, where $\mathrm{TV}$ is the total variation distance.
	When $\df \mu/\df \rho$ exists, one can derive \citep[see, e.g.,][]{roberts1997geometric}
	\begin{equation} \label{ine:mixing}
	t_*(K,\mu,\epsilon) \leq \frac{\log \|\df\mu/\df\rho-1\|_{\rho}- \log \epsilon}{-\log \|K\|_{\rho}} \leq \frac{\log \|\df\mu/\df\rho-1\|_{\rho}- \log \epsilon}{1 - \|K\|_{\rho}}.
	\end{equation}
	The denominator $1 - \|K\|_{\rho}$ can be replaced by $\mathrm{Gap}(K)$ if~$K$ is positive semi-definite.
	}

	Assume that $K$ is reversible with respect to~$\rho$.
	For $f \in L_0^2(\rho)$, it holds that
	\[
	\mathcal{E}_K(f) =  \frac{1}{2} \int_{\Omega^2} \rho(\df w) \, K(w, \df w') \, [f(w') - f(w)]^2.
	\]

	Dirichlet forms are connected to the asymptotic variance of ergodic averages.
	For a function $f \in L^2(\rho)$ and a stationary Markov chain $(X(t))_{t=0}^{\infty}$ with Mtk~$K$, the asymptotic variance of $n^{-1/2} \sum_{t=1}^n f(X(t))$ is $\mathrm{var}_K(f) = \|f - \rho f\|_{\rho}^2 + 2 \sum_{t=1}^{\infty} \langle f - \rho f, K^t (f - \rho f) \rangle_{\rho}$, assuming that $\sum_{t=0}^{\infty} \|K^t (f - \rho f) \|_{\rho} < \infty$ \citep[see, e.g.,][Theorem 21.2.6]{douc2018markov}.
	Let $K_1$ and $K_2$ be Mtks that are reversible with respect to $\rho$, and assume that, for $f \in L_0^2(\rho)$, $\mathcal{E}_{K_1}(f) \geq c \, \mathcal{E}_{K_2}(f)$, where $c$ is a positive constant.
	Then, when $\sum_{t=0}^{\infty} \|K^t (f - \rho f) \|_{\rho} < \infty$, for $f \in L^2(\rho)$,
	\begin{equation} \label{ine:asymptotic-variance}
		\mathrm{var}_{K_1}(f) \leq c^{-1} \, \mathrm{var}_{K_2}(f) + (c^{-1} - 1) \|f - \rho f\|_{\rho}^2.
	\end{equation}
	See Lemma 33 in \cite{andrieu18supplement}; see also \cite{caracciolo1990nonlocal}.

	\section{Examples of Spectral Gap Decomposition, Old and New} \label{sec:examples}
	
	In this section, we shall describe five existing and two new examples of spectral gap decomposition.
	Decomposition results involving Dirichlet forms are also provided when it is convenient.
	In Section~\ref{sec:special} and Appendix~\ref{app:derive}, we will show that they can be derived in a unified framework, constructed in Section~\ref{sec:unified}.

	\subsection{Markov chain decomposition based on state space partitioning} \label{ssec:madras}
	
	Markov chain decomposition is a series of spectral gap decomposition results that factor the dynamics of a Markov chain into global and local components, based on a finite partition or covering of the state space.
	
	Let $(\X,\A, \pi)$ be a probability space.
	Suppose that $\X$ has a partition $\X = \bigcup_{z = 1}^k \X_z$, where $k$ is a positive integer, and the $\X_z$'s are non-overlapping measurable subsets such that $\pi(\X_z) > 0$.
	Let $M$ and $N$ be Mtks that are reversible with respect to~$\pi$, and assume that~$M$ is positive semi-definite.

	Let $[k] = \{1,\dots,k\}$, and let $2^{[k]}$ be the power set of $[k]$.
	Let $\varpi: 2^{[k]} \to [0,1]$ be such that $\varpi(\{z\}) = \pi(\X_z)$ for $z \in [k]$.
	Define an Mtk $M_0: [k] \times 2^{[k]} \to [0,1]$ as follows:
	For $z, z' \in [k]$, let 
	\begin{equation} \nonumber
	M_0(z,\{z'\}) = \frac{1}{\pi(\X_z)} \int_{\X_z} \pi(\df x) \, M(x, \X_{z'}).
	\end{equation}
	Then $M_0$ describes the global dynamics across the partition of a stationary chain associated with~$M$.
	It is reversible with respect to $\varpi$.
	
	For $z \in [k]$, let $\A_z$ be $\A$ restricted to $\X_z$, and let $\omega_z: \A_z \to [0,1]$ be such that $\omega_z(\mathsf{A}) = \pi(A)/\pi(\X_z)$.
	Let $H_z: \X_z \times \A_z \to [0,1]$ be such that 
	\begin{equation} \nonumber
	H_z(x, \mathsf{A}) = N(x, \mathsf{A}) + N(x, \X_z^c) \, \delta_x(\mathsf{A}).
	\end{equation}
	Here, $\delta_x$ is the point mass (Dirac measure) at $x$, and $\X_z^c$ is the complement of $\X_z$.
	Then $H_z$ depicts the local dynamics of a chain associated with~$N$ within the subset $\X_z$.
	It can be checked that the Mtk $H_z$ is reversible with respect to $\omega_z$.
	
	The following decomposition result was discovered by \cite{caracciolo1992two} and documented in \cite{madras2002markov}:
	\begin{proposition} \citep{caracciolo1992two} \label{pro:madras}
		In the context of this subsection,
		\[
		\mathrm{Gap}(M^{1/2}N M^{1/2}) \geq \left[ \min_z \mathrm{Gap}(H_z) \right] \mathrm{Gap}(M_0).
		\]
	\end{proposition}
	
	Here, $M^{1/2}$ is the positive square root of the operator~$M$ \citep[see, e.g., ][\S 32, Theorem~4]{helmberg2014introduction}.
	When $M = N$, Proposition \ref{pro:madras} decomposes the right spectral gap of $M^2$ into global and local components.
	
	As we will show later, the sandwich structure $M^{1/2} N M^{1/2}$ alludes to an important feature that leads to spectral gap decomposition in more general settings.
	
	\subsection{Markov chain decomposition, a second form} \label{ssec:madras-1}
	
	\cite{madras2002markov} proposed another form of Markov chain decomposition while assuming the $\X_z$'s may overlap.
	See also \cite{madras1996factoring}.
	
	Let $(\X,\A, \pi_0)$ be a probability space.
	Suppose that $\X = \bigcup_{z=1}^k \X_z$, where $k$ is a positive integer.
	Here, $\X_z \in \A$, $\pi_0(\X_z) > 0$, and the $\X_z$'s may overlap.
	Let $N: \X \times \A \to [0,1]$ be an Mtk that is reversible with respect to $\pi_0$.
	
	Let $\varpi: 2^{[k]} \to [0,1]$ be such that $\varpi(\{z\}) = \pi_0(\X_z)/\Delta$ for $z \in [k]$, where $\Delta = \sum_{z=1}^k \pi_0(\X_z)$.
	Let $\Theta = \max_x \sum_{z=1}^k \ind_{\X_z}(x)$.
	For $z, z' \in [k]$, let
	\[
	\Pi_0(z, \{z'\}) = \frac{\pi_0(\X_z \cap \X_{z'})}{\Theta \, \pi_0(\X_z)}
	\]
	if $z \neq z'$, and let $\Pi_0(z, \{z\}) = 1 - \sum_{z'' \neq z} \Pi_0(z, \{z''\})$.
	Then $\Pi_0$ is an Mtk that is reversible with respect to $\varpi$.
	
	For $z \in [k]$, let $\A_z$ be $\A$ restricted to $\X_z$, and let $\omega_z: \A_z \to [0,1]$ be such that $\omega_z(\mathsf{A}) = \pi_0(A)/\pi_0(\X_z)$. 
	Let $H_z: \X_z \times \A_z \to [0,1]$ be such that 
	\begin{equation} \nonumber
		H_z(x, \mathsf{A}) = N(x, \mathsf{A}) + N(x, \X_z^c) \, \delta_x(\mathsf{A}).
	\end{equation}
	Then $H_z$ is reversible with respect to $\omega_z$.
	
	The following result was established in \cite{madras2002markov}.
	
	\begin{proposition} \label{pro:madras-1} \citep{madras2002markov}
		In the context of this subsection, 
		\[
		\mathrm{Gap}(N) \geq \Theta^{-2} \left[ \min_z \mathrm{Gap}(H_z) \right] \, \mathrm{Gap}(\Pi_0).
		\]
	\end{proposition}
	
	{ The two factors of $\Theta^{-1}$ in Proposition~\ref{pro:madras-1} have different origins.  In the derivation in Appendix~\ref{ssec:decomposition-1}, the first comes from the approximate $z$-invariance condition \ref{H4}, for which $c_0=\Theta^{-1}$; the second comes from the Dirichlet-form and stationary-measure comparison between the original kernel~$N$ and the reweighted kernel~$S$ in \eqref{ine:overlap-3}.}
	
	
	\subsection{Hybrid data augmentation algorithms with one intractable conditional} \label{ssec:andrieu}
	
	Spectral gap decomposition can appear in situations that seem quite distinct from the ones in Sections \ref{ssec:madras} and \ref{ssec:madras-1}.
	The next few examples involve Gibbs-like algorithms and their hybrid variants.
	
	Let $(\X,\A, \pi)$ and $(\Z,\C, \varpi)$ be probability spaces.
	Let $\tilde{\pi}: \A \times \C \to [0,1]$ be a probability measure such that $\tilde{\pi}(\mathsf{A} \times \Z) = \pi(\mathsf{A})$ for $\mathsf{A} \in \A$ and $\tilde{\pi}(\X \times \mathsf{C}) = \varpi(\mathsf{C})$ for $\mathsf{C} \in \C$.
	Assume that $\tilde{\pi}$ has the decomposition
	\[
	\tilde{\pi}(\df (x,z)) = \pi(\df x) \, \pi_x(\df z) = \varpi(\df z) \, \varpi_z(\df x).
	\]
	This means that, if $(X,Z) \sim \tilde{\pi}$, then $X \sim \pi$, $Z \sim \varpi$, $Z \mid X = x \sim \pi_x$, and $X \mid Z = z \sim \varpi_z$.
	
	An ideal data augmentation algorithm \citep{tanner1987calculation} targeting~$\pi$ is defined by the Mtk
	\[
	\bar{S}(x, \df x') = \int_{\Z} \pi_x(\df z) \, \varpi_z(\df x').
	\]
	Given the current state $x$, one obtains the next state~$x'$ by sampling~$z$ from $\pi_x$, and drawing~$x'$ from $\varpi_z$.
	This Mtk is reversible with respect to $\pi$, and is positive semi-definite.
	
	In practice, $\varpi_z$ may be difficult to sample from directly, and one may consider replacing an exact draw from $\varpi_z$ with a Markovian step, e.g., a Metropolis-Hastings step.
	Suppose that, for each $z \in \Z$, there is an Mtk $Q_z: \X \times \A \to [0,1]$ that is reversible with respect to $\varpi_z$.
	Assume that $(x,z) \mapsto Q_z(x, \mathsf{A})$ is measurable for $\mathsf{A} \in \A$.
	Then a hybrid data augmentation chain can be defined by the Mtk
	\[
	S(x, \df x') = \int_{\Z} \pi_x(\df z) \,  Q_z(x, \df x').
	\]
	It can be shown that~$S$ is reversible with respect to~$\pi$.
	To simulate the chain, given the current state $x$, one draws $z$ from $\pi_x$, then draws the next state $x'$ from the distribution $Q_z(x, \cdot)$.
	
	The following decomposition result was established in \cite{andrieu2018uniform}.
	
	\begin{proposition} \citep{andrieu2018uniform} \label{pro:andrieu}
		In the context of this subsection, for $f \in L_0^2(\pi)$,
		\[
		\mathcal{E}_S(f) \geq \left[ \essinf_z \mathrm{Gap}(Q_z) \right] \, \mathcal{E}_{\bar{S}}(f), \quad
		\mathrm{Gap}(S) \geq \left[ \essinf_z \mathrm{Gap}(Q_z) \right] \, \mathrm{Gap}(\bar{S}),
		\]
		where the essential infimum is defined with respect to the measure~$\varpi$.
	\end{proposition}

	\subsection{Random-scan hybrid Gibbs samplers} \label{ssec:qin}
	
	We now investigate random-scan hybrid Gibbs samplers, which are commonly used for sampling from multi-dimensional distributions.
	
	Consider the case where $(\X, \A)$ is a product of measurable spaces.
	To be specific, let $\X = \X_1 \times \cdots \times \X_k$, $\A = \A_1 \times \cdots \times \A_k$, where $k$ is a positive integer, and each $\X_i$ is a Polish space equipped with Borel algebra $\A_i$.
	For a vector $x = (x_1, \dots, x_k) \in \X$, let $x_{-i} = (x_1, \dots, x_{i-1}, x_{i+1}, \dots, x_k)$.
	Let $\X_{-i} = \{x_{-i}: \, (x_1, \dots, x_k) \in \X\}$.
	Let $\pi: \A \to [0,1]$ be a probability distribution.
	Let $(X_1, \dots, X_k)$ be a random vector distributed as~$\pi$.
	For $i \in [k]$ and $u \in \X_{-i}$, let $\varphi_{i,u}: \A_i \to [0,1]$ be the conditional distribution of $X_i \mid X_{-i} = u$.
	Let $p_1, \dots, p_k$ be positive constants satisfying $\sum_{i=1}^k p_i = 1$.
	
	An ideal random-scan Gibbs sampler is defined by the Mtk
	\[
	\bar{S}(x, \df x') = \sum_{i=1}^{k} p_i \, \varphi_{i, x_{-i}}(\df x'_i) \, \delta_{x_{-i}}(\df x'_{-i}).
	\]
	Given the current state $x$, to sample the next state~$x'$, pick a random index~$i$ according to the probability vector $(p_1, \dots, p_k)$, draw $x_i'$ from $\varphi_{i,x_{-i}}$, and set $x_{-i}' = x_{-i}$.
	The Mtk is reversible with respect to $\pi$.
	
	In practice, it may be difficult to make exact draws from $\varphi_{i,u}$, and one may replace the exact draw with a Markovian step.
	Let $H_{i,u}: \X_i \times \A_i \to [0,1]$ be an Mtk that is reversible with respect to $\varphi_{i,u}$, and assume that $x \mapsto H_{i,x_{-i}}(x_i, \mathsf{A}_i)$ is measurable for $\mathsf{A}_i \in \A_i$.
	A random-scan hybrid Gibbs sampler is defined by the Mtk
	\[
	S(x, \df x') = \sum_{i=1}^{k} p_i \, H_{i, x_{-i}}(x_i, \df x'_i) \, \delta_{x_{-i}}(\df x'_{-i}).
	\]
	Then~$S$ is reversible with respect to $\pi$.
	Given the current state $x$, to sample the next state~$x'$, one picks a random index~$i$ according to the probability vector $(p_1, \dots, p_k)$, draws $x_i'$ from $H_{i,x_{-i}}(x_i, \cdot)$, and sets $x_{-i}' = x_{-i}$.
	
	The following result was proved in \cite{qin2024spectral}.
	
	\begin{proposition} \citep{qin2024spectral} \label{pro:qin}
		In the context of this subsection, for $f \in L_0^2(\pi)$,
		\[
		\mathcal{E}_S(f) \geq \left[ \min_i \essinf_u \mathrm{Gap}(H_{i,u}) \right] \, \mathcal{E}_{\bar{S}}(f), \quad \mathrm{Gap}(S) \geq \left[ \min_i \essinf_u \mathrm{Gap}(H_{i,u}) \right] \, \mathrm{Gap}(\bar{S}),
		\]
		where the essential infimum is taken with respect to the marginal distribution of $X_{-i}$ when $X \sim \pi$.
	\end{proposition}
	We shall make the observation that a random-scan (hybrid) Gibbs sampler may be understood as a (hybrid) data augmentation algorithm \citep[see, e.g.,][]{andersen2007hit}.
	Thus, Propositions \ref{pro:andrieu} and \ref{pro:qin} can be unified.
	
	\subsection{Approximate variance conservation in localization schemes} \label{ssec:chen}
	
	Localization schemes are a framework proposed by \cite{chen2022localization} that can be used for analyzing Markov chains.
	This framework builds on the idea of stochastic localization \citep{eldan2013thin}, and extends techniques based on spectral independence \citep{anari2021spectral}, which have been widely applied to analyze Markov chains in theoretical computer science \citep{chen2021rapid,feng2022rapid,chen2023rapid,qin2024spectrala,chen2024rapid}.
	
	We shall review a key result in \cite{chen2022localization} involving variance conservation, and explain how it can be unified with the other decomposition results in this section.
	We formalize localization schemes in a way slightly different from \cite{chen2022localization}, but the fundamental idea remains the same.
	
	Let $(\X, \A, \pi)$ be a probability space.
	Localization schemes consider Mtks defined by a random sequence of probability measures $(\nu_t)_{t=0}^{\infty}$ on $(\X,\A)$ initialized at~$\pi$.
	
	For $s \in \mathbb{N} := \{0,1,2,\dots\}$, let $\mathsf{W}_s$ be a Polish space and let $\mathcal{D}_s$ be its Borel algebra.
	Let $(W_s)_{s=0}^{\infty}$ be a sequence of random elements, with $W_s$ taking values in $\mathsf{W}_s$.
	For $t \in \mathbb{N}$, let $\Z_t = \mathsf{W}_0 \times \mathsf{W}_1 \times \cdots \times \mathsf{W}_t$, and denote by $\mu_t$ the distribution of $(W_s)_{s=0}^t$.

	For $t \in \mathbb{N}$, let $v_t: \Z_t \times \X \to [0,\infty]$ be a measurable function such that the following conditions hold:
	\begin{itemize}
		\item  [\namedlabel{A1}{(A1)}] For $t \in \mathbb{N}$, almost surely,
		$
		\mathsf{A} \mapsto \int_{\mathsf{A}} v_t((W_s)_{s=0}^t, x) \, \pi(\df x)
		$
		is a probability measure on $\A$.

		\item [\namedlabel{A2}{(A2)}]
		For $x \in \X$, the sequence $v_t((W_s)_{s=0}^t, x)$, $t \in \mathbb{N}$, forms a martingale initialized at~1: $v_t(W_0, x) = 1$ almost surely, and,
		for $t \in \mathbb{N}$,
		\[
		 \mathrm{E}[v_{t+1}((W_s)_{s=0}^{t+1}, x) \mid (W_s)_{s=0}^t] = v_t((W_s)_{s=0}^t,x) \quad \text{almost surely.}
		\]
		
	\end{itemize}
	
	For $t \in \mathbb{N}$ and $\mathsf{A} \in \A$, let 
	\[
	\nu_t(\mathsf{A}) = \int_{\mathsf{A}} v_t((W_s)_{s=0}^t, x) \, \pi(\df x).
	\]
	Then $(\nu_t)_{t=0}^{\infty}$ can be seen as a sequence of random probability measures initialized at $\pi$.
	
	For $t \in \mathbb{Z}_+$, define the kernel
	\[
	K_t(x, \mathsf{A}) = \mathrm{E}\left[ \frac{\df \nu_t}{\df \pi}(x) \, \nu_t( \mathsf{A}) \right] =  \int_{\Z_t} v_t(z, x) \left[ \int_{\mathsf{A}} v_t(z, x') \, \pi(\df x') \right]  \mu_t(\df z).
	\]
	By \ref{A1} and \ref{A2}, this is an Mtk reversible with respect to~$\pi$, and it is positive semi-definite.
	
	In \cite{chen2022localization}, it is shown that some important Markov chains in theoretical computer science, e.g., Glauber dynamics, can be formulated in the above manner.
	We refer readers to Section~2 of that work for a variety of examples.
	
	
	\cite{chen2022localization} derived the following simple proposition relating $\mathrm{Gap}(K_t)$ to $(\nu_s)_{s=0}^t$.
	
	\begin{proposition} \citep{chen2022localization} \label{pro:chen}
		Suppose that \ref{A1} and \ref{A2} hold. Then each of the following holds.
		\begin{enumerate}[label=(\alph*)]
			\item For $t \in \mathbb{N}$, 
			\begin{equation} \label{eq:GapKt}
				\mathrm{Gap}(K_t) = \inf_{f \in L_0^2(\pi) \setminus \{0\}} \frac{\mathrm{E} [ \mathrm{var}_{\nu_t}(f) ] }{\mathrm{var}_{\pi}(f)},
			\end{equation}
			where, for a probability measure $\nu: \A \to [0,1]$, $\mathrm{var}_{\nu}(f)$ means $\int_{\X} [f(x) - \nu f]^2 \, \nu(\df x)$.
			
			\item Assume further that, for some $t \in \mathbb{Z}_+$, there exist positive constants $\kappa_1, \dots, \kappa_t$ such that, for $s \in \{0,\dots,t-1\}$, the following ``approximate conservation of variance" holds almost surely:
			For $f \in L_0^2(\nu_s)$,
			\begin{equation} \label{ine:specind}
				\mathrm{E}[\mathrm{var}_{\nu_{s+1}}(f) \mid (W_i)_{i=0}^s] \geq \kappa_{s+1} \mathrm{var}_{\nu_s}(f).
			\end{equation}
			Then 
			\begin{equation} \label{ine:chen-decomposition}
				\mathrm{Gap}(K_t) \geq \prod_{s=1}^t \kappa_s
			\end{equation}
		\end{enumerate}

	\end{proposition}

	Proposition~\ref{pro:chen} allows one to bound the spectral gap of an Mtk $K_t$ by studying the behavior of $\nu_s$ as $s$ goes from 0 to $t$.
	In \cite{chen2022localization}, this result yielded a rather simple proof of the main theorem in \cite{anari2021spectral}, which was in turn a breakthrough in the analysis of Glauber dynamics in spin systems.
	In particular, it was shown that the approximate conservation of variance encompasses the spectral independence condition from \cite{anari2021spectral}.
 	
 	In Appendix~\ref{app:derive}, we demonstrate how Proposition \ref{pro:chen} can be unified with other results in this section.
 	It is shown that, for $s \in \mathbb{N}$, $K_s$ can be viewed as the Mtk of a data augmentation algorithm, i.e., $\bar{S}$ in Section~\ref{ssec:andrieu}, while $K_{s+1}$ can be seen as the Mtk of a hybrid data augmentation algorithm with one conditional distribution replaced by a Markovian approximation, i.e., $S$ in Section~\ref{ssec:andrieu}.
 	It is then demonstrated that \eqref{ine:specind} is equivalent to $\mathrm{Gap}(Q_{s,z}) \geq \kappa_{s+1}$, where $Q_{s,z}$ corresponds to the Markovian approximation in $S = K_{s+1}$.
 	Thus, \eqref{ine:chen-decomposition} is a consequence of Proposition~\ref{pro:andrieu}, which states that $\mathrm{Gap}(K_{s+1}) \geq [ \essinf_z \mathrm{Gap}(Q_{s,z}) ] \, \mathrm{Gap}(K_s)$.

	\subsection{Hybrid hit-and-run samplers} \label{ssec:hit-0}
	
	In this and the next subsection, we give two new examples of spectral gap decomposition.
	
	It is well-known that hit-and-run samplers are similar to data augmentation and Gibbs samplers \citep{andersen2007hit}.
	Thus, in light of the results in Sections \ref{ssec:andrieu} and \ref{ssec:qin}, we also expect spectral gap decomposition to hold for hybrid hit-and-run samplers.
	Consider a multi-dimensional hit-and-run sampler described in \cite{ascolani2024entropy}.
	
	Let $\X = \mathbb{R}^k$, where $k$ is a positive integer, and let $\A$ be the Borel algebra.
	Assume that $\pi: \A \to [0,1]$ admits a density function $\dot{\pi}: \X \to [0,\infty]$ with respect to the Lebesgue measure on~$\X$.
	
	Let $\ell$ be a positive integer no greater than~$k$.
	Let $\mathsf{W}$ be the following Stiefel manifold:
	\[
	\mathsf{W} = \left\{(w_1, \dots, w_{\ell}) \in (\mathbb{R}^k)^{\ell}: \, w_i^{\top} w_i = 1 \text{ for each } i, \; w_i^{\top} w_j = 0 \text{ if } i \neq j \right\}.
	\]
	Then each element of $\mathsf{W}$ is an ordered orthonormal basis for some $\ell$-dimensional subspace of $\mathbb{R}^k$.
	Let $\nu: \mathcal{D} \to [0,\infty)$ be a probability measure on $(\mathsf{W}, \mathcal{D})$, where $\mathcal{D}$ is the Borel algebra of the Stiefel manifold.
	For $x \in \X$ and $w = (w_1, \dots, w_{\ell}) \in \mathsf{W}$, let $\theta(x,w) = x - \sum_{i=1}^{\ell} (w_i^{\top} x) w_i$.
	Then $w$ and $\theta(x,w)$ identify a hyperplane $\mathsf{L}_{w, \theta(x,w)}$ passing through~$x$ that is parallel to the span of $\{w_1, \dots, w_{\ell}\}$, where, for $x' \in \mathbb{R}^k$,
	\[
	\mathsf{L}_{w, x'} = \left\{ x' + \sum_{i=1}^{\ell} u_i w_i: \, u  \in \mathbb{R}^{\ell} \right\}.
	\]
	Define the following probability measure on $\mathbb{R}^{\ell}$:
	\[
	\varphi_{w, \theta(x,w)}(\df u) = \frac{\dot{\pi}( \theta(x,w) + \sum_i^{\ell} u_i w_i ) \, \df u}{ \int_{\mathbb{R}^{\ell}} \dot{\pi}( \theta(x,w) + \sum_i^{\ell} u'_i w_i ) \, \df u'},
	\]
	provided that the denominator is in $(0,\infty)$, which is the case for each $w \in \mathsf{W}$ and $\pi$-a.e. $x \in \X$.
	Let
	\begin{equation} \nonumber
	\varpi_{w, \theta(x,w)}(\mathsf{A}) = \int_{\mathbb{R}^{\ell}} \varphi_{w, \theta(x,w)}(\df u) \, \ind_{\mathsf{A}}\left(\theta(x,w) + \sum_{i=1}^{\ell} u_i w_i \right), \quad \mathsf{A} \in \A.
	\end{equation}
	That is, $\varpi_{w, \theta(x,w)}$ is the distribution of $\theta(x,w) + \sum_{i=1}^{\ell} u_i w_i$ if $u \sim \varphi_{w, \theta(x,w)}$.
	Equivalently, $\varpi_{w, \theta(x,w)}$ is the distribution of $x + \sum_{i=1}^{\ell} u_i w_i$ if $u$ is distributed according to the density $\dot{\pi}(x + \sum_{i=1}^{\ell} u_i w_i )$.
	This distribution can be understood as $\pi$ restricted to the hyperplane $\mathsf{L}_{w, \theta(x,w)}$.

	An ideal $\ell$-dimensional hit-and-run algorithm targeting~$\pi$ simulates a Markov chain through the following procedure:
	Given the current state $x \in \X$, an orthonormal basis $(w_1, \dots, w_{\ell})$ is drawn from the distribution~$\nu$; 
	then the new state is drawn from $\varpi_{w, \theta(x,w)}$.
	The corresponding Mtk is given by the formula
	\[
	\begin{aligned}
		\bar{S}(x, \mathsf{A}) &= \int_{\mathsf{W}} \nu(\df w) \int_{\mathbb{R}^{\ell}} \varphi_{w, \theta(x,w)}(\df u) \, \ind_{\mathsf{A}}\left(\theta(x,w) + \sum_{i=1}^{\ell} u_i w_i \right).
	\end{aligned}
	\]
	
	
	It may be difficult to make exact draws from $\varpi_{w, \theta(x,w)}$ when $\varphi_{w, \theta(x,w)}$ has an intricate form.
	Let $H_{w, \theta(x,w)}$ be an Mtk that is reversible with respect to $\varphi_{w, \theta(x,w)}$.
	We may then define a hybrid hit-and-run algorithm with Mtk
	\[
	S(x, \mathsf{A}) = \int_{\mathsf{W}} \nu(\df w) \int_{\mathbb{R}^{\ell}} H_{w,\theta(x,w)}\left((w_i^{\top} x)_{i=1}^{\ell}, \df u \right) \, \ind_{\mathsf{A}}\left(\theta(x,w) + \sum_{i=1}^{\ell} u_i w_i \right) .
	\]
	To simulate this chain, one draws~$w$ from~$\nu$, then changes the current state from $x = \theta(x,w) + \sum_{i=1}^{\ell} (w_i^{\top} x) w_i$ to $\theta(x,w)  + \sum_{i=1}^{\ell} u_i w_i$, where $u \in \mathbb{R}^{\ell}$ is drawn from  the distribution $H_{w, \theta(x,w)}((w_i^{\top} x)_{i=1}^{\ell}, \cdot)$.
	To ensure $S$ is well-defined, it is assumed that $(x,w,u) \mapsto H_{w,\theta(x,w)}(u, \mathsf{B})$ is measurable for any measurable $\mathsf{B} \subset \mathbb{R}^{\ell}$.
	It will be shown that~$S$ is reversible with respect to~$\pi$.
	See Section \ref{sec:example-hybrid} for a concrete example of~$S$ when $H_{w,\theta(x,w)}$ corresponds to a random-walk Metropolis algorithm.
	
	While many works have focused on the convergence properties of the ideal hit-and-run sampler \citep[see, e.g.,][]{lovasz1999hit,lovasz2004hit,chewi2022query,ascolani2024entropy,bou2025accelerated}, there is a lack of study on its hybrid counterpart.
	In Appendix \ref{ssec:hit-1}, we establish the following decomposition result for hybrid hit-and-run samplers:
	\begin{proposition} \label{pro:hit}
		In the context of this subsection, $S$ and $\bar{S}$ are both reversible with respect to $\pi$, and, for $f \in L_0^2(\pi)$,
		\[
		\mathcal{E}_S(f) \geq \left[ \essinf_{w,x} \mathrm{Gap}(H_{w, \theta(x,w)}) \right] \mathcal{E}_{\bar{S}}(f), \quad \mathrm{Gap}(S) \geq \left[ \essinf_{w,x} \mathrm{Gap}(H_{w, \theta(x,w)}) \right] \mathrm{Gap}(\bar{S}),
		\]
		where the essential infimum is taken with respect to $\nu \times \pi$.
	\end{proposition}
	

	\subsection{Hybrid data augmentation algorithms with two intractable conditionals} \label{ssec:doubly-0}
	
	In this subsection, we study a scenario involving a non-reversible chain.
	
	Recall the setting in Section \ref{ssec:andrieu}, but with a different set of notations.
	Let $(\X_1, \A_1, \varphi_1)$ and $(\X_2, \A_2, \varphi_2)$ be probability spaces.
	Suppose that $\varphi: \A_1 \times \A_2 \to [0,1]$ is a probability measure with the disintegration
	\[
	\varphi(\df (x_1, x_2)) = \varphi_1(\df x_1) \, \varphi_{2, x_1}(\df x_2) = \varphi_2(\df x_2) \, \varphi_{1, x_2}(\df x_1).
	\]
	That is, if $(X_1, X_2) \sim \varphi$, then $\varphi_i$ is the marginal distribution of $X_i$, and $\varphi_{i, x_{-i}}$ is the conditional distribution of $X_i$ given $X_{-i} = x_{-i}$.
	The Mtk of a data augmentation algorithm targeting $\varphi_2$ is of the form
	\[
	\bar{S}(x_2, \df x_2') = \int_{\X_1} \varphi_{1,x_2}(\df x_1) \, \varphi_{2,x_1}(\df x_2').
	\]
	
	In Section \ref{ssec:andrieu}, it was assumed that $\varphi_{2,x_1}$ is intractable.
	Consider now the doubly intractable case, where $\varphi_{1,x_2}$ is also intractable.
	Suppose that one can simulate $H_{2,x_1}$ and $H_{1,x_2}$, where, for $i \in \{1,2\}$, $H_{i,x_{-i}}: \X_i \times \A_i \to [0,1]$ is an Mtk that is reversible with respect to $\varphi_{i,x_{-i}}$.
	It is assumed that $(x_1, x_2) \mapsto H_{i,x_{-i}}(x_i, \mathsf{A})$ is measurable if $\mathsf{A} \in \A_i$.
	
	We are interested in the Mtk
	\[
	T((x_1, x_2), \df (x_1', x_2')) = H_{1,x_2}(x_1, \df x_1') \, H_{2,x_1'}(x_2, \df x_2').
	\]
	One can check that $\varphi T = \varphi$, but $T$ is not necessarily reversible.

	The norm of $T$ will be compared to the right spectral gaps of $\bar{S}$ and the following Mtks:
	\[
	\hat{S}_1(x_2, \df x_2') = \int_{\X_1} \varphi_{1,x_2}(\df x_1) \, H_{2,x_1}(x_2, \df x_2'), \quad \hat{S}_2(x_2, \df x_2') = \int_{\X_1} \varphi_{1,x_2}(\df x_1) \, H_{2,x_1}^2(x_2, \df x_2').
	\]
	Note that $\hat{S}_1$ and $\hat{S}_2$ can be viewed as hybrid data augmentation algorithms, each with one intractable conditional, and are subject to Proposition \ref{pro:andrieu}.
	Thus,
	\[
	\mathrm{Gap}(\hat{S}_1) \geq \left[\essinf_{x_1} \mathrm{Gap}(H_{2,x_1}) \right] \, \mathrm{Gap}(\bar{S}), \quad \mathrm{Gap}(\hat{S}_2) \geq \left( 1 - \esssup_{x_1} \| H_{2, x_1} \|_{\varphi_{2,x_1}}^2 \right) \mathrm{Gap}(\bar{S}).
	\]
	Here, we have used the fact that $H_{2,x_1}$ is reversible, which implies that $\mathrm{Gap}(H_{2,x_1}^2) = 1 - \|H_{2,x_1}\|_{\varphi_{2,x_1}}^2$.
	
	In Appendix~\ref{ssec:doubly-1}, we establish the following proposition relating $T$ to the other Mtks.
	
	\begin{proposition} \label{pro:doubly}
		In the context of this subsection,
		\[
		\begin{aligned}
			1 - \|T\|_{\varphi}^2 &\geq \left( 1 - \esssup_{x_2} \| H_{1, x_2} \|_{\varphi_{1, x_2}}^2 \right) \mathrm{Gap}(\hat{S}_2) ,
		\end{aligned}
		\]
		where the essential supremum is taken with respect to $\varphi_2$.
		If, furthermore, $H_{2,x_1}$ is positive semi-definite for $\varphi_1$-a.e. $x_1 \in \X_1$, then $\mathrm{Gap}(\hat{S}_2) \geq \mathrm{Gap}(\hat{S}_1)$.
	\end{proposition}
	
	{  Subsequent to the initial release of this article, \cite{gao2026weak} followed up on the above result and Proposition~\ref{pro:weak} below by developing complementary weak Poincar\'e comparison theory for hybrid data augmentation samplers.
	Their analysis yields explicit subgeometric convergence bounds for $T$ when $H_{i,x_{-i}}$ and $\bar{S}$ satisfy suitable weak Poincar\'{e} inequalities.}

%

	\section{A Unified Framework} \label{sec:unified}
	
	\subsection{A sandwich structure} \label{ssec:sandwich}
	
	We now describe a general framework in which spectral gap decomposition can be derived.
	We shall begin by introducing some notations.
	See Table \ref{tab:notation} for a summary.

	\begin{table} [t]
		\caption{Notations used in Section~\ref{ssec:sandwich}} \label{tab:notation}
		\begin{tabular}{>{ }l >{ }l}
			\hline 
			Notation & Meaning  \\
			\hline
			$\pi$ & a distribution on $\X$ \\  $\tilde{\pi}$ & a distribution on $\Y \times \Z$ \\
			$\varpi$ & marginal distribution of $\tilde{\pi}$ on $\Z$ \\ $\varpi_z$ & a conditional distribution on $\Y$ defined via $\tilde{\pi}(\df (y,z)) = \varpi(\df z) \, \varpi_z(\df y)$ \\
			$P$ & a contraction mapping a function on $\X$ to one on $\Y \times \Z$ \\
			$P^*$ & the adjoint of $P$, mapping a function on $\Y \times \Z$ to one on $\X$ \\
			$Q$ & an Mtk on $\Y \times \Z$ that is reversible with respect to $\tilde{\pi}$ \\
			$S$ & $P^*QP$ \\
			$Q_z$ & an Mtk on $\Y$ indexed by $z \in \Z$, related to~$Q$ through \ref{H4} \\
			$\hat{Q}$ & an Mtk on $\Y \times \Z$ given by $\hat{Q}((y,z), \df (y',z')) = Q_z(y, \df y') \, \delta_z(\df z')$ \\
			$c_0$ & a number in $(0,1]$ used in \ref{H4} \\
			$\Proj$ & a functional mapping $h: \Y \times \Z \to \mathbb{R}$ to $z \mapsto \int_{\Y} h(y,z) \, \varpi_z(\df y)$ \\
			$\Proj^*$ & a functional mapping $g: \Z \to \mathbb{R}$ to $(y,z) \mapsto g(z)$ \\
			$\kappa^{\dagger}$ & a uniform lower bound on $\mathrm{Gap}(Q_z)$, used in Theorem \ref{thm:uniform} \\
			$\kappa^{\ddagger}$ & a uniform lower bound on $\mathrm{Gap}_-(Q_z)$, used in Proposition \ref{pro:left-gap} \\
			\hline
		\end{tabular}
	\end{table}

	Let $(\X,\A)$, $(\Y,\B)$, and $(\Z,\C)$ be measurable spaces.
	Let $\pi: \A \to [0,1]$ and $\tilde{\pi}: \B \times \C \to [0,1]$ be probability measures.
	Let $\varpi(\mathsf{C}) = \tilde{\pi}(\Y \times \mathsf{C})$ for $\mathsf{C} \in \C$, and assume that $\tilde{\pi}(\df (y, z)) = \varpi_z(\df y) \, \varpi(\df z)$.
	In other words, if $(Y,Z) \sim \tilde{\pi}$, then $Z \sim \varpi$, and $Y \mid Z = z \sim \varpi_z$.

	Let $S: \X \times \A \to [0,1]$ be an Mtk that is reversible with respect to~$\pi$, so that it can be regarded as a self adjoint linear operator on $L_0^2(\pi)$.
	We say $S$ has a sandwich structure with an approximate $z$-invariant core if the following holds:
	
	\begin{itemize}
		\item [\namedlabel{H1}{(H1)}] The operator $S$ can be written into the form $S = P^*QP$, where $P: L_0^2(\pi) \to L_0^2(\tilde{\pi})$ is a linear transformation, and $Q$ corresponds to an Mtk on $(\Y \times \Z) \times (\B \times \C)$ that is reversible with respect to~$\tilde{\pi}$.
		
		\item [\namedlabel{H2}{(H2)}] The linear transformation $P$ satisfies $\|P^*P\|_{\pi} \leq 1$.
		Equivalently, $\|Pf\|_{\tilde{\pi}} \leq \|f\|_{\pi}$ for $f \in L_0^2(\pi)$.
		
		\item [\namedlabel{H3}{(H3)}] There exists a set of Mtks $\{Q_z\}_{z \in \Z}$ on $\Y \times \B$ that satisfy the following:
		\begin{enumerate}
			\item[(i)] For $\varpi$-a.e. $z \in \Z$, $Q_z: \Y \times \B \to [0,1]$ is reversible with respect to $\varpi_z$.
			\item[(ii)] For $\mathsf{B} \in \B$, the function $(y,z) \mapsto Q_z(y, \mathsf{B})$ is measurable.
		\end{enumerate}
		
		\item [\namedlabel{H4}{(H4)}]
		The Mtk $Q$ is approximately $z$-invariant in the following sense:
		For $\tilde{\pi}$-a.e. $(y,z) \in \Y \times \Z$ and $\mathsf{B} \in \B$, 
		\[
		Q((y,z), (\mathsf{B} \setminus \{y\}) \times \{z\}) \geq c_0 \, Q_z(y, \mathsf{B} \setminus \{y\}),
		\]
		where $c_0 \in (0,1]$ is a constant.
	\end{itemize}
	
	{ In most applications herein, $c_0 = 1$ in \ref{H4}. The construction in Appendix~\ref{ssec:decomposition-1}, which concerns the example in Section~\ref{ssec:madras-1}, is the only application herein for which $c_0<1$.  The abstraction is retained because it also allows future applications based on approximate Peskun comparisons, in which the available within-$z$ flow dominates that of $Q_z$ only up to a constant.}
	
	\begin{remark}
		If  $P$ is associated with an Mtk $P: (\Y \times \Z) \times \A \to [0,1]$ such that $$\int_{\Y \times \Z} \tilde{\pi}(\df (y,z)) \, P((y,z), \mathsf{A}) = \pi(\mathsf{A})$$ for $\mathsf{A} \in \A$, then \ref{H2} holds by the Cauchy-Schwarz inequality.
	\end{remark}
	

	Under \ref{H3}, we may define an Mtk $\hat{Q}: (\Y \times \Z) \times (\B \times \C) \to [0,1]$ as follows:
	\[
	\hat{Q}((y,z), \df (y', z')) = Q_z(y, \df y') \, \delta_z(\df z').
	\] 
	Note that $\hat{Q}$ is reversible with respect to $\tilde{\pi}$.
	\ref{H4} can be viewed as some form of Peskun ordering \citep{peskun1973optimum,tierney1998note}.
	It leads to the following comparison result regarding the self-adjoint operators $P^*QP$ and $P^*\hat{Q}P$, which is proved in Appendix \ref{app:peskun}.
	
	\begin{lemma} \label{lem:peskun}
		Assume that \ref{H1} to \ref{H4} hold.
		Then, for $f \in L_0^2(\pi)$, $\mathcal{E}_S(f) \geq c_0 \, \mathcal{E}_{P^*\hat{Q}P}(f)$.
		In particular, $\mathrm{Gap}(S) \geq c_0 \, \mathrm{Gap}(P^* \hat{Q} P)$.
	\end{lemma}
	
	\begin{remark}
		Note that $\hat{Q}$ itself is $z$-invariant with $c_0 = 1$.
		We will see, in many applications, one can actually take $Q = \hat{Q}$, so $S = P^* \hat{Q} P$.
	\end{remark}
	
	{ 
	\begin{remark} \label{rem:H4}
		In all our results regarding $S$, conditions \ref{H1} and \ref{H4} are solely used to establish the conclusion of Lemma~\ref{lem:peskun}:
		\begin{itemize}
			\item [\namedlabel{H5}{(H5)}] For $f \in L_0^2(\pi)$, $\mathcal{E}_S(f) \geq c_0 \mathcal{E}_{P^*\hat{Q}P}(f)$.
		\end{itemize}
		Therefore, in Theorem~\ref{thm:uniform}, Corollary~\ref{cor:uniform}, and Lemma~\ref{lem:weak-new}, stated later in this article, one may replace \ref{H1} and \ref{H4} with the weaker (but less explicit) condition \ref{H5}.
		In other words,~$S$ does not have to possess a sandwich structure, but one must be able to compare it to a kernel $P^*\hat{Q}P$ that has the sandwich structure.
		See Section~\ref{app:beyond} for a more detailed discussion.
	\end{remark}
	}
	
	We will examine the spectral gap of $P^* \hat{Q} P$. 
	This operator has the integral kernel form
	\[
	P^* \hat{Q} P(x, \df x') = \int_{\Y \times \Z \times \Y} P^*(x, \df (y,z)) \, Q_z(y, \df y') \, P((y',z), \df x'),
	\]
	assuming that $P$ and $P^*$ themselves can be written into kernel forms.

	For $h \in L_0^2(\tilde{\pi})$, let $\Proj h \in L_0^2(\varpi)$ be such that
	\begin{equation} \label{eq:Proj}
	\Proj h(z) = \int_{\Y} h(y,z) \, \varpi_z(\df y), \quad z \in \Z.
	\end{equation}
	Then $\Proj: L_0^2(\tilde{\pi}) \to L_0^2(\varpi)$ is a bounded linear transformation.
	Note also that, for $g \in L_0^2(\varpi)$, $\Proj^*g(y,z) = g(z)$ for $(y,z) \in \Y \times \Z$.
	Define the positive semi-definite operator $P^* \Proj^* \Proj P$.
	This operator has kernel form
	\[
	P^* \Proj^* \Proj P (x, \df x') = \int_{\Y \times \Z \times \Y} P^*(x, \df (y,z)) \, \varpi_z(\df y') \, P((y',z), \df x')
	\]
	if $P$ and $P^*$ can be represented as integral kernels.
	
	We shall compare Dirichlet forms associated with $P^* \hat{Q} P$ to those associated with the operator $P^* \Proj^* \Proj P$.
	A spectral decomposition formula of the form \eqref{ine:generic} with $\bar{S} = P^* \Proj^* \Proj P$ would arise naturally from this comparison.
	
	\subsection{Decomposition of spectral gaps and norms}
	
	In this subsection, we derive spectral gap decomposition for $S$ satisfying \ref{H1} to \ref{H4}.
	
	First note that, if $f \in L_0^2(\pi)$, then for $\varpi$-a.e. $z \in \Z$, the function $y \mapsto Pf_z(y) := Pf(y,z)$ is in $L^2(\varpi_z)$, and $Pf_z - \varpi_z Pf_z \in L_0^2(\varpi_z)$, where $\varpi_z Pf_z = \int_{\Y} Pf_z(y) \, \varpi_z(\df y)$.
	Then we have the following lemma.
	
	\begin{lemma} \label{lem:main}
		Assume that \ref{H2} and \ref{H3} hold.
		Then, for $f \in L_0^2(\pi)$,
		\[
		\mathcal{E}_{P^*\hat{Q}P}(f) = \|f\|_{\pi}^2 - \|Pf\|_{\tilde{\pi}}^2 + \int_{\Z} \varpi(\df z) \, \mathcal{E}_{Q_z}(Pf_z - \varpi_z Pf_z).
		\]
	\end{lemma}
	
	\begin{proof}
		Let $f \in L_0^2(\pi)$ be arbitrary.
		It holds that
		\[
		\mathcal{E}_{P^*\hat{Q}P}(f) = \|f\|_{\pi}^2 - \|Pf\|_{\tilde{\pi}}^2 + \mathcal{E}_{\hat{Q}}(Pf),
		\]
		Hence, it suffices to show that
		\begin{equation} \nonumber
			\mathcal{E}_{\hat{Q}}(Pf) = \int_{\Z} \varpi(\df z) \, \mathcal{E}_{Q_z}(Pf_z - \varpi_z Pf_z).
		\end{equation}
		
		Now,
		\[
		\begin{aligned}
			\mathcal{E}_{\hat{Q}}(Pf) 
			&= \frac{1}{2} \int_{\Y \times \Z} \tilde{\pi}(\df y, \df z) \int_{\Y} Q_z(y, \df y') \, [Pf(y',z) - Pf(y,z)]^2 \\
			&= \frac{1}{2} \int_{\Z} \varpi(\df z) \int_{\Y^2} \varpi_z(\df y) \, Q_z(y, \df y') \, [Pf(y',z) - Pf(y,z)]^2 \\
			&=  \int_{\Z} \varpi(\df z) \, \mathcal{E}_{Q_z}(Pf_z - \varpi_z Pf_z)
		\end{aligned}
		\]
		as desired.
	\end{proof}
	
	We can then derive the following theorem involving Dirichlet forms.

	\begin{theorem} \label{thm:uniform}
		Assume that \ref{H1} to \ref{H4} hold.
		Suppose that there is a constant $\kappa_{\dagger} \in [0,1]$ such that $\mathrm{Gap}(Q_z) \geq \kappa_{\dagger}$ for $\varpi$-a.e. $z \in \Z$.
		Then, for $f \in L_0^2(\pi)$,
		\[
		\mathcal{E}_S(f) \geq c_0 \, \mathcal{E}_{P^*\hat{Q}P}(f) \geq c_0 \, \kappa_{\dagger} \, \mathcal{E}_{P^*\Proj^*\Proj P}(f).
		\]
	\end{theorem}
	
	\begin{proof}
		Let $f \in L_0^2(\pi)$ be arbitrary.
		By Lemma \ref{lem:main}, 
		\[
		\begin{aligned}
			\mathcal{E}_{P^* \hat{Q} P}(f) \geq
			& \|f\|_{\pi}^2 - \|Pf\|_{\tilde{\pi}}^2 + \kappa_{\dagger} \, \int_{\Z} \varpi(\df z) \, \|Pf_z - \varpi_z Pf_z\|^2_{\varpi_z} \\
			=& \|f\|_{\pi}^2 - \|Pf\|_{\tilde{\pi}}^2 + \kappa_{\dagger} \, \int_{\Z} \varpi(\df z) \, [ \|Pf_z\|^2_{\varpi_z} - (\varpi_z Pf_z)^2 ] \\
			=& \|f\|_{\pi}^2 - \|Pf\|_{\tilde{\pi}}^2 + \kappa_{\dagger} (\|Pf\|_{\tilde{\pi}}^2 - \langle \Proj Pf, \Proj Pf \rangle_{\varpi} ) \\
			=& \|f\|_{\pi}^2 - \|Pf\|_{\tilde{\pi}}^2 + \kappa_{\dagger} \, (\|Pf\|_{\tilde{\pi}}^2 - \langle f, P^* \Proj^* \Proj P f \rangle_{\pi} ) \\
			=& (1-\kappa_{\dagger}) (\|f\|_{\pi}^2 - \|Pf\|_{\tilde{\pi}}^2) + \kappa_{\dagger} \, \mathcal{E}_{P^*\Proj^*\Proj P}(f).
		\end{aligned}
		\]
		By \ref{H2}, $\|Pf\|_{\tilde{\pi}} \leq \|f\|_{\pi}$.
		Then
		\[
		\mathcal{E}_{P^* \hat{Q} P}(f) \geq \kappa_{\dagger} \, \mathcal{E}_{P^*\Proj^*\Proj P}(f).
		\]
		The desired result then follows from Lemma \ref{lem:peskun}.
	\end{proof}
	
	Based on Theorem \ref{thm:uniform}, we may immediately derive the following spectral gap decomposition result.
	
	\begin{corollary} \label{cor:uniform}
		Assume that \ref{H1} to \ref{H4} hold.
		Suppose that there is a constant $\kappa_{\dagger} \in [0,1]$ such that $\mathrm{Gap}(Q_z) \geq \kappa_{\dagger}$ for $\varpi$-a.e. $z \in \Z$.
		Then 
		\[
		\mathrm{Gap}(S) \geq c_0 \, \mathrm{Gap}(P^* \hat{Q} P)  \geq c_0 \, \kappa_{\dagger} \, \mathrm{Gap}(P^*\Proj^*\Proj P).
		\]
	\end{corollary}
	
	Corollary \ref{cor:uniform} gives \eqref{ine:generic} with $\bar{S} = P^*\Proj^*\Proj P$.
	
	{ 
	\begin{remark} \label{rem:tightness}
		The lower bound in Corollary~\ref{cor:uniform} can be tight.  For example, in Section~\ref{ssec:da}, it is shown that the data augmentation setting of Section~\ref{ssec:andrieu} falls into the general framework, with $S = P^* \hat{Q} P$ and $P^* P = I$.
		In that scenario, let $\kappa\in[0,1]$ and take
		\[
		Q_z(y, d y')=(1-\kappa) \delta_y(d y')+\kappa \varpi_z(d y'),
		\]
		so that $\mathrm{Gap}(Q_z) = \kappa$.
		Then $S=(1-\kappa)I+\kappa\bar S$, and therefore
		\[
		\mathrm{Gap}(S)=\kappa\,\mathrm{Gap}(\bar S).
		\]
		Thus, Corollary~\ref{cor:uniform} is attained with equality in this special case.
	\end{remark}
	}
	
	\begin{remark}
		Sometimes, it may be more convenient to compare $S$ to $\Proj PP^*\Proj ^*$, which has the same right spectral gap as $P^*\Proj^*\Proj P$.
	\end{remark}

	The following result complements Theorem~\ref{thm:uniform} and Corollary~\ref{cor:uniform}.
	
	\begin{proposition} \label{pro:left-gap}
		Assume that \ref{H2} and \ref{H3} hold.
		Suppose that there is a constant $\kappa_{\ddagger} \in [0,1]$ such that $\mathrm{Gap}_-(Q_z) \geq \kappa_{\ddagger}$ for $\varpi$-a.e. $z \in \Z$.
		Then, for $f \in L_0^2(\pi)$,
		\[
		\mathcal{E}_{P^*\hat{Q}P}(f) \leq  (2-\kappa_{\ddagger}) \, \mathcal{E}_{P^* \Proj^* \Proj P}(f).
		\]
		As a result, 
		\[
		\mathrm{Gap}(P^* \hat{Q} P) \leq (2- \kappa_{\ddagger}) \, \mathrm{Gap}(P^* \Proj^* \Proj P), \quad \mathrm{Gap}_-(P^* \hat{Q} P) \geq \kappa_{\ddagger}.
		\]
		In particular, if $Q_z$ is positive semi-definite for $\varpi$-a.e. $z \in \Z$ (so that one can take $\kappa_{\ddagger} = 1$), then $\mathrm{Gap}(P^* \hat{Q} P) \leq \mathrm{Gap}(P^* \Proj^* \Proj P)$, and $P^* \hat{Q} P$ is positive semi-definite.
	\end{proposition}
	
	{ 
		Proposition~\ref{pro:left-gap} implies that $P^*\Proj^*\Proj P$ is never much worse than $P^*\hat{Q}P$ in terms of the spectral gap. The factor $2 - \kappa_{\ddagger}$ approaches its maximum value~2 only when $Q_z$ has spectral components close to -1. 
		}
	
	\begin{proof} [Proof of Proposition~\ref{pro:left-gap}]
		Recall that , for a $\rho$-reversible Mtk~$J$, $$\mathrm{Gap}_-(J) = 2 - \sup_{f \in L_0^2(\rho) \setminus \{0\}} \mathcal{E}_J(f)/\|f\|_{\rho}^2.$$
		By Lemma \ref{lem:main} and \ref{H2}, for $f \in L_0^2(\pi)$,
		\[
		\begin{aligned}
			\mathcal{E}_{P^*\hat{Q}P}(f) &\leq \|f\|_{\pi}^2 - \|Pf\|_{\tilde{\pi}}^2 + (2 - \kappa_{\ddagger}) \int_{\Z} \varpi(\df z) \, \|Pf_z - \varpi Pf_z\|_{\varpi_z}^2 \\
			&= \|f\|_{\pi}^2 - \|Pf\|_{\tilde{\pi}}^2 + (2 - \kappa_{\ddagger}) (\|Pf\|_{\tilde{\pi}}^2 - \langle f, P^*\Proj^*\Proj Pf \rangle_{\pi} ) \\
			&= \|f\|_{\pi}^2 + (1-\kappa_{\ddagger})\|Pf\|_{\tilde{\pi}}^2 - (2 - \kappa_{\ddagger}) \langle f, P^*\Proj^*\Proj Pf \rangle_{\pi} \\
			&\leq  (2-\kappa_{\ddagger}) \, \mathcal{E}_{P^* \Proj^* \Proj P}(f).
		\end{aligned}
		\]
		Then $\mathrm{Gap}(P^* \hat{Q} P) \leq (2- \kappa_{\ddagger}) \, \mathrm{Gap}(P^* \Proj^* \Proj P)$.
		Since $P^* \Proj^* \Proj P$ is positive semi-definite, $\mathcal{E}_{P^* \Proj^* \Proj P}(f) \leq \|f\|_{\pi}^2$.
		It is then easy to see from the display that $\mathrm{Gap}_-(P^* \hat{Q} P) \geq \kappa_{\ddagger}$.
	\end{proof}

	Combining Corollary \ref{cor:uniform} and Proposition \ref{pro:left-gap} yields a decomposition result involving norms.
	
	\begin{corollary}
		Assume that \ref{H2} and \ref{H3} hold.
		Suppose that there is a constant $\kappa_0 \in [0,1]$ such that $1 - \|Q_z\|_{\varpi_z} \geq \kappa_0$ for $\varpi$-a.e. $z \in \Z$.
		Then
		\[
		1 - \|P^* \hat{Q} P \|_{\pi} \geq \min\{ \kappa_0 \, \mathrm{Gap}(P^* \Proj^* \Proj P), \, \kappa_0 \} =  \kappa_0 \, \mathrm{Gap}(P^* \Proj^* \Proj P),
		\]
		where the equality holds because $P^* \Proj^* \Proj P$ is positive semi-definite, so $\mathrm{Gap}(P^* \Proj^* \Proj P) \leq 1$.
	\end{corollary}

	\subsection{Non-reversible chains} \label{ssec:nonreveresible}

	It is straightforward to extend Theorem \ref{thm:uniform} and Corollary \ref{cor:uniform} to a non-reversible setting.
	
	Let $T: \X \times \A \to [0,1]$ be a possibly non-reversible Mtk such that $\pi T = \pi$, so that $T$ can be regarded as a bounded linear operator on $L_0^2(\pi)$.
	The following condition will be imposed:
	
	\begin{itemize}
		\item [\namedlabel{C1}{(C1)}]
		One can find (i) a linear transformation $P: L_0^2(\pi) \to L_0^2(\tilde{\pi})$ satisfying \ref{H2}, (ii) a collection of Mtks $\{Q_z\}_{z \in \Z}$ satisfying \ref{H3}, (iii) an Mtk $Q: (\Y \times \Z) \times (\B \times \C) \to [0,1]$ reversible with respect to $\tilde{\pi}$ satisfying \ref{H4}, and (iv) a linear transformation $R: L_0^2(\tilde{\pi}) \to L_0^2(\pi)$ such that $Q = R^*R$ and that $T = RP$.
	\end{itemize}

	Under \ref{C1}, for $f \in L_0^2(\pi)$, $\|f\|_{\pi}^2 - \|Tf\|_{\pi}^2 = \mathcal{E}_{P^* Q P}(f)$.
	Consequently, $1 - \|T\|_{\pi}^2 = \mathrm{Gap}(P^*QP)$.
	One may then immediately obtain the following result from Theorem \ref{thm:uniform}.
	
	\begin{corollary} \label{cor:nonreversible}
		Assume that \ref{C1} holds.
		Suppose that there is a constant $\kappa_{\dagger} \in [0,1]$ such that $\mathrm{Gap}(Q_z) \geq \kappa_{\dagger}$.
		Then, for $f \in L_0^2(\pi)$,
		\[
		\|f\|_{\pi}^2 - \|Tf\|_{\pi}^2 \geq c_0 \, \kappa_{\dagger} \, \mathcal{E}_{P^* \Proj^* \Proj P}(f).
		\]
		Consequently,
		\[
		1 - \|T\|_{\pi}^2 \geq c_0 \, \kappa_{\dagger} \mathrm{Gap}(P^* \Proj^* \Proj P).
		\]
	\end{corollary}
	
	\section{Special Cases} \label{sec:special}
	
	\subsection{Overview}
	
	We demonstrate that all the spectral gap decomposition results listed in Section \ref{sec:examples} can be derived from Corollaries \ref{cor:uniform} or \ref{cor:nonreversible}, and are thus unified within the framework in Section~\ref{sec:unified}.
	The following derivations are performed:
	\begin{itemize}
		\item Propositions~\ref{pro:madras} and \ref{pro:madras-1} (concerning decomposition over a finite partition or cover), along with Proposition~\ref{pro:andrieu} (concerning data augmentation algorithms with one intractable conditional) are established directly from Corollary~\ref{cor:uniform}.
		\item Propositions~\ref{pro:qin} (concerning random-scan Gibbs samplers), \ref{pro:chen} (concerning localization schemes), \ref{pro:hit} (concerning hit-and-run algorithms) are derived from Proposition~\ref{pro:andrieu}.
		\item Proposition~\ref{pro:doubly} (concerning data augmentation algorithms with two intractable conditionals) is derived from Corollary \ref{cor:nonreversible}.
	\end{itemize}

	
	To derive the desired results within the framework herein, it suffices to identify the objects in Section \ref{sec:unified} and check the relevant conditions, i.e., \ref{H1} to \ref{H4}.
	Table \ref{tab:summary} summarizes the identification for each example.
	In Section \ref{sec:special}, we provide the detailed constructions for deriving Propositions~\ref{pro:madras} and \ref{pro:andrieu}.
	The constructions for other examples are relegated to Appendix~\ref{app:derive}.
	
	\begin{table} \caption{Important elements from Section \ref{sec:unified} in the context of the examples found in Sections \ref{ssec:madras} to \ref{ssec:doubly-0}.
	The first column lists the relevant proposition numbers.
	The last column indicates whether $Q = \hat{Q}$ in the sandwich construction.
	} \label{tab:summary}
		\centering
		\begin{tabular}{llllll}
			\hline 
			\makecell[l]{Propo- \\sition	} & $P^* Q P$ & \makecell[l]{$P^* \Proj^* \Proj P$ \\ or \\ $\Proj P P^* \Proj^*$}  & $\Z$ & $Q_z$ & $Q = \hat{Q}$? \\
			\hline
			\ref{pro:madras} & $M^{1/2} N M^{1/2}$ & $M_0$ & $[k]$ & \makecell[l]{$Q_z(x, \mathsf{A}) = $  \\
			$N(x, \mathsf{A} \cap \X_z) +$ \\ $N(x, \X_z^c) \, \delta_x(\mathsf{A})$ } & No \\
			
			\vspace{0.3cm}
			
			\ref{pro:madras-1} & see note & see note & $[k]$ & same as above & No \\
			
			\vspace{0.3cm}
			
			\ref{pro:andrieu} & $S$ & $\bar{S}$ & $\Z$ & $Q_z$ & Yes \\

			\vspace{0.3cm}
			
			\ref{pro:qin} & $S$ & $\bar{S}$ & $\bigcup_{i=1}^k \{i\} \times \X_{-i}$ & 
		\makecell[l]{$Q_{i,u}(x, \df x') =$  \\ $H_{i,u}(x_i, \df x_i') \, \delta_u(\df x_{-i}')$}  & Yes \\
			
			\vspace{0.3cm}
			
			\ref{pro:chen} & $K_{t+1}$ & $K_t$ & $\mathsf{W}_0 \times \cdots \times \mathsf{W}_t$ & \makecell[l]{$Q_z(x, \mathsf{A}) = $ \\ $\mathrm{E} \left[ \frac{\df \nu_{s+1}}{\df \nu_s}(x) \, \nu_{s+1}(\mathsf{A}) \mid \right.$ \\
				$\left.  (W_s)_{s=0}^t = z \right]$} & Yes \\
			
			\vspace{0.3cm}
			
			\ref{pro:hit} & $S$ & $\bar{S}$ & $\mathsf{W} \times \X$ & \makecell[l]{$Q_{w,x'}(x, \mathsf{A}) = $ \\
			$ \int_{\mathbb{R}^{\ell}} H_{w, \theta(x',w)} \left( (w_i^{\top} x)_{i=1}^{\ell}, \df u \right)  $ \\
			$ \ind_{\mathsf{A}} \left( \theta(x',w) + \sum_{i=1}^{\ell} u_i w_i \right) $} & Yes \\
			
			\vspace{0.3cm}
			
			\ref{pro:doubly} & $T^* T$ & $\hat{S}_2$ & $\X_2$ & $H_{1,x_2}^2$ & Yes \\
			
			\hline
			
			\multicolumn{6}{l}{ \makecell[l]{Note: For Proposition \ref{pro:madras-1}, \\
			$P^* Q P(x, \df x') = \Theta^{-1} N(x, \df x') \, L(x') + [1 - \Theta^{-1} \int_{\X} N(x, \df x'') \, L(x'') ] \, \delta_x(\df x')$, \\
			$\Proj P P^* \Proj^*(z, \{z'\}) = \pi_0(\X_z)^{-1} \int_{\X_z \cap \X_{z'}} L(x)^{-1} \pi_0(\df x)$, \\
			where $L(x) = \sum_{z=1}^k \ind_{\X_z}(x)$.
			} } \\
			\hline
		\end{tabular}
	\end{table}

	\subsection{Markov chain decomposition based on state space partitioning} \label{ssec:decomposition}
	
	Recall the setting in Section \ref{ssec:madras}:
	The space $\X$ has a partition $\X = \bigcup_{z=1}^k \X_z$, where $\pi(\X_z) > 0$.
	$M$ and $N$ are Mtks that are reversible with respect to~$\pi$, and $M$ is positive semi-definite.
	Moreover, the Mtks $M_0: [k] \times 2^{[k]} \to [0,1]$ and $H_z: \X_z \times \A_z \to [0,1]$ (where $\A_z$ is $\A$ restricted to $\X_z$) are defined via the formulas below:
	\[
		M_0(z,\{z'\}) = \frac{1}{\pi(\X_z)} \int_{\X_z} \pi(\df x) \, M(x, \X_{z'}), \quad H_z(x, \mathsf{A}) = N(x, \mathsf{A}) + N(x, \X_z^c) \, \delta_x(\mathsf{A}).
	\]
	$M_0$ is reversible with respect to $\varpi$, where $\varpi(\{z\}) = \pi(\X_z)$ for $z \in [k]$.
	$H_z$ is reversible with respect to $\omega_z$, where $\omega_z(\mathsf{A}) = \pi(\mathsf{A})/\pi(\X_z)$ for $\mathsf{A} \in \A_z$.
	
	Our goal is to derive Proposition \ref{pro:madras} using Corollary~\ref{cor:uniform}.
	To this end, we identify the elements in Section \ref{ssec:sandwich} as follows:
	
	\begin{itemize}
		\item $(\Y, \B) = (\X, \A)$, $(\Z, \C) = ([k], 2^{[k]})$.
		
		\item $\tilde{\pi}: \A \times 2^{[k]} \to [0,1]$ is such that
		\[
		\tilde{\pi}(\mathsf{A} \times \{z\}) = \pi(\mathsf{A} \cap \X_z) = \int_{\mathsf{A}} \pi(\df x) \, \ind_{\X_z}(x).
		\]
		In other words, $\tilde{\pi}$ is the distribution of an $\X \times [k]$-valued random element $(X,Z)$ where $X \sim \pi$ and $Z$ satisfies $X \in \X_Z$.
		
		\item For $z \in [k]$ and $\mathsf{A} \in \A$, $\varpi_z(\mathsf{A}) = \omega_z(\mathsf{A} \cap \X_z) = \pi(\mathsf{A} \cap \X_z)/\pi(\X_z)$.
		Note that, if $(X, Z) \sim \tilde{\pi}$, then $Z \sim \varpi$, and $X \mid Z = z \sim \varpi_z$ as required.
		
		\item 
		$S = M^{1/2} N M^{1/2}$.
		
		\item For $f \in L_0^2(\pi)$ and $(x,z) \in \X \times [k]$, $P f(x,z) = M^{1/2} f(x)$.
		It is easy to check $P$ is a linear transformation from $L_0^2(\pi)$ to $L_0^2(\tilde{\pi})$.
		
		\item For $h \in L_0^2(\tilde{\pi})$ and $x \in \X$, $P^*h(x) = M^{1/2} \Projj h(x)$, where $\Projj h(x') = \sum_{z=1}^k h(x',z) \, \ind_{\X_z}(x')$ for $x' \in \X$.
		
		\item The Mtk $Q: (\X \times [k]) \times (\A \times 2^{[k]}) \to [0,1]$ is of the form
		\[
		Q((x,z), \mathsf{A} \times \{z'\}) = N(x, \X_{z'} \cap \mathsf{A}) .
		\]
		It is easy to verify that $Q$ is reversible with respect to $\tilde{\pi}$.

		\item For $z \in [k]$, $x \in \X$, and $\mathsf{A} \in \A$,
		\[
		Q_z(x, \mathsf{A}) = N(x, \mathsf{A} \cap \X_z) + N(x, \X_z^c) \, \delta_x(\mathsf{A}).
		\]
		Note that, if $x \in \X_z$ and $\mathsf{A} \in \A_z$, then $Q_z(x, \mathsf{A}) = H_z(x, \mathsf{A})$.
		
		\item For $h \in L_0^2(\tilde{\pi})$ and $z \in [k]$,
		\[
		\Proj h (z) = \frac{1}{\pi(\X_z)} \int_{\X_z} h(x,z) \, \pi(\df x).
		\]
		
	\end{itemize}
	
	One may then verify the following:
	
	\begin{itemize}
		\item $P^* Q P = M^{1/2} N M^{1/2}$, and \ref{H1} holds.
		
		\item $P^* P = M$, and \ref{H2} holds.
		
		\item $Q_z$ is reversible with respect to $\varpi_z$, so \ref{H3} holds.
		
		\item For $z \in [k]$ and $x \in \X_z$ (which account for $\tilde{\pi}$-a.e. possible value of $(x,z)$) and $\mathsf{A} \in \A$,
		\[
		\begin{aligned}
			Q((x,z), (\mathsf{A} \setminus \{x\}) \times \{z\}) &= N(x, \X_z \cap \mathsf{A} \setminus \{x\}) = Q_z(x, \mathsf{A} \setminus \{x\}).
		\end{aligned}
		\]
		Thus,~$Q$ is approximately $z$-invariant, and \ref{H4} holds with $c_0 = 1$.
		
		\item 
		For $g \in L_0^2(\varpi)$ and $z \in [k]$,
		\[
		\begin{aligned}
			\Proj P P^* \Proj^* g (z) &= \int_{\X} \varpi_z(\df x) \, \int_{\X} M(x, \df x') \sum_{z'=1}^k g(z') \, \ind_{\X_{z'}}(x') \\
			&= \sum_{z' = 1}^k \frac{1}{\pi(\X_z)} \int_{\X_z} \pi(\df x) \, M(x, \X_{z'}) \, g(z') \\
			&= M_0 \, g(z).
		\end{aligned}
		\]
		Thus, $ \Proj P P^* \Proj^* = M_0$.
		
	\end{itemize}
	
	It holds that $\mathrm{Gap}(P^* \Proj^* \Proj P) = \mathrm{Gap}(\Proj P P^* \Proj^*)$.
	By Corollary~\ref{cor:uniform},
	\[
	\mathrm{Gap}(M^{1/2} N M^{1/2}) \geq \left[ \essinf_z \mathrm{Gap}(Q_z) \right] \, \mathrm{Gap}(M_0).
	\]
	A simple argument shows that $\mathrm{Gap}(Q_z) = \mathrm{Gap}(H_z)$; see Lemma~\ref{lem:iso-1} in Appendix~\ref{app:lemma}.
	Proposition~\ref{pro:madras} then follows.

	\subsection{Hybrid data augmentation algorithms with one intractable conditional} \label{ssec:da}
	
	Recall the setting of Section \ref{ssec:andrieu}:
	$(\Z,\C)$ is a measurable space.
	On the space $\X \times \Z$, a joint distribution has the form
	\[
	\tilde{\pi}(\df (x,z)) = \pi(\df x) \, \pi_x(\df z) = \varpi(\df z) \, \varpi_z(\df x).
	\]
	$\{Q_z\}_{z \in \Z}$ is a collection of Mtks satisfying \ref{H3}.
	The Mtks $\bar{S}$ and $S$ are defined via the following formulas:
	\[
	\bar{S}(x, \df x') = \int_{\Z} \pi_x(\df z) \, \varpi_z(\df x'), \quad S(x, \df x') = \int_{\Z} \pi_x(\df z) \,  Q_z(x, \df x').
	\]
	Our goal is to derive Proposition \ref{pro:andrieu} using Corollary \ref{cor:uniform}.
	
	The elements of Section \ref{ssec:sandwich} that are not already specified above are identified as follows:
	\begin{itemize}
		\item $(\Y, \B) = (\X, \A)$.
		
		\item For $f \in L_0^2(\pi)$ and $(x,z) \in \X \times \Z$, let $Pf(x,z) = f(x)$.
		
		\item For $h \in L_0^2(\tilde{\pi})$ and $x \in \X$, $P^* h(x) = \int_{\Z} h(x, z) \, \pi_x(\df z)$.
		
		\item Let $Q: (\X \times \Z) \times (\A \times \C) \to [0,1]$ be an Mtk such that
		\[
		Q((x,z), \df (x',z')) = Q_z(x, \df x') \, \delta_z(\df z').
		\]
		
		\item For $h \in L_0^2(\tilde{\pi})$ and $z \in \Z$, $\Proj h(z) = \int_{\X} h(x,z) \, \varpi_z(\df x)$, as defined in \eqref{eq:Proj}.
	\end{itemize}
	
	Then one may easily verify:
	\begin{itemize}
		\item $S = P^*QP$, and \ref{H1} holds.
		\item $\|Pf\|_{\tilde{\pi}} = \|f\|_{\pi}$, so \ref{H2} holds.
		\item \ref{H3} and \ref{H4} are satisfied with $c_0 = 1$.
		
		\item $P^* \Proj^* \Proj P = \bar{S}$.
	\end{itemize}

	Proposition \ref{pro:andrieu} then follows from Theorem \ref{thm:uniform} and Corollary \ref{cor:uniform}.

	\section{Some Properties of the Sandwich Structure} \label{sec:simple}
	
	Going back to the generic setting in Section \ref{sec:unified}, we now establish some secondary results within this general framework.
	
	\subsection{Some simple results involving $P^*\hat{Q}P$}
	
	As indicated in Table~\ref{tab:summary}, in many situations, $Q = \hat{Q}$ in the sandwich structure, so $S = P^* \hat{Q} P$.
	This is the case for hybrid data augmentation, random-scan Gibbs, and hit-and-run samplers.
	In this subsection, we provide some simple results regarding $P^* \hat{Q} P$.
	
	\subsubsection{Loewner ordering involving different choices of $Q_z$}
	
	Assume that, given $z \in \Z$, there are two choices of the Mtk $Q_z$, say $Q_z^{(1)}$ and $Q_z^{(2)}$, both reversible with respect to $\varpi_z$.
	For $j = 1,2$, let
	\[
	Q_{(j)}((y,z), \df (y',z')) = Q_z^{(j)}(y, \df y') \, \delta_z(\df z').
	\]
	Let $S_{(j)} = P^* Q_{(j)} P$.
	We see that $S_{(j)}$ can be viewed as a version of $P^* \hat{Q} P$ by taking $\hat{Q} = Q_{(j)}$.
	In hybrid Gibbs-like algorithms (e.g., data augmentation, random-scan Gibbs, and hit-and-run algorithms), having two choices of $Q_z$ corresponds to having two choices of Markovian approximations to approximate conditional distributions.
	See Table~\ref{tab:summary}.
	
	For two self-adjoint linear operators $K_{(1)}$ and $K_{(2)}$ defined on the same Hilbert space, write $K_{(1)} \leq K_{(2)}$ if $K_{(2)} - K_{(1)}$ is positive semi-definite.
	This is called a Loewner ordering.
	It is worth mentioning that Loewner ordering of Mtks is implied by Peskun ordering \citep{peskun1973optimum,tierney1998note}.
	
	The following result is easy to establish.
	
	\begin{proposition} \label{pro:loewner}
		Assume that \ref{H2} and \ref{H3} hold.
		Assume further that, for $\varpi$-a.e. $z \in \Z$, $Q_z^{(1)} \leq Q_z^{(2)}$.
		Then $S_{(1)} \leq S_{(2)}$.
		Equivalently, for $f \in L_0^2(\pi)$, $\mathcal{E}_{S_{(1)}}(f) \geq \mathcal{E}_{S_{(2)}}(f)$.
		In particular, $\mathrm{Gap}(S_{(1)}) \geq \mathrm{Gap}(S_{(2)})$.
	\end{proposition}
	
	\subsubsection{Iterated implementation of $Q_z$} \label{sssec:iteration}

	Since $Q_z$ is an Mtk, it makes sense to consider iterating it in each step of an algorithm.
	Note that the number of iterations may depend on~$z$.

	Let $\lambda: \Z \to \mathbb{Z}_+$ be a measurable function.
	Let 
	\[
	Q_{\lambda}((y,z), \df (y',z')) = Q_z^{\lambda(z)}(y, \df y') \, \delta_z(\df z').
	\]
	We briefly investigate the behavior of the operator $S_{\lambda} = P^* Q_{\lambda} P$.
	
	Intuitively, it may be beneficial to iterate $Q_z$ more times for values of~$z$ such that $\mathrm{Gap}(Q_z)$ is small, provided that it is not too costly to do so.
	The following simple corollary of Theorem \ref{thm:uniform} is established in Appendix~\ref{app:proofs}.
	
	\begin{proposition} \label{pro:iteration-1}
		Suppose that \ref{H2} and \ref{H3} hold.
		Suppose further that $\mathrm{Gap}(Q_z) > 0$ for $\varpi$-a.e. $z \in \Z$.
		Let $\kappa_{\dagger} \in (0,1)$ be an arbitrary constant, and let $\lambda: \Z \to \mathbb{Z}_+$ be such that
		\begin{equation} \nonumber
			\lambda(z) \geq \frac{\log (1-\kappa_{\dagger})}{\log[1 - \min\{1,\mathrm{Gap}(Q_z)\}]}, \quad \text{which holds if} \quad \lambda(z) \geq - \frac{\log(1-\kappa_{\dagger})}{\min\{1,\mathrm{Gap}(Q_z)\}},
		\end{equation}
		for $\varpi$-a.e. $z \in \Z$.
		(We use the convention that $1/\log 0 = 0$.)
		Finally, assume that, for $\varpi$-a.e. $z \in \Z$, either $\mathrm{Gap}_-(Q_z) \geq \mathrm{Gap}(Q_z)$ (which holds if $Q_z$ is positive semi-definite), or $\lambda(z)$ is odd.
		Then, for $f \in L_0^2(\pi)$,
		\[
		\mathcal{E}_{S_{\lambda}}(f) \geq \kappa_{\dagger} \, \mathcal{E}_{P^* \Proj^* \Proj P}(f).
		\]
		In particular,
		\[
		\mathrm{Gap}(S_{\lambda}) \geq \kappa_{\dagger} \mathrm{Gap}(P^* \Proj^* \Proj P).
		\]
	\end{proposition}
	
	{  In the statement of Proposition~\ref{pro:iteration-1}, we have used the fact that $-\log(1-x) \geq x$ for $x \in (0,1)$ to argue that $\lambda(z) \geq -\log(1-\kappa^{\dagger})/\min\{1,\mathrm{Gap}(Q_z)\}$ leads to the condition on $\lambda(z)$.} 
	
	{ The condition involving odd $\lambda(z)$ is needed only to control the negative spectrum of $Q_z$.  An odd power leaves a negative eigenvalue negative, so it cannot move to the upper spectral edge and reduce the right spectral gap; an even power could turn an eigenvalue near $-1$ into one near~$1$, causing $\mathrm{Gap}(Q_z^{\lambda(z)})$ to become very small even when $\mathrm{Gap}(Q_z)$ is not.}
	
	Proposition~\ref{pro:iteration-1} shows that even if $\essinf_z \mathrm{Gap}(Q_z) = 0$, as long as $\mathrm{Gap}(Q_z) > 0$ for $\varpi$-a.e. $z \in \Z$, one may construct an algorithm $S_{\lambda}$ such that $\mathrm{Gap}(S_{\lambda})/\mathrm{Gap}(P^* \Proj^* \Proj P)$ has a nonzero lower bound.
	Loosely speaking, letting $\lambda(z) \propto 1/\mathrm{Gap}(Q_z)$ would suffice.

	\subsection{Weak Poincar\'{e} inequalities} \label{ssec:weak}
	
	In the spectral decomposition formula, if $\mathrm{Gap}(Q_z) = 0$ for a nonzero measure set of~$z$, then even iterating $Q_z$ would not prevent $\mathrm{Gap}(Q_z^{\lambda(z)})$ from having a vanishing essential infimum.

	Relaxing the requirement of $\mathrm{Gap}(Q_z) > 0$ is challenging in general.
	Here we discuss one possible strategy based on techniques from \cite{andrieu2022comparison} involving weak Poincar\'{e} inequalities.
	This technique was used in \cite{power2024weak} to study hybrid slice samplers.
	For similar techniques, see \cite{ascolani2024scalability}, which studied hybrid random-scan Gibbs samplers via the $s$-conductance, and \cite{atchade2021approximate}, which extended decomposition results from \cite{madras2002markov} using approximate spectral gaps.
	
	Let $K: \Omega \times \mathcal{F} \to [0,1]$ be an Mtk that is reversible with respect to a probability measure~$\rho$.
	We say it satisfies a weak Poincar\'{e} inequality if there exists a non-increasing function $\beta: (0,\infty) \to [0,\infty)$ satisfying $\beta(u) \downarrow 0$ as $u \to \infty$ such that, for each $f \in L_0^2(\rho)$ and $u > 0$,
	\begin{equation} \nonumber
		\|f\|_{\rho}^2 \leq u \, \mathcal{E}_K(f) + \beta(u) \|f\|_{\scriptsize\mathrm{osc}}^2,
	\end{equation}
	where $\|f\|_{\scriptsize\mathrm{osc}} = \esssup_w f(w) - \essinf_w f(w)$, with essential supremum and infimum defined in terms of~$\rho$.
	
	An Mtk may satisfy a weak Poincar\'{e} inequality even if it does not admit a positive right spectral gap.
	Moreover, a weak Poincar\'{e} inequality leads to a quantitative convergence bound, often of subgeometric nature.
	We refer readers to \cite{andrieu2022comparison} and \cite{power2024weak}.
	
	The following two results, established in Appendix~\ref{app:proofs}, can be regarded as analogues of Theorem \ref{thm:uniform} in terms of weak Poincar\'{e} inequalities.
	
	\begin{lemma} \label{lem:weak-new}
		Assume that \ref{H1} to \ref{H4} hold.
		Assume also that, for $f \in L_0^2(\pi)$, it holds that $\|Pf\|_{\scriptsize\mathrm{osc}} \leq \|f\|_{\scriptsize\mathrm{osc}}$.
		Suppose further that there is a measurable function $\alpha: \Z \times (0,\infty) \to [0,\infty)$ such that, for $\varpi$-a.e. $z \in \Z$, $u > 0$, and $g \in L_0^2(\varpi_z)$, 
		\[
		\|g\|_{\varpi_z}^2 \leq u \, \mathcal{E}_{Q_z}(g) + \alpha(z,u) \, \|g\|_{\scriptsize\mathrm{osc}}^2.
		\]
		Then, for $f \in L_0^2(\pi)$ and $u \geq 1$,
		\[
		\mathcal{E}_{P^* \Proj^* \Proj P}(f) \leq u \, c_0^{-1} \mathcal{E}_S(f) + \bar{\alpha}(u) \, \|f\|_{\scriptsize\mathrm{osc}}^2.
		\]
		where $\bar{\alpha}(u) = \int_{\Z} \alpha(z,u) \, \varpi(\df z)$.
	\end{lemma}

	\begin{remark}
		If $P$ is associated with an Mtk defined on $(\Y \times \Z) \times \A \to [0,1]$, then obviously, $\|Pf\|_{\scriptsize\mathrm{osc}} \leq \|f\|_{\scriptsize\mathrm{osc}}$ for $f \in L_0^2(\pi)$.
		{  In several of our examples, $P$ is indeed an Mtk.
			This includes the hybrid data augmentation algorithm with one intractable conditional distribution (see Section~\ref{ssec:da}), and by extension, the random-scan Gibbs algorithm, the approximate variance conservation in localization schemes, and the hybrid hit-and-run sampler.
			This also includes the Markov chain decomposition over an overlapping cover (see Section~\ref{ssec:decomposition-1}).
			The assumption is not automatic
			for the partition construction in (see Section~\ref{ssec:decomposition}), where
			$P$ is defined through an operator square root rather than through an Mtk.}
	\end{remark}
	
	Based on Lemma \ref{lem:weak-new}, one may use Theorem 33 from \cite{andrieu2022comparison} to derive a weak Poincar\'{e} inequality for~$S$ based on one for $P^* \Proj^* \Proj P$.
	This is given the following result, whose proof is provided in Appendix \ref{app:proofs} for completeness.
	
	\begin{proposition} \label{pro:weak}
		In addition to the assumptions in Lemma~\ref{lem:weak-new}, assume that $P^* \Proj^* \Proj P$ satisfies a weak Poincar\'{e} inequality with a non-increasing function $\beta: (0,\infty) \to [0,\infty)$ such that $\beta(u) \downarrow 0$ as $u \to \infty$.
		Then, for $f \in L_0^2(\pi)$ and $u > 0$,
		\[
		\|f\|_{\pi}^2 \leq u \, \mathcal{E}_S(f) + \tilde{\beta}(u) \, \|f\|_{\scriptsize\mathrm{osc}}^2,
		\]
		where $\tilde{\beta}(u) = \max\{1/4, \beta(c_0) + \bar{\alpha}(1) \}$ when $u < 1$, and $\tilde{\beta}(u) = \inf\{ \beta(u_1) + u_1 \, \bar{\alpha}(u_2): \, u_1 u_2 = c_0 u, \, u_1 > 0, \, u_2 > 1\}$ when $u \geq 1$.
		Furthermore, if $\bar{\alpha}(u)$ is a non-increasing function such that $\bar{\alpha}(u) \downarrow 0$ as $u \to \infty$, then so is $\tilde{\beta}(u)$.
	\end{proposition}

	\section{Hybrid Hit-and-Run for Well-Conditioned Distributions} \label{sec:example-hybrid}
	
	In this section, we illustrate spectral gap decomposition via a more concrete example.
	We study a version of the hybrid hit-and-run sampler described in Section \ref{ssec:hit-0} when the target density~$\dot{\pi}$ is well-conditioned.
	
	\subsection{Quantitative bounds} \label{ssec:hit-quant}

	In Section \ref{ssec:hit-0}, take $\ell = 1$ for simplicity.
	
	The hit-and-run sampler has Mtk
	\[
	\begin{aligned}
		\bar{S}(x, \mathsf{A}) &= \int_{\mathsf{W}} \nu(\df w) \int_{\mathbb{R}} \varphi_{w, \theta(x,w)}(\df u) \, \ind_{\mathsf{A}}\left(\theta(x,w) + uw \right).
	\end{aligned}
	\]
	One may make draws from the one dimensional distribution $\varphi_{w, \theta(x,w)}$ using, say, rejection sampling; see \cite{chewi2022query} for a discussion on the matter.
	
	{ 
	Alternatively, even in dimension~1, it is often simpler to use a Metropolis-Hastings-type algorithm to sample from $\varphi_{w, \theta(x,w)}$.
	We consider a multi-scale Gaussian random-walk Metropolis-Hastings algorithm.
	Let $\sigma_1, \dots, \sigma_J$ be positive numbers, where~$J$ is a positive integer.
	Let $(p_1,\dots,p_J)$ be probability vector.
	Given the current state $u \in \mathbb{R}$, the algorithm chooses an index $j \in \{1,\dots,J\}$ with probability $p_j$.
	Then, a new state is proposed according to the distribution $\mathrm{N}(u, \sigma_j^2)$.
	The new state is accepted or rejected via a Metropolis step, and the acceptance probability $a(u,u')$ is specified below.
	This algorithm has Mtk
	\[
	\begin{aligned}
		&H_{w, \theta(x,w)}(u, \df u') \\
		=& \sum_{j=1}^J p_j H_{w, \theta(x,w),j}(u, \df u') \\
		=:& \sum_{j=1}^J p_j \left\{ \frac{1}{\sigma_j} \phi \left( \frac{u' - u}{\sigma_j} \right) a(u,u') \, \df u' + \left[ 1 - \int_{\mathbb{R}} \frac{1}{\sigma_j} \phi \left( \frac{u'' - u}{\sigma_j} \right) a(u,u'') \, \df u '' \right] \, \delta_u(\df u') \right\},
	\end{aligned}
	\]
	where $\phi$ is the probability density function of the $\mbox{N}(0, 1^2)$ distribution, and
	\[
	a(u,u') = \min \left\{ 1, \frac{\dot{\pi}( \theta(x,w) + u'w ) }{ \dot{\pi}( \theta(x,w) + u w ) } \right\}.
	\]
	The resultant hybrid hit-and-run sampler has Mtk
	\[
	S(x, \mathsf{A}) = \int_{\mathsf{W}} \nu(\df w) \int_{\mathbb{R}} H_{w, \theta(x,w)}(w^{\top} x, \df u) \, \ind_{\mathsf{A}}\left(\theta(x,w) + uw \right).
	\]
	
	
	}
	
	Using Proposition \ref{pro:hit} and recent results from \cite{ascolani2024entropy,secchi2025spectral}, one can establish the following quantitative bound on $\mathrm{Gap}(S)$ when $\dot{\pi}$ is in some sense well-conditioned.
	
	{ 
	
	\begin{proposition} \label{pro:hit-concave}
		Assume that $x \mapsto -\log \dot{\pi}(x)$ is twice differentiable, and the eigenvalues of its Hessian are bounded between two positive numbers, $m$ and $L$, where $m \leq L$.
		In the multiscale Gaussian random-walk Metropolis-Hastings algorithm targeting $\varphi_{w, \theta(x,w)}$, let $J = 1 + \lceil \log_2 (L/m) \rceil$, $p_j = 1/J$, and $\sigma_j^2 = 2^{j-1}/L$ for $j = 1,\dots,J$. 
		Then, for $f \in L_0^2(\pi)$,
		\begin{equation} \label{ine:hit-relation}
		 \frac{ C_* }{ 1 + \lceil \log_2 (L/m) \rceil } \mathcal{E}_{\bar{S}}(f) \leq \mathcal{E}_S(f) \leq \mathcal{E}_{\bar{S}}(f).
		\end{equation}
		where $C_* \approx 10^{-14}$ is some universal constant.
		In particular,
		\begin{equation} \nonumber
		\frac{ C_* }{ 1 + \lceil \log_2 (L/m) \rceil }\, \mathrm{Gap}(\bar{S}) \leq \mathrm{Gap}(S) \leq \mathrm{Gap}(\bar{S}).
		\end{equation}
		Moreover, if $\nu$ is the uniform distribution on the unit sphere $\mathsf{W}$, then
		\begin{equation} \label{ine:hit-gap}
		\mathrm{Gap}(S) \geq \frac{C_*}{2 (1 + \lceil \log_2 (L/m) \rceil)} \frac{m}{L} \frac{1}{k}.
		\end{equation}
	\end{proposition}
	
	\begin{proof}
		For $w \in \mathsf{W}$, $x \in \mathbb{R}^k$, the negative log-density of $\varphi_{w, \theta(x,w)}$ is twice-differentiable, and its second derivative always lies in $[m,L]$.
		By the Brascamp--Lieb inequality and the Cram\'{e}r--Rao lower bound \citep[see, e.g.,][Lemma~7]{chewi2023entropic}, the variance of $\varphi_{w, \theta(x,w)}$, denoted by $v_{w,\theta(x,w)}$, lies in $[1/L,1/m]$.
		Then, for some $i \in \{1,\dots,J\}$, 
		\[
		v_{w,\theta(x,w)} \leq \sigma_i^2 \leq 2 v_{w,\theta(x,w)}.
		\]
		By Theorem~1 of \cite{secchi2025spectral} and Cheeger's inequality \citep{lawler1988bounds}, $\mathrm{Gap}(H_{w, \theta(x,w),i})$ is no less than a positive universal constant $C_*$, where
		\[
		C_* = \frac{1}{2} \inf_{c \in [1,2]} \min\left\{ \frac{b(c)}{8}, \frac{b(c)^2 \sqrt{c}}{96\sqrt{3}} \right\}^2 \approx 10^{-14}, \quad b(c) = \frac{1}{2} e^{-2.6} \frac{1}{\sqrt{2\pi c}} \exp \left( - \frac{1.5^2}{2c} \right).
		\]
		On the other hand, for $f \in L_0^2( \varphi_{w, \theta(x,w)} )$,
		\[
		\mathcal{E}_{H_{w, \theta(x,w)}}(f) = \frac{1}{J} \sum_{j=1}^J \mathcal{E}_{H_{w, \theta(x,w), j}}(f) \geq \frac{1}{J} \mathcal{E}_{H_{w,\theta(x,w),i}}(f).
		\]
		Taking the infimum with respect to~$f$ shows that
		\[
		\mathrm{Gap}(H_{w, \theta(x,w)}) \geq \frac{C_*}{J} = \frac{C_*}{1 + \lceil \log_2(L/m) \rceil}.
		\]
		By Proposition \ref{pro:hit}, for $f \in L_0^2(\pi)$,
		\[
		\mathcal{E}_S(f) \geq \left[ \inf_{w,x} \mathrm{Gap}(H_{w, \theta(x,w)} ) \right] \, \mathcal{E}_{\bar{S}}(f).
		\]
		Combining the two most recent displays yields the lower bound in \eqref{ine:hit-relation}.

		By Lemma 3.1 of \cite{baxendale2005renewal}, Gaussian random-walk Metropolis algorithms have positive semi-definite kernels. 
		Thus, $H_{w, \theta(x,w),j}$, $j \in \{1,\dots,J\}$, is positive semi-definite.
		Then $H_{w, \theta(x,w)}$ is positive semi-definite.
		Applying Proposition~\ref{pro:left-gap} with $P^* \hat{Q} P = S$ and $ P^* \Proj^* \Proj P = \bar{S}$ (see Table~\ref{tab:summary}) yields the upper bound in \eqref{ine:hit-relation}.

		In a recent work by \cite{ascolani2024entropy}, it was shown that, when $\nu$ is the uniform distribution on $\mathsf{W}$, $\mathrm{Gap}(\bar{S}) \geq (1/2) (m/L) (1/k)$. 
		Then \eqref{ine:hit-gap} follows from \eqref{ine:hit-relation}.
	\end{proof}

	\begin{remark}
		The very small universal constant $C_*$ in Proposition~\ref{pro:hit-concave} should be viewed as a certificate of uniform positivity of the ratio 
		\[
		\frac{[1 + \lceil \log_2(L/m) \rceil] \, \mathrm{Gap}(S)}{\mathrm{Gap}(\bar{S})},
		\]
		not as a realistic estimate of this quantity.
		The precise magnitude of this constant is overwhelmingly an artifact of the proof of Theorem~1 from \cite{secchi2025spectral}.
		It mainly originates from a conservative estimate of the worst-case acceptance rate of a one-dimensional Metropolis-Hastings algorithm. 
		Very small constants like $C_*$ are common when MCMC convergence proofs track every universal factor; see \cite{andrieu2024explicit,liu2026spectral}.
	\end{remark}
	
	Based on Proposition~\ref{pro:hit}, \eqref{ine:mixing}, and \eqref{ine:asymptotic-variance}, we have the following mixing time and asymptotic variance bounds.
	\begin{corollary}
		In the context of Proposition~\ref{pro:hit-concave}, assume that $\nu$ is uniform on $\mathsf{W}$.
		Then each of the following statements holds.
		For $\epsilon > 0$ and initial distribution~$\mu$ such that $\df \mu/\df \pi \in L^2(\pi)$, we have the mixing time bound
		\[
		t_*(S, \mu,\epsilon) \leq \frac{2 k (L/m) [1 + \lceil \log_2 (L/m) \rceil] (\log\|\df \mu/\df \pi\|_{\pi} - \log \epsilon)}{C_*}.
		\]
		Moreover, for $f \in L^2(\pi)$, 
		\[
		\begin{aligned}
			\mathrm{var}_S(f) &\leq \left\{ \frac{4}{C_*} k (L/m) [1 + \lceil \log_2 (L/m) \rceil] - 1 \right\} \|f - \pi f\|_{\pi}^2.
		\end{aligned}
		\]
	\end{corollary}
	
	A recent complementary analysis by \cite{bou2025accelerated} derives sharp, geometry-sensitive Wasserstein contraction rates for exact hit-and-run on Gaussian targets.  In some examples, their rates improve from the diffusive order $(L/m)^{-1}$ to ballistic order.  Their result is sharper for those Gaussian geometries, whereas Proposition~\ref{pro:hit-concave} applies to all $m$-strongly log-concave and $L$-log-smooth targets and, through \eqref{ine:hit-relation}, transfers the convergence of exact hit-and-run to an implementable Metropolis-within-hit-and-run kernel.  The two bounds also concern different metrics: Wasserstein contraction in their work and $L^2$ spectral gaps and Dirichlet forms here.
	}

	To illustrate Proposition \ref{pro:hit-concave}, we let $\dot{\pi}$ be the posterior density for a Bayesian logistic regression model.
	To be specific, we consider logistic regression with a design matrix $\Xi \in \mathbb{R}^{n \times k}$, and assume that the prior of the regression coefficient~$\beta$ is the $\mbox{N}_k(0, I_k)$ distribution, where $I_k$ is a $k \times k$ identity matrix.
	When the response vector is $\gamma \in \mathbb{R}^n$, the posterior density is
	\[
	\dot{\pi}(\beta) \propto \exp\left( - \frac{\|\beta\|^2}{2} \right) \prod_{i=1}^n \frac{\exp[-(1-\gamma_i)\Xi_i^{\top}\beta]}{1+\exp(-\Xi_i^{\top} \beta)},
	\]
	where $\Xi_i^{\top}$ is the $i$'th row of~$\Xi$.
	For this model, the Hessian of $- \log \dot{\pi}(\beta)$ satisfies the conditions of Proposition \ref{pro:hit-concave} with $m = 1$ and $L = \|\Xi^{\top} \Xi /4 + I_k \|$, where $\|\cdot\|$ returns the spectral norm of a matrix \citep[see, e.g.,][]{lee2024fast}.
	Then $L/m$ is no greater than $\xi = \|\Xi^{\top} \Xi /4 + I_k \|$.
	
	{ 
	Consider applying the multiscale Metropolis-within-hit-and-run algorithm with uniform~$\nu$ to the Bayesian logistic model.
	By Proposition \ref{pro:hit-concave}, for $f \in L_0^2(\pi)$,
	\begin{equation} \label{ine:hit-Dirichlet}
		\mathcal{E}_S(f) \geq \frac{C_*}{1 + \lceil \log_2 \xi \rceil} \mathcal{E}_{\bar{S}}(f), \quad \mathcal{E}_S(f) \geq \frac{C_*}{2 (1 + \lceil \log_2 \xi \rceil)} \frac{1}{k \xi} \|f - \pi f\|_{\pi}^2.
	\end{equation}
	which implies
	\[
	\mathrm{Gap}(S) \geq \frac{C_*}{1 + \lceil \log_2 \xi \rceil} \mathrm{Gap}(\bar{S}), \quad \mathrm{Gap}(S) \geq \frac{C_*}{2 (1 + \lceil \log_2 \xi \rceil)} \frac{1}{k \xi}.
	\]
	}

	\subsection{Simulated experiments}

	{ 
		We conduct a simulation experiment to examine the Dirichlet-form comparisons in
		Proposition~\ref{pro:hit-concave}, and test the asymptotic behavior of bounds in \eqref{ine:hit-Dirichlet}.
		A supplementary archive
		contains the Python script, numerical output, configuration file, and figures; it can be found in the GitHub repository \cite{qin2026simulation}.
		
		In the Bayesian logistic regression model, we use the following synthetic data.
		Let $k\geq 8$.  For a positive integer~$r$, take
		$n=16r$, and write
		\[
		a_{\xi,r}=\left[\frac{2(\xi-1)}{r}\right]^{1/2}.
		\]
		Let the design matrix and the response vector be, respectively,
		\[
		\Xi = a_{\xi,r} \left( \begin{array}{ccccc}
			\mathbf{1}_{2r} & \mathbf{0}_{2r} & \hdots & \mathbf{0}_{2r} & \mathbf{0}_{2r \times (k-8)} \\
			\mathbf{0}_{2r} & \mathbf{1}_{2r} & \hdots & \mathbf{0}_{2r} & \mathbf{0}_{2r \times (k-8)}  \\
			\vdots & \vdots & \ddots & \vdots & \vdots \\
			\mathbf{0}_{2r} & \mathbf{0}_{2r} & \hdots & \mathbf{1}_{2r} & \mathbf{0}_{2r \times (k-8)}
		\end{array} \right) \in \mathbb{R}^{n \times k}, \;\; \gamma = \left( \begin{array}{c}
		1 \\ 0 \\ \vdots \\ 1 \\ 0
		\end{array} \right) \in \{0,1\}^{n}.
		\]
		The first
		$8$ coordinates receive likelihood information, whereas the remaining
		coordinates retain their prior-scale curvature.  Moreover,
		\[
		I_k+\frac14\Xi^\top\Xi
		=\operatorname{diag}(\xi I_{8},I_{k-8}).
		\]
		Thus the posterior is $1$-strongly log-concave and $\xi$-log-smooth for every
		value of $n$ and $k$ in the experiment.
		
		We vary $\xi$, $k$, and $n$ separately along the three paths
		\begin{equation} \label{eq:paths}
		\begin{alignedat}{2}
			\xi&\in\{2^4,2^5,\ldots,2^{12}\}
			&&\text{ with $(k,n)=(16,128)$},\\
			k&\in\{2^4,2^5,\ldots,2^{12}\}
			&&\text{ with $(\xi,n)=(8,128)$},\\
			n&\in\{2^4,2^5,\ldots,2^{12}\}
			&&\text{ with $(\xi,k)=(8,16)$}.
		\end{alignedat}
		\end{equation}
		Along the $k$-path, the data-informed component of the target remains fixed
		and the target is extended by additional standard-normal coordinates.  Along
		the $n$-path, $r=n/16$ changes while $a_{\xi,r}$ is rescaled so that $\xi$
		remains fixed.  This separates dimension, sample size, and conditioning.
		For each configuration, we estimate the ratios
		$\mathcal{E}_S(f)/\mathcal{E}_{\bar{S}}(f)$ and
		$\mathcal{E}_S(f)/\|f-\pi f\|_{\pi}^2$ for the two coordinate test functions
		described below.

		We consider the two test functions
		\[
		f_1(\beta)=\beta_1,\quad f_k(\beta)=\beta_k.
		\]
		The design is deliberately made simple so that we can sample from the posterior~$\pi$ exactly.
		Indeed, the posterior factorizes: its first $8$ coordinates have density
		proportional to
		\[
		t\longmapsto\phi(t)\operatorname{sech}^{2r}(a_{\xi,r}t/2),
		\]
		while its remaining coordinates are independent $\mathrm{N}(0,1^2)$ variables.
		Thus, one exact method is to propose each of the first $8$ coordinates from
		$\mathrm{N}(0,1^2)$ and accept a proposal~$t$ with probability
		$\operatorname{sech}^{2r}(a_{\xi,r}t/2)$.
		Consequently, for $f$ equal to either $f_1$ or $f_k$, we may estimate
		$\mathcal{E}_S(f)$, $\mathcal{E}_{\bar{S}}(f)$, and
		$\|f-\pi f\|_{\pi}^2$ using classical Monte Carlo.
		To approximate the Dirichlet form ratio $\mathcal{E}_S(f)/\mathcal{E}_{\bar{S}}(f)$, we generate $N = 100,000$ random vectors from~$\pi$, and call them $\beta^{(1)}, \dots \beta^{(N)}$.
		Then, for each $\beta^{(i)}$, generate $\beta^{(i)}_S$ from $S(\beta^{(i)},\cdot)$, and independently generate $\beta^{(i)}_{\bar{S}}$ from $\bar{S}(\beta^{(i)},\cdot)$.
		For the ideal update $\bar S$, the unique mode of the one-dimensional conditional line density is located using 60 bisection iterations.  If $u_{\mathrm{mode}}$ denotes the resulting mode, we propose $u'$ from $\mathrm{N}(u_{\mathrm{mode}},1^2)$ and accept it with probability
		\[
		\exp\left[\log\dot{\pi}(\theta+u'w)-\log\dot{\pi}(\theta+u_{\mathrm{mode}}w)+\frac{(u'-u_{\mathrm{mode}})^2}{2}\right],
		\]
		where $\theta=\theta(\beta^{(i)},w)$.
		The proposal is a valid rejection envelope because the line density is $1$-strongly log-concave owing to the Gaussian prior.  Thus, apart from floating-point evaluation and the fixed 60-step mode calculation, this produces an exact draw from the ideal conditional distribution.
		
		
		Then the ratio of Dirichlet forms is estimated via the natural Monte Carlo estimator 
		\[
		\frac{(2N)^{-1} \sum_{i=1}^N [f(\beta^{(i)}) - f(\beta^{(i)}_S)]^2}{(2N)^{-1} \sum_{i=1}^N [f(\beta^{(i)}) - f(\beta^{(i)}_{\bar{S}})]^2}.
		\]
		To approximate the Rayleigh quotient $\mathcal{E}_S(f)/\|f - \pi f\|_{\pi}^2$, we replace the denominator with $N^{-1} \sum_{i=1}^N f(\beta^{(i)})^2$ --- this is reasonable since the expectation $\pi f$ happens to vanish for the current design and test functions.

		Pointwise \(95\%\) error bars are constructed for each estimated Dirichlet ratio.
		For instance, to obtain a 95\% confidence interval for $\mathcal{E}_S(f)/\mathcal{E}_{\bar{S}}(f)$, we compute the asymptotic variance of
		\[
		\log \frac{(2N)^{-1} \sum_{i=1}^N [f(\beta^{(i)}) - f(\beta^{(i)}_S)]^2}{(2N)^{-1} \sum_{i=1}^N [f(\beta^{(i)}) - f(\beta^{(i)}_{\bar{S}})]^2}.
		\]
		Per the delta method, the variance of this quantity can be estimated via $N^{-1}$ times the sample variance of
		\[
		\frac{[f(\beta^{(i)}) - f(\beta^{(i)}_S)]^2}{\sum_{j=1}^N N^{-1} [ f(\beta^{(j)}) - f(\beta^{(j)}_S)]^2} - \frac{[f(\beta^{(i)}) - f(\beta^{(i)}_{\bar{S}})]^2}{\sum_{j=1}^N N^{-1} [f(\beta^{(j)}) - f(\beta^{(j)}_{\bar{S}})]^2}, \quad i \in \{1,\dots,N\}.
		\]

		\begin{figure}[t]
			\centering
			\includegraphics[width=0.32\textwidth]
			{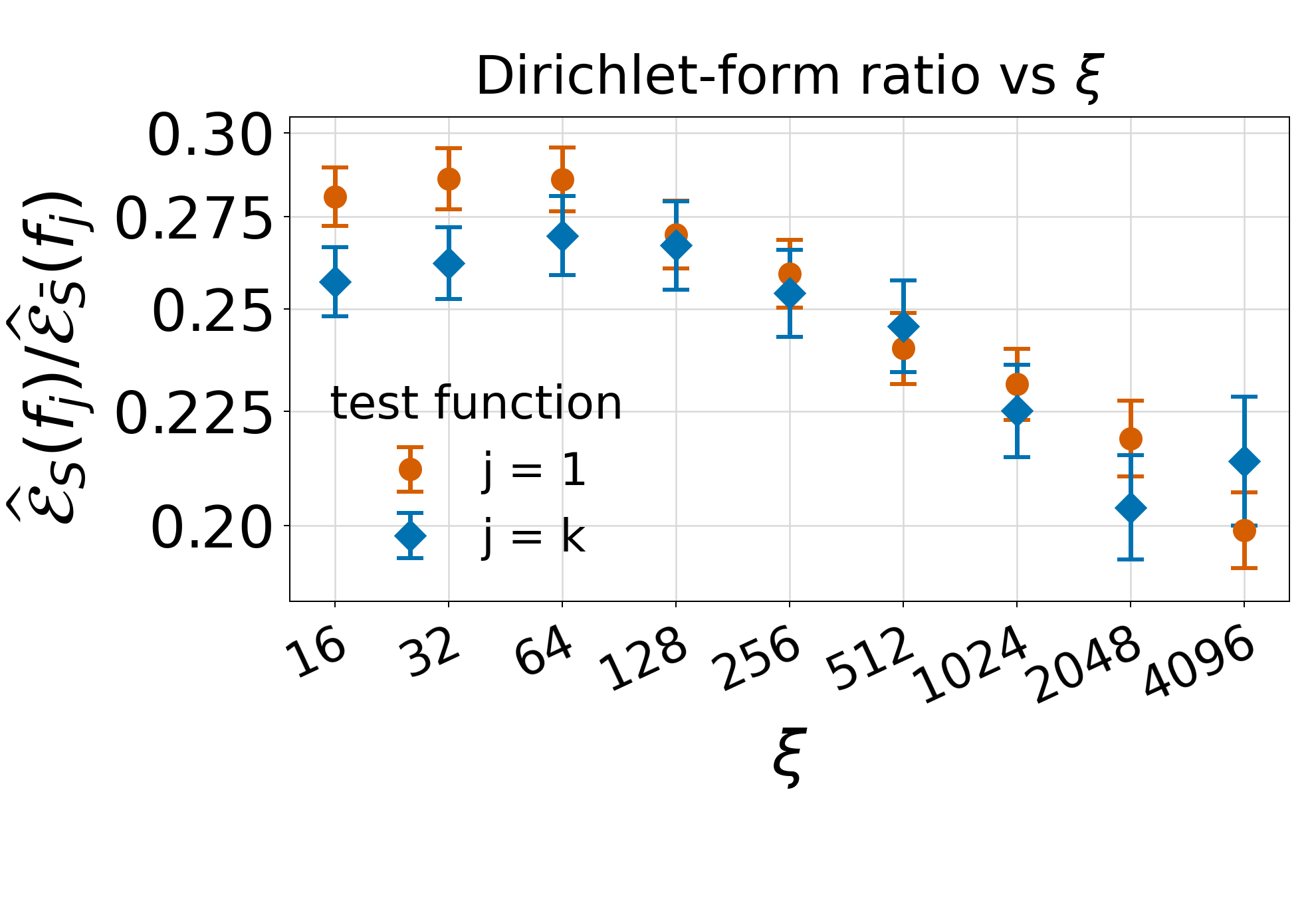}
			\includegraphics[width=0.32\textwidth]
			{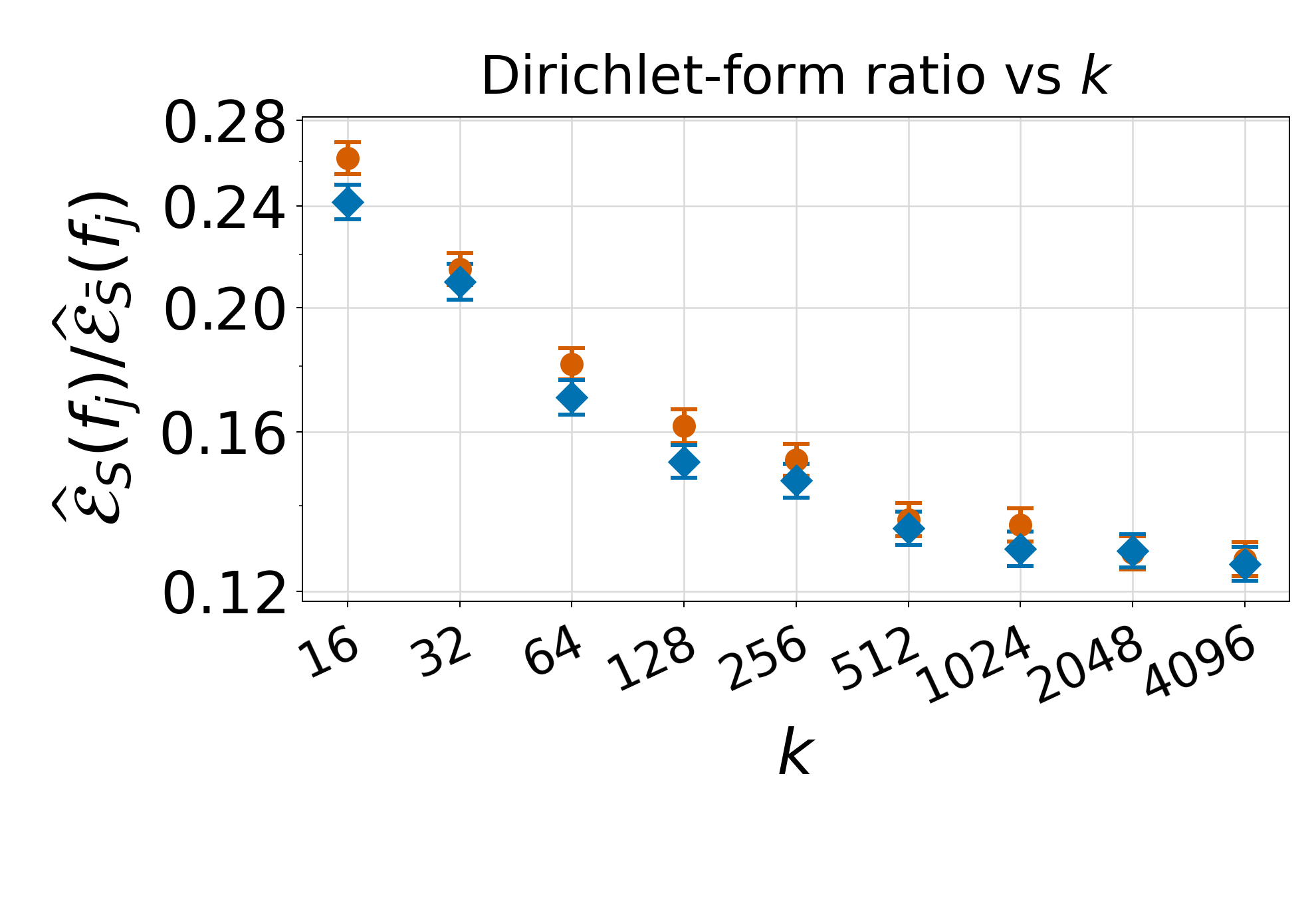}
			\includegraphics[width=0.32\textwidth]
			{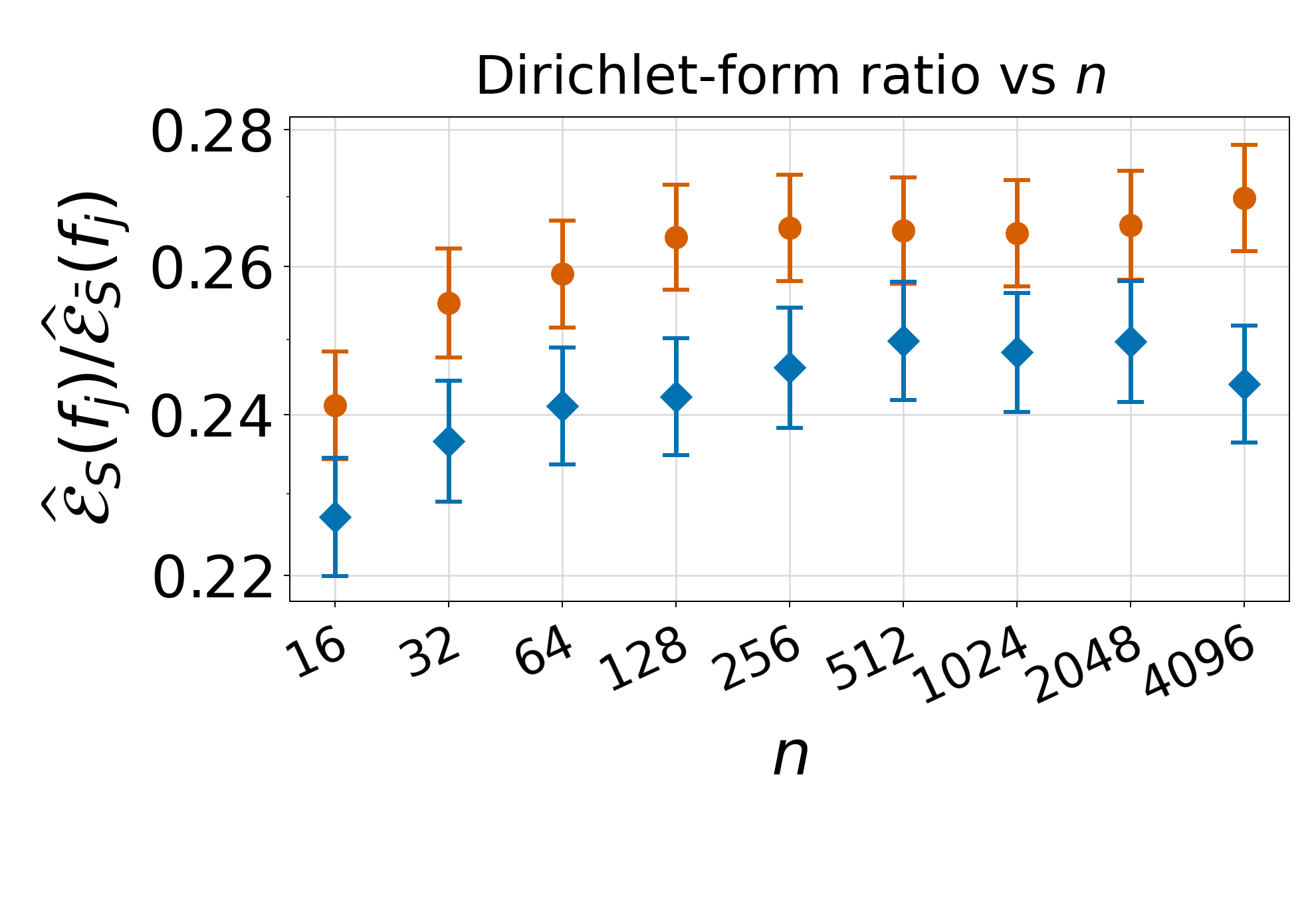}
			\caption{  Estimated Dirichlet-form ratios
				\(\mathcal E_S(f_j)/\mathcal E_{\bar S}(f_j)\) for \(j=1,k\) along the paths \eqref{eq:paths}.
				Error bars are pointwise 95\% paired
				log-delta Monte Carlo intervals.}
			\label{fig:dirichlet-S-barS}
		\end{figure}
		
		\begin{figure}[t]
			\centering
			\includegraphics[width=0.32\textwidth]
			{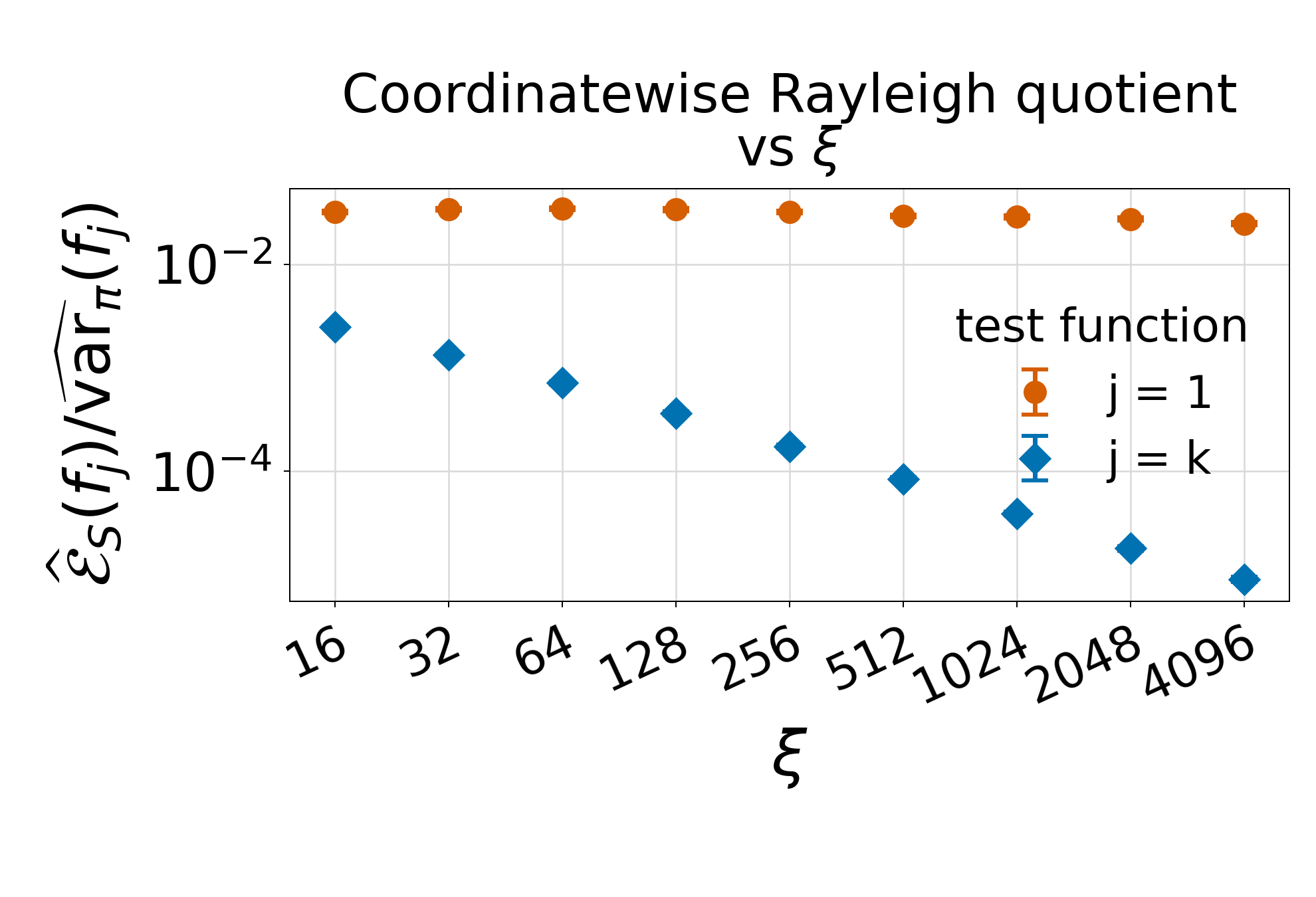}
			\includegraphics[width=0.32\textwidth]
			{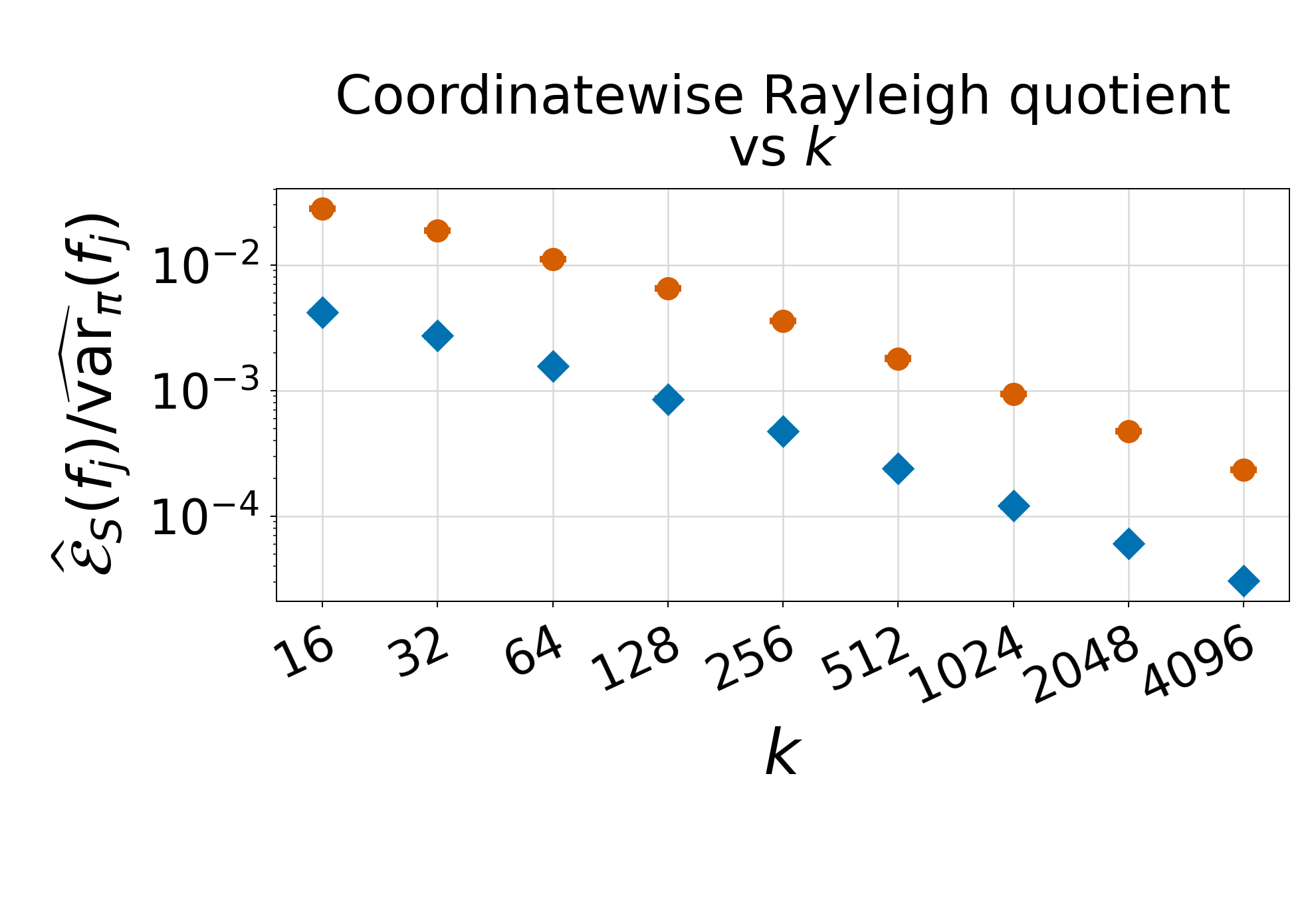}
			\includegraphics[width=0.32\textwidth]
			{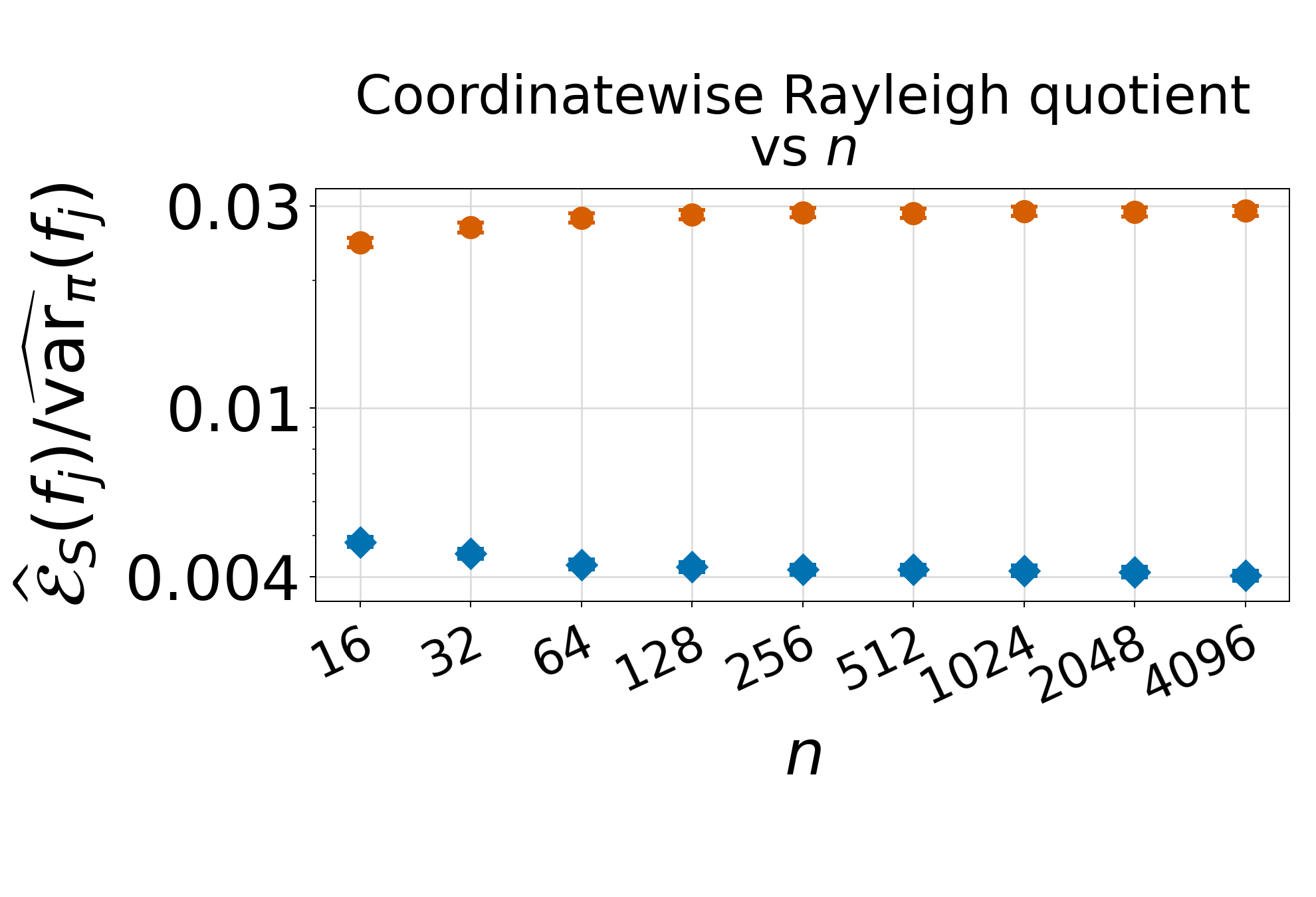}
			\caption{  Estimated coordinatewise Rayleigh quotients
				\(\mathcal E_S(f_j)/\| f_j - \pi f_j \|_{\pi}^2\) for \(j=1,k\) along the paths \eqref{eq:paths}.
				 Error bars are pointwise 95\% paired
				log-delta Monte Carlo intervals.}
			\label{fig:dirichlet-S-pi}
		\end{figure}
		
		The simulation results are plotted in Figures \ref{fig:dirichlet-S-barS} and \ref{fig:dirichlet-S-pi}.
		The axes of the plots are in log-scale.
		
		By \eqref{ine:hit-Dirichlet}, it holds that, for $f \in L_0^2(\pi)$,
		\begin{equation} \label{eq:hit-order}
			\log \frac{\mathcal{E}_S(f)}{\mathcal{E}_{\bar{S}}(f)} \gtrsim -\log (1 + \lceil \log_2 \xi \rceil), \quad \log \frac{\mathcal{E}_S(f)}{\|f-\pi f\|_{\pi}^2} \gtrsim -\log k - \log \xi - \log(1 + \lceil\log_2 \xi \rceil),
		\end{equation}
		where $\gtrsim$ means that the left-hand side is greater than the right-hand
		side plus a constant for large values of $\xi$ and $k$.
		
		Figure~\ref{fig:dirichlet-S-barS} gives the point estimates of
		$\mathcal{E}_S(f_j)/\mathcal{E}_{\bar S}(f_j)$.
		A slow deterioration with $\xi$ is visible.
		Regressing $\log[\mathcal{E}_S(f_j)/\mathcal{E}_{\bar S}(f_j)]$ against $\log \xi$ yields powers (slopes) $-0.066$ ($R^2 = 0.919$) and $-0.047$ ($R^2 = 0.746$) for $j = 1$ and $j = k$, respectively.
		Deterioration with $k$ appears to stop after a certain point.
		No deterioration
		with $n$ is visible after $k$ and $\xi$ are held fixed.
		
		Figure~\ref{fig:dirichlet-S-pi} gives a more direct view of the dimensional and
		conditioning factors in the second relation of~\eqref{eq:hit-order}.  Regressing
		$\log [\mathcal{E}_S(f_j)/\|f_j-\pi f_j\|_{\pi}^2]$ against $\log\xi$ yields
		powers $-0.054$ ($R^2=0.788$) and $-1.025$ ($R^2=0.999$) for $j=1$ and $j=k$,
		respectively.  The corresponding powers along the $k$-path are $-0.876$
		($R^2=0.994$) and $-0.901$ ($R^2=0.996$).  These fitted powers are consistent with the worst-case $\xi^{-1}$
		and $k^{-1}$ factors in~\eqref{ine:hit-Dirichlet}.  Finally, along the $n$-path, no significant deterioration is observed.
		The fitted powers in this case are $0.024$ ($R^2 = 0.664$) and $-0.027$ ($R^2 = 0.775$) for $j = 1$ and $j = k$, respectively.
		
		Overall, the experiment supports the asymptotic relations in
		\eqref{ine:hit-Dirichlet}.  We emphasize that the bounds in
		\eqref{ine:hit-Dirichlet} cover the worst-case scenario over test functions.
		It is therefore unsurprising that, for some test functions, the observed dependence on $(\xi,k,n)$ is milder than the bounds allow.
		The experiment also provides no evidence that the universal constant $C_*$ must be as small as its theoretical value suggests.
		Indeed, \eqref{ine:hit-Dirichlet} is best interpreted as a worst-case scaling guarantee rather than as a numerically sharp bound over the finite parameter ranges considered here. For the finite parameter values considered here, the resulting lower bounds are much smaller than the observed Dirichlet-form ratios.
	}

	{ 

	\section{Hybrid Data Augmentation for Well-Conditioned Distributions} \label{sec:da-example}
	
	As a second example, we study the convergence rate of hybrid data augmentation algorithms with two intractable conditional distributions, as described in Section~\ref{ssec:doubly-0}, when the target distribution~$\varphi$ is well-conditioned.
	
	For simplicity, we first take $k_1=k_2=1$, so that $(x_1,x_2)\in\mathbb{R}^2$.
	The multi-dimensional case is studied in Remark~\ref{rem:da-multi}.
	Recall that the target distribution $\varphi(\df(x_1,x_2))$ has disintegrations
	\[
	\varphi_1(\df x_1)\,\varphi_{2,x_1}(\df x_2)
	=
	\varphi_2(\df x_2)\,\varphi_{1,x_2}(\df x_1).
	\]
	We assume that these distributions admit density functions, denoted by $\dot{\varphi}$, with the appropriate subscripts for the marginal and conditional densities.
	
	Assume that the two conditional distributions are both difficult to sample from exactly.
	Instead, let $H_{i,x_{-i}}$ be reversible with respect to $\varphi_{i,x_{-i}}$, $i\in\{1,2\}$, where $x_{-1}=x_2$ and $x_{-2}=x_1$.
	The hybrid data augmentation algorithm of interest has Mtk
	\[
	T((x_1,x_2),\df(x_1',x_2'))
	=H_{1,x_2}(x_1,\df x_1')\,H_{2,x_1'}(x_2,\df x_2').
	\]

	We use the multi-scale Gaussian random-walk Metropolis-Hastings algorithm described in Section~\ref{sec:example-hybrid} to sample from the conditional distributions.
	Let $J$ be a positive integer, let $\sigma_1,\dots,\sigma_J$ be positive numbers, and let $(p_1,\dots,p_J)$ be a probability vector.
	For each $(i,x_{-i})$, take
	\[
	H_{i,x_{-i}}(x_i,\df x_i')
	=\sum_{j=1}^J p_j H_{i,x_{-i},j}(x_i,\df x_i'),
	\]
	where
	\[
	\begin{aligned}
		H_{i,x_{-i},j}(x_i,\df x_i')
		={}&\frac{1}{\sigma_j}\phi\left(\frac{x_i'-x_i}{\sigma_j}\right)
		a_{i,x_{-i}}(x_i,x_i')\,\df x_i' \\
		&+\left[1-\int_{\mathbb{R}}\frac{1}{\sigma_j}
		\phi\left(\frac{x_i''-x_i}{\sigma_j}\right)
		a_{i,x_{-i}}(x_i,x_i'')\,\df x_i''\right]\delta_{x_i}(\df x_i'),
	\end{aligned}
	\]
	$\phi$ is the density of the $\mbox{N}(0,1^2)$ distribution, and
	\[
	a_{i,x_{-i}}(x_i,x_i')
	=\min\left\{1,
	\frac{\dot{\varphi}_{i,x_{-i}}(x_i')}
	{\dot{\varphi}_{i,x_{-i}}(x_i)}\right\}.
	\]
	
	Recall that $T$ is in general non-reversible, making it challenging to analyze.
	We use Proposition~\ref{pro:doubly} to bound $\|T\|_{\varphi}$.
	
	\begin{proposition}\label{pro:da-concave}
		Assume that $(x_1,x_2)\mapsto-\log\dot{\varphi}(x_1,x_2)$ is twice differentiable, and that the eigenvalues of its Hessian are bounded between two positive numbers $m$ and $L$, where $m\leq L$.
		In the multi-scale Metropolis-Hastings algorithms targeting $\varphi_{i,x_{-i}}$, let
		\[
		J=1+\lceil\log_2(L/m)\rceil,\qquad
		p_j=\frac{1}{J},\qquad
		\sigma_j^2=\frac{2^{j-1}}{L},\quad j=1,\dots,J.
		\]
		Then
		\[
		1-\|T\|_{\varphi}^2
		\geq
		\left(\frac{C_*}{1+\lceil\log_2(L/m)\rceil}\right)^2\frac{m}{L},
		\]
		where $C_*\approx10^{-14}$ is the universal constant in Proposition~\ref{pro:hit-concave}.
	\end{proposition}
	
	\begin{proof}
		Fix $i\in\{1,2\}$ and $x_{-i}\in\mathbb{R}$.
		The negative log-density of $\varphi_{i,x_{-i}}$ is twice differentiable, and its second derivative lies in $[m,L]$.
		By the Brascamp--Lieb inequality and the Cramer--Rao lower bound, the variance $v_{i,x_{-i}}$ of $\varphi_{i,x_{-i}}$ lies in $[1/L,1/m]$.
		Thus, for some $j\in\{1,\dots,J\}$,
		\[
		v_{i,x_{-i}}\leq\sigma_j^2\leq2v_{i,x_{-i}}.
		\]
		By Theorem~1 of \cite{secchi2025spectral} and Cheeger's inequality,
		\[
		\mathrm{Gap}(H_{i,x_{-i},j})\geq C_*.
		\]
		It follows that
		\[
		\begin{aligned}
			\mathcal{E}_{H_{i,x_{-i}}}(f)
			&=\frac{1}{J}\sum_{\ell=1}^J
			\mathcal{E}_{H_{i,x_{-i},\ell}}(f) \\
			&\geq\frac{1}{J}\mathcal{E}_{H_{i,x_{-i},j}}(f),
		\end{aligned}
		\]
		and hence
		\[
		\mathrm{Gap}(H_{i,x_{-i}})
		\geq\frac{C_*}{J}
		=\frac{C_*}{1+\lceil\log_2(L/m)\rceil}.
		\]
		
		By Lemma~3.1 of \cite{baxendale2005renewal}, each $H_{i,x_{-i},j}$ is positive semi-definite.
		Therefore, $H_{i,x_{-i}}$ is positive semi-definite and
		\[
		1-\|H_{i,x_{-i}}\|_{\varphi_{i,x_{-i}}}^2
		\geq \mathrm{Gap}(H_{i,x_{-i}})
		\geq\frac{C_*}{J}.
		\]
		
		Let
		\[
		\bar{S}(x_2,\df x_2')
		=\int_{\mathbb{R}}\varphi_{1,x_2}(\df x_1)\,
		\varphi_{2,x_1}(\df x_2')
		\]
		be the ideal data augmentation Mtk on the second coordinate.
		By Corollary~3.7 of \cite{ascolani2024entropy},
		\[
		\mathrm{Gap}(\bar{S})\geq\frac{m}{L}.
		\]
		Applying Proposition~\ref{pro:doubly}, its positive semi-definiteness assertion, and Proposition~\ref{pro:andrieu}, we obtain
		\[
		\begin{aligned}
			1-\|T\|_{\varphi}^2
			&\geq
			\left(1-\esssup_{x_2}
			\|H_{1,x_2}\|_{\varphi_{1,x_2}}^2\right)
			\left[\essinf_{x_1}\mathrm{Gap}(H_{2,x_1})\right]
			\mathrm{Gap}(\bar{S}) \\
			&\geq
			\left(\frac{C_*}{J}\right)^2\frac{m}{L}.
		\end{aligned}
		\]
		This proves the result.
	\end{proof}
	
	As a simple but concrete example, suppose that
	\[
	\dot{\varphi}(x_1,x_2)
	\propto
	\exp\left[-\frac{x_1^2+x_2^2}{2}
	-a\log\cosh(x_1-x_2)\right],
	\qquad (x_1,x_2)\in\mathbb{R}^2,
	\]
	where $a\geq0$.
	The full conditional densities are
	\[
	\begin{aligned}
		\dot{\varphi}_{1,x_2}(x_1)
		&\propto \exp\left(-\frac{x_1^2}{2}\right)
		\cosh(x_1-x_2)^{-a}, \\
		\dot{\varphi}_{2,x_1}(x_2)
		&\propto \exp\left(-\frac{x_2^2}{2}\right)
		\cosh(x_1-x_2)^{-a},
	\end{aligned}
	\]
	which are both non-standard.
	
	The Hessian matrix of $-\log\dot{\varphi}(x_1,x_2)$ is
	\[
	I_2+\frac{4a}{[\exp(x_1-x_2)+\exp(x_2-x_1)]^2}
	\left(\begin{array}{cc}
		1&-1\\
		-1&1
	\end{array}\right).
	\]
	Its eigenvalues lie in $[m,L]$, where $m=1$ and $L=1+2a$.
	Thus, in each multi-scale Metropolis-Hastings algorithm, we take
	\[
	J=1+\lceil\log_2(1+2a)\rceil,\qquad
	p_j=\frac{1}{J},\qquad
	\sigma_j^2=\frac{2^{j-1}}{1+2a},\quad j=1,\dots,J.
	\]
	Proposition~\ref{pro:da-concave} gives
	\[
	1-\|T\|_{\varphi}^2
	\geq
	\left(\frac{C_*}{1+\lceil\log_2(1+2a)\rceil}\right)^2
	\frac{1}{1+2a}.
	\]
	Consequently, by \eqref{ine:mixing}, for $\epsilon>0$ and an initial distribution~$\mu$,
	\[
	\begin{aligned}
		t_*(T,\mu,\epsilon)
		&\leq
		\frac{2\{\log\|\df\mu/\df\varphi-1\|_{\varphi}
			-\log\epsilon\}}{1-\|T\|_{\varphi}^2} \\
		&\leq
		\frac{2(1+2a)[1+\lceil\log_2(1+2a)\rceil]^2}
		{C_*^2}
		\left(\log\|\df\mu/\df\varphi-1\|_{\varphi}
		-\log\epsilon\right).
	\end{aligned}
	\]
	Hence, up to a squared logarithmic factor, this bound scales linearly with~$a$.
	
	\begin{remark} \label{rem:da-multi}
		The restriction $k_1=k_2=1$ allows us to use the one-dimensional, variance-adaptive result of \cite{secchi2025spectral}.
		A multivariate extension can instead be obtained from the random-walk Metropolis bounds of \cite{andrieu2024explicit}.
		Suppose that $(x_1,x_2)\in\mathbb{R}^{k_1}\times\mathbb{R}^{k_2}$ and retain the Hessian assumption in Proposition~\ref{pro:da-concave}.
		For the update of coordinate block~$i$, use a random-walk Metropolis kernel with proposal increments distributed as
		\[
		\mbox{N}_{k_i}\left(0,\frac{1}{2Lk_i}I_{k_i}\right).
		\]
		Theorem~1 of \cite{andrieu2024explicit} gives, uniformly in $x_{-i}$,
		\[
		\mathrm{Gap}(H_{i,x_{-i}})
		\geq C_{\mathrm{RWM}}\frac{m}{L}\frac{1}{k_i},
		\qquad
		C_{\mathrm{RWM}}=\frac{1.972\times10^{-4}}{2e}.
		\]
		Since Gaussian random-walk Metropolis kernels are positive semi-definite, the preceding argument yields
		\[
		1-\|T\|_{\varphi}^2
		\geq
		\frac{C_{\mathrm{RWM}}^2}{k_1k_2}
		\left(\frac{m}{L}\right)^3.
		\]
		Thus, the framework extends directly to multivariate blocks, although this bound has an additional $1/(k_1k_2)$ factor and a less favorable dependence on $m/L$ than the scalar multi-scale bound.
	\end{remark}
	}

	\bigskip
	
    \noindent{\bf Acknowledgment:}
	The author was partially supported by NSF.
	The author thanks Guanyang Wang for insightful discussions, especially regarding the analysis of hybrid hit-and-run samplers.

	\vspace{1.5cm}
	
	\appendix

	{\noindent \bf \Large Appendix}

	\section{Technical Lemmas} \label{app:lemma}

	\begin{lemma} \label{lem:iso-1}
		Let $(\Omega, \mathcal{F}, \rho)$ be a probability space.
		Let $K: \Omega \times \mathcal{F} \to [0,1]$ be an Mtk.
		Assume that there is a set $\Omega_0 \in \mathcal{F}$ such that $\rho(\Omega_0) = 1$.
		Let $\mathcal{F}_0$ be $\mathcal{F}$ restricted to $\Omega_0$.
		Let $\rho_0: \mathcal{F}_0 \to [0,1]$ be $\rho$ restricted to $\mathcal{F}_0$.
		Let $K_0: \Omega_0 \times \mathcal{F}_0 \to [0,1]$ be an Mtk satisfying $K_0(w, \mathsf{F}) = K(w, \mathsf{F})$ for $(w, \mathsf{F}) \in \Omega_0 \times \mathcal{F}_0$.
		Assume that $\rho_0 K_0 = \rho_0$.
		Then $\rho K = \rho$, and $\|K\|_{\rho} = \|K_0\|_{\rho_0}$.
		If, furthermore, $K_0$ is reversible with respect to $\rho_0$, then $K$ is reversible with respect to $\rho$, and $\mathrm{Gap}(K) = \mathrm{Gap}(K_0)$.
	\end{lemma}
	
	\begin{proof}
		It is easy to see that $\rho K = \rho$, and if $K_0$ is reversible with respect to $\rho_0$, then $K$ is reversible with respect to~$\rho$.
		We will focus on showing that the norms and spectral gaps of the two operators are equal.
		
		Define the map $U: L_0^2(\rho) \to L_0^2(\rho_0)$ as follows: $U f(w) = f(w)$ for $w \in \Omega_0$.
		Then $U$ is invertible, and $U^{-1} g(w) = g(w) \, \ind_{\Omega_0}(w)$ for $g \in L_0^2(\rho_0)$ and $w \in \Omega$.
		(Note that two functions in $L_0^2(\rho)$ are equal as long as they coincide $\rho_0$-a.e. on $\Omega_0$.)
		In fact, $U$ is unitary.
		Note that $K = U^{-1} K_0 U$, so $K$ and $K_0$ are unitarily equivalent.
		It is then easy to see that $K$ and $K_0$ are equal in norm and spectral gap.
	\end{proof}
	
	The next lemma is an easy corollary of Lemma~\ref{lem:iso-1}.
	\begin{lemma} \label{lem:iso-2}
		Let $(\Omega, \mathcal{F}, \rho)$ be a probability space.
		Let $K: \Omega \times \mathcal{F} \to [0,1]$ be an Mtk.
		Assume that there is a set $\Omega_0 \in \mathcal{F}$ such that $\rho(\Omega_0) = 1$.
		Let $\mathcal{F}_0$ be $\mathcal{F}$ restricted to $\Omega_0$.
		Let $\rho_0: \mathcal{F}_0 \to [0,1]$ be $\rho$ restricted to $\mathcal{F}_0$.
		Let $(\Omega_1, \mathcal{F}_1, \rho_1)$ be another probability space, and let $K_1: \Omega_1 \times \mathcal{F}_1 \to [0,1]$ be an Mtk such that $\rho_1 K_1 = \rho_1$.
		Assume that there is a one-to-one measurable function $f: \Omega_1 \to \Omega_0$ such that $\rho_0 = \rho_1 \circ f^{-1}$ and that $K(w, \mathsf{A}) = K_1(f^{-1}(w), f^{-1}( \mathsf{A} ))$ for $(w, \mathsf{A}) \in \Omega_0 \times \mathcal{F}_0$.
		Then $\rho K = \rho$, and $\|K\|_{\rho} = \|K_1\|_{\rho_1}$.
		If, furthermore, $K_1$ is reversible with respect to $\rho_1$, then $K$ is reversible with respect to~$\rho$, and $\mathrm{Gap}(K) = \mathrm{Gap}(K_1)$.
		
	\end{lemma}

	\section{Technical Proofs} \label{app:proofs}
	
%

	\subsection{Proof of Lemma \ref{lem:peskun}} \label{app:peskun}
	
	\begin{proof}
		Let $f \in L_0^2(\pi)$ be arbitrary.
		Then
		\[
		\mathcal{E}_S(f) = \|f\|_{\pi}^2 - \|Pf\|_{\tilde{\pi}}^2 + \|Pf\|_{\tilde{\pi}}^2 - \langle f, P^*QPf \rangle_{\pi} = \|f\|_{\pi}^2 - \|Pf\|_{\tilde{\pi}}^2 + \mathcal{E}_Q(Pf), 
		\]
		so, by \ref{H2} and the assumption that $c_0 \in (0,1]$,
		\[
		\mathcal{E}_S(f) \geq c_0 \|f\|_{\pi}^2 - c_0 \|Pf\|_{\tilde{\pi}}^2 + \mathcal{E}_Q(Pf).
		\]
		Similarly,
		\[
		\mathcal{E}_{P^*\hat{Q}P}(f) = \|f\|_{\pi}^2 - \|Pf\|_{\tilde{\pi}}^2 + \mathcal{E}_{\hat{Q}}(Pf).
		\]
		To prove the lemma, it suffices to show that $\mathcal{E}_Q(Pf) \geq c_0 \, \mathcal{E}_{\hat{Q}}(Pf)$.
		{ In the first inequality below, the integrand is nonnegative, so we retain only transitions with $z'=z$ and $y'\neq y$.  The transitions with $z'\neq z$ are simply discarded for this lower bound, while the transitions with $y'=y$ make zero contribution because the squared increment vanishes.}
		By \ref{H4},
		\[
		\begin{aligned}
			\mathcal{E}_Q(Pf) &= \frac{1}{2} \int_{\Y \times \Z} \tilde{\pi}(\df (y,z)) \int_{\Y \times \Z} Q((y,z), \df (y',z')) \, [Pf(y',z') - Pf(y,z)]^2 \\
			&\geq \frac{1}{2} \int_{\Y \times \Z} \tilde{\pi}(\df (y,z)) \int_{\{y\}^c \times \{z\}} Q((y,z), \df (y',z')) \, [Pf(y',z') - Pf(y,z)]^2 \\
			&\geq \frac{c_0}{2} \int_{\Y \times \Z} \tilde{\pi}(\df (y,z)) \int_{\{y\}^c} Q_z(y, \df y') \, [Pf(y',z) - Pf(y,z)]^2 \\
			&= \frac{c_0}{2} \int_{\Y \times \Z} \tilde{\pi}(\df (y,z)) \int_{\Y} Q_z(y, \df y') \, [Pf(y',z) - Pf(y,z)]^2 \\
			&= c_0 \, \mathcal{E}_{\hat{Q}}(Pf).
		\end{aligned}
		\]
		Thus, the lemma holds.
	\end{proof}
	
	\subsection{Proof of Proposition \ref{pro:iteration-1}} \label{app:iteration-1}
	
	\begin{proof}
			Fix $z\in\Z$ for which the assumptions hold. Since
			$Q_z:L_0^2(\varpi_z)\to L_0^2(\varpi_z)$ is self-adjoint, the
			spectrum of $Q_z^{\lambda(z)}$ is obtained by raising the elements
			of the spectrum of $Q_z$ to the power $\lambda(z)$
			\citep[see, e.g.,][\S 31, Theorem~2]{helmberg2014introduction}.
			Moreover,
			\[
			\sup\operatorname{Spec}(Q_z)
			=1-\mathrm{Gap}(Q_z).
			\]
			
			If $\lambda(z)$ is odd, then the function
			$t\mapsto t^{\lambda(z)}$ is increasing on $[-1,1]$.
			Note that $1 - \mathrm{Gap}(Q_z) \in [-1,1]$.
			Then it holds that, for $\lambda(z)$ odd,
			\[
			1-\mathrm{Gap}\big(Q_z^{\lambda(z)}\big) = \sup \operatorname{Spec}(Q_z^{\lambda(z)})
			=
			[1-\mathrm{Gap}(Q_z)]^{\lambda(z)}
			\leq
			\big[1-\min\{1,\mathrm{Gap}(Q_z)\} \big]^{\lambda(z)}.
			\]
			
			Suppose instead that
			$\mathrm{Gap}_-(Q_z)\geq\mathrm{Gap}(Q_z)$. Writing
			\[
			a_z=\inf\operatorname{Spec}(Q_z) = -1 + \mathrm{Gap}_-(Q_z),
			\qquad
			b_z=\sup\operatorname{Spec}(Q_z)
			=1-\mathrm{Gap}(Q_z),
			\]
			the assumed inequality gives
			\[
			1+a_z\geq1-b_z,
			\]
			and therefore $a_z\geq-b_z$. Since $a_z\leq b_z$, this also implies
			$b_z\geq0$ and
			\[
			\|Q_z\|_{\varpi_z}
			=\max\{|a_z|,|b_z|\}
			=b_z.
			\]
			Consequently,
			\[
			\begin{aligned}
				1-\mathrm{Gap}\big(Q_z^{\lambda(z)}\big)
				&\leq \big\|Q_z^{\lambda(z)}\big\|_{\varpi_z} \\
				&=\|Q_z\|_{\varpi_z}^{\lambda(z)} \\
				&=[1-\mathrm{Gap}(Q_z)]^{\lambda(z)} \\
				&= [1-\min\{1, \mathrm{Gap}(Q_z)\}]^{\lambda(z)}
			\end{aligned}
			\]
			
			Thus, in either case, by the assumption on $\lambda(z)$,
			\[
			1-\mathrm{Gap}\big(Q_z^{\lambda(z)}\big)
			\leq
			\big[1-\min\{1,\mathrm{Gap}(Q_z)\}\big]^{\lambda(z)}
			\leq1-\kappa_{\dagger}.
			\]
			It follows that
			\[
			\mathrm{Gap}\big(Q_z^{\lambda(z)}\big)
			\geq\kappa_{\dagger}.
			\]
			The desired conclusions now follow from
			Theorem~\ref{thm:uniform}.
	\end{proof}
	
	\subsection{Proof of Lemma~\ref{lem:weak-new}}
	
	\begin{proof}
		By Lemmas \ref{lem:peskun} and \ref{lem:main} and the weak Poincar\'{e} inequality for $Q_z$, for $f \in L_0^2(\pi)$ and $u > 0$,
		\[
		\begin{aligned}
			c_0^{-1} \mathcal{E}_S(f) &\geq \|f\|_{\pi}^2 - \|P f \|_{\tilde{\pi}}^2 + \int_{\Z} \varpi(\df z) \, \mathcal{E}_{Q_z}(Pf_z - \varpi_z Pf_z) \\
			& \geq \|f\|_{\pi}^2 - \|P f \|_{\tilde{\pi}}^2 + u^{-1}  \int_{\Z} \varpi(\df z) \, \left[ \|Pf_z - \varpi_z Pf_z\|_{\varpi_z}^2 -  \alpha(z,u) \, \|Pf_z - \varpi_z Pf_z\|_{\scriptsize\mathrm{osc}}^2 \right] \\
			& = \|f\|_{\pi}^2 - \|P f \|_{\tilde{\pi}}^2 + u^{-1}  \int_{\Z} \varpi(\df z) \, \left[ \|Pf_z - \varpi_z Pf_z\|_{\varpi_z}^2 -  \alpha(z,u) \, \|Pf_z\|_{\scriptsize\mathrm{osc}}^2 \right] ,
		\end{aligned}
		\]
		where $Pf_z(y) = Pf(y,z)$ for $(y,z) \in \Y \times \Z$.
		It is not difficult to see that
		\[
		\int_{\Z} \varpi(\df z) \, \|Pf_z - \varpi_z Pf_z\|_{\varpi_z}^2 = \|Pf\|_{\tilde{\pi}}^2 - \langle P^* \Proj^* \Proj P f, f \rangle_{\pi},
		\]
		so
		\[
		\begin{aligned}
			c_0^{-1} \mathcal{E}_S(f) \geq& \|f\|_{\pi}^2 - \|P f \|_{\tilde{\pi}}^2 + u^{-1} \left( \|Pf\|_{\tilde{\pi}}^2 - \langle P^* \Proj^* \Proj P f, f \rangle_{\pi} \right) \\ &- u^{-1} \int_{\Z} \varpi(\df z) \,  \alpha(z,u) \, \|Pf_z\|_{\scriptsize\mathrm{osc}}^2.
		\end{aligned}
		\]
		{  On the other hand, for any measurable $g: \Y \times \Z \to \mathbb{R}$ and $\varpi$-a.e. $z \in \Z$, $\essinf_y g(y,z) \geq \essinf_{y,z'} g(y,z')$, $\esssup_y g(y,z) \leq \esssup_{y,z'} g(y,z')$,
		so letting $g = Pf$ gives $\|Pf_z\|_{\scriptsize\mathrm{osc}} \leq \|Pf\|_{\scriptsize\mathrm{osc}}$ for $\varpi$-a.e. $z \in \Z$.
		} 
		Moreover, by assumption, $\|Pf\|_{\scriptsize\mathrm{osc}} \leq \|f\|_{\scriptsize\mathrm{osc}}$.
		It follows that 
		\[
		c_0^{-1} \mathcal{E}_S(f) \geq \|f\|_{\pi}^2 - \|P f \|_{\tilde{\pi}}^2 + u^{-1} [ \|Pf\|_{\tilde{\pi}}^2 - \langle P^* \Proj^* \Proj P f, f \rangle_{\pi} - \bar{\alpha}(u) \|f\|_{\scriptsize\mathrm{osc}}^2 ].
		\]
		Then, by \ref{H2}, for $f \in L_0^2(\pi)$ and $u \geq 1$,
		\[
		\begin{aligned}
			u \, c_0^{-1} \mathcal{E}_S(f) &\geq u \|f\|_{\pi}^2 - (u-1)\|Pf\|_{\tilde{\pi}}^2 - \langle P^* \Proj^* \Proj P f, f \rangle_{\pi} - \bar{\alpha}(u) \|f\|_{\scriptsize\mathrm{osc}}^2 \\
			&\geq \mathcal{E}_{P^* \Proj^* \Proj P}(f) -  \bar{\alpha}(u) \|f\|_{\scriptsize\mathrm{osc}}^2.
		\end{aligned}
		\]
		This establishes the desired result.
	\end{proof}
	
	\subsection{Proof of Proposition \ref{pro:weak}}
		
	\begin{proof}
		By Popoviciu's inequality \citep{popoviciu1935equations}, when $u < 1$, for $f \in L_0^2(\pi)$,
		\[
		\|f\|_{\pi}^2 \leq \frac{1}{4} \|f\|_{\scriptsize\mathrm{osc}}^2 \leq u \, \mathcal{E}_S(f) + \tilde{\beta}(u) \,  \|f\|_{\scriptsize\mathrm{osc}}^2.
		\]
		Assume that $u \geq 1$.
		Then, by Lemma \ref{lem:weak-new} and the weak Poincar\'{e} inequality for $P^* \Proj^* \Proj P$, for $f \in L_0^2(\pi)$, $u_1 > 0$, and $u_2 > 1$,
		\[
		\begin{aligned}
			\|f\|_{\pi}^2 &\leq u_1 \, \mathcal{E}_{P^* \Proj^* \Proj P}(f) + \beta(u_1) \, \|f\|_{\scriptsize\mathrm{osc}}^2 \\
			&\leq u_1 \, u_2 \, c_0^{-1} \mathcal{E}_S(f) + u_1 \, \bar{\alpha}(u_2) \, \|f\|_{\scriptsize\mathrm{osc}}^2 + \beta(u_1) \, \|f\|_{\scriptsize\mathrm{osc}}^2.
		\end{aligned}
		\]
		Let $u_1$ and $u_2$ be such that $u_1 \, u_2 = c_0 \, u$.
		Then it holds that
		\[
		\|f\|_{\pi}^2 \leq u \, \mathcal{E}_S(f) + \tilde{\beta}(u) \,  \|f\|_{\scriptsize\mathrm{osc}}^2.
		\]
		That $u \mapsto \tilde{\beta}(u)$ is a non-increasing function that goes to 0 as $u \to \infty$ is established in Section~3 of \cite{andrieu2022comparison}.
	\end{proof}

	\section{Derivations of Proposition \ref{pro:madras-1}, \ref{pro:qin}, \ref{pro:chen}, \ref{pro:hit}, and \ref{pro:doubly}} \label{app:derive}
	
		\subsection{Markov chain decomposition, a second form} \label{ssec:decomposition-1}
	
	Recall the setting of Section \ref{ssec:madras-1}:
	$(\X, \A, \pi_0)$ is a probability space.
	$\X= \bigcup_{z=1}^k \X_z$, where $\pi_0(\X_z) > 0$, and the $\X_z$'s may overlap.
	The Mtk $N: \X\times \A \to [0,1]$ is reversible with respect to $\pi_0$.
	Moreover, the Mtks $\Pi_0: [k] \times 2^{[k]} \to [0,1]$ and $H_z: \X_z \times \A_z \to [0,1]$ (where $\A_z$ is $\A$ restricted to $\X_z$) are defined via the formulas below:
	\[
	\begin{aligned}
		& \Pi_0(z,\{z'\}) = \frac{\pi_0(\X_z \cap \X_{z'})}{\Theta \, \pi_0(\X_z)} \, \ind_{\{z\}^c}(z') + \left[ 1 - \frac{\sum_{z'' \neq z} \pi_0(\X_z \cap \X_{z''})}{\Theta \, \pi_0(\X_z)} \right]\, \ind_{\{z\}}(z') , \\
		& H_z(x, \mathsf{A}) = N(x, \mathsf{A}) + N(x, \X_z^c) \, \delta_x(\mathsf{A}),
	\end{aligned}
	\]
	where $\Theta = \max_x \sum_{z=1}^k \ind_{\X_z}(x)$.
	
	The goal is to establish Proposition \ref{pro:madras-1} using Corollary \ref{cor:uniform}.
	To this end, we identify the elements in Section \ref{ssec:sandwich} as follows:
	\begin{itemize}
		\item $\pi: \A \to [0,1]$ is such that $\pi(\df x) = \pi_0(\df x) \, L(x) / \Delta$, where $L(x) = \sum_{z=1}^k \ind_{\X_z}(x)$, and $\Delta = \sum_{z=1}^k \pi_0(\X_z)$.
		
		\item $(\Y, \B) = (\X, \A)$, $(\Z, \C) = ([k], 2^{[k]})$.
		
		\item $\tilde{\pi}: \A \times 2^{[k]} \to [0,1]$ is such that
		\[
		\tilde{\pi}(\mathrm{A} \times \{z\}) = \frac{\pi_0(\mathsf{A} \cap \X_z)}{ \Delta },
		\]
		It is easy to check that $\tilde{\pi}(\mathsf{A} \times [k]) = \pi(\mathsf{A})$.
		
		\item For $z \in [k]$, $\varpi(z) = \pi_0(\X_z)/
		\Delta$, and $\varpi_z(\mathsf{A}) = \pi_0(\mathsf{A} \cap \X_z) / \pi_0(\X_z)$ if $\mathsf{A} \in \A$.
		It is easy to see that $\tilde{\pi}(\mathsf{A} \times \{z\}) = \varpi(\{z\}) \, \varpi_z(\mathsf{A})$.
		
		\item The Mtk $S: \X \times \A \to [0,1]$ has the form
		\[
		S(x, \df x') = \frac{N(x, \df x') \, L(x')}{\Theta} + \left( 1 - \frac{\int_{\X} N(x, \df x'') \, L(x'') }{\Theta} \right) \delta_x(\df x').
		\]
		
		\item For $f \in L_0^2(\pi)$ and $(x,z) \in \X \times [k]$, $Pf(x,z) = f(x)$.
		For $h \in L_0^2(\tilde{\pi})$ and $x \in \X$, $P^*h(x) = L(x)^{-1} \sum_{z=1}^k \ind_{\X_z}(x) \, h(x,z)$.
		
		\item The Mtk $Q: (\X \times [k]) \times (\mathcal{A} \times 2^{[k]}) \to [0,1]$ is of the form
		\[
		Q((x,z), \mathsf{A} \times \{z'\}) = \frac{N(x, \mathsf{A} \cap \X_{z'} )  }{\Theta} + \left( 1 - \frac{\int_{\X} N(x, \df x') \, L(x') }{\Theta} \right) \delta_x(\mathsf{A}) \, \ind_{\{z\}}(z').
		\]
		It is easy to see that $Q$ is reversible with respect to $\tilde{\pi}$.

		\item For $z \in [k]$, $x \in \X$, and $\mathsf{A} \in \A$,
		\[
		Q_z(x, \mathsf{A}) = N(x, \mathsf{A} \cap \X_z) + N(x, \X_z^c) \, \delta_x(\mathsf{A}).
		\]
		
		\item For $h \in L_0^2(\tilde{\pi})$ and $z \in [k]$, 
		\[
		\Proj h(z) = \frac{1}{\pi_0(\X_z)} \int_{\X_z} h(x,z) \, \pi_0(\df x).
		\]
		For $g \in L_0^2(\varpi)$ and $(x,z) \in \X \times [k]$, $\Proj^* g (x,z) = g(z)$.
	\end{itemize}
	
	One may then verify the following:
	\begin{itemize}
		\item $S = P^* Q P$, and \ref{H1} is satisfied.
		
		\item \ref{H2} and \ref{H3} are satisfied.
		
		\item \ref{H4} holds with $c_0 = 1/\Theta$.
		
		\item $\Proj P P^* \Proj^*$ corresponds to the Mtk
		\[
		\Proj P P^* \Proj^*(z, \{z'\}) = \frac{1}{\pi_0(\X_z)} \int_{\X_z \cap \X_{z'}} \frac{\pi_0(\df x)}{L(x)}.
		\]
		This Mtk is reversible with respect to $\varpi$.
	\end{itemize}
	
	By Corollary \ref{cor:uniform} and Lemma \ref{lem:iso-1},
	\begin{equation} \label{ine:overlap-1}
		\mathrm{Gap}(S) \geq \frac{1}{\Theta} \left[ \min_z \mathrm{Gap}(Q_z) \right] \, \mathrm{Gap}( P^* \Proj^* \Proj P) = \frac{1}{\Theta} \left[ \min_z \mathrm{Gap}(H_z) \right] \, \mathrm{Gap}(\Proj P P^* \Proj^*).
	\end{equation}
	
	The above display is almost the desired spectral gap decomposition.
	To derive Proposition \ref{pro:madras-1} precisely, one can utilize some simple arguments from \cite{madras2002markov}, which we include for completeness.
	Note that, for $g \in L_0^2(\varpi)$,
	\[
	\begin{aligned}
		\mathcal{E}_{\Proj P P^* \Proj^*}(g) &= \frac{1}{2 \Delta} \sum_{z,z'} \int_{\X_z \cap \X_{z'}} \frac{\pi_0(\df x)}{L(x)} [g(z') - g(z)]^2 \\
		&\geq \frac{1}{2 \Delta \Theta} \sum_{z,z'} \int_{\X_z \cap \X_{z'}} \pi_0(\df x) \, [g(z') - g(z)]^2 \\
		&= \mathcal{E}_{\Pi_0}(g).
	\end{aligned}
	\]
	Thus, 
	\begin{equation} \label{ine:overlap-2}
		\mathrm{Gap}(\Proj P P^* \Proj^*) \geq \mathrm{Gap}(\Pi_0).
	\end{equation}
	
	Next, we use a standard argument \citep[see, e.g.,][]{diaconis1996logarithmic,madras1999importance,madras2002markov}.
	Observe that $L^2(\pi) = L^2(\pi_0)$.
	Then, for $f \in L^2(\pi_0)$,
	\[
	\begin{aligned}
		\mathcal{E}_N(f - \pi_0 f) &= \frac{1}{2} \int_{\X^2} \pi_0(\df x) \, N(x, \df x') \, [f(x') - f(x)]^2 \\
		&= \frac{\Delta \Theta}{2} \int_{\X^2}  \pi(\df x) \, \frac{1}{L(x) \, L(x')} \, S(x, \df x') \, [f(x') - f(x)]^2 \\
		& \geq \frac{\Delta}{\Theta} \, \mathcal{E}_S(f - \pi f).
	\end{aligned}
	\]
	Moreover, since $\pi_0(\cdot) \leq \Delta \, \pi(\cdot)$,
	\[
	\begin{aligned}
		\|f - \pi_0 f \|_{\pi_0}^2 &= \int_{\X} \pi_0(\df x) \, [f(x) - \pi_0 f]^2 \\
		&\leq   \int_{\X} \pi_0(\df x)  \, [f(x) - \pi f]^2 \\
		&\leq \Delta \|f - \pi f\|_{\pi}^2.
	\end{aligned}
	\]
	As a consequence,
	\begin{equation} \label{ine:overlap-3}
		\mathrm{Gap}(N) = \inf_{f} \frac{\mathcal{E}_N(f - \pi_0 f)}{\|f - \pi_0 f\|_{\pi_0}^2} \geq \frac{1}{\Theta} \inf_{f} \frac{\mathcal{E}_S(f - \pi f)}{\|f - \pi f\|_{\pi}^2} = \frac{1}{\Theta} \mathrm{Gap}(S).
	\end{equation}
	
	Proposition \ref{pro:madras-1} then follows from \eqref{ine:overlap-1}, \eqref{ine:overlap-2}, and \eqref{ine:overlap-3}.
	
		\subsection{Random-scan hybrid Gibbs samplers} \label{ssec:gibbs}
	
	Recall the setting of Section \ref{ssec:qin}:
	$\X = \X_1 \times \cdots \times \X_k$, where each $\X_i$ is a Polish space with Borel algebra $\A_i$.
	The distribution of interest $\pi$ is defined on $(\X, \A)$, where $\A = \A_1 \times \cdots \times \A_k$.
	If $(X_1, \dots, X_k) \sim \pi$, $\varphi_{i,u}$ denotes the conditional distribution of $X_i \mid X_{-i} = u$ for $i \in [k]$ and $u \in \X_{-i}$.
	$H_{i,u}: \X_i \times \A_i \to [0,1]$ is an Mtk that is reversible with respect to $\varphi_{i,u}$ such that $x \mapsto H_{i,x_{-i}}(x_i, \mathsf{A}_i)$ is measurable for $\mathsf{A}_i \in \A_i$.
	Finally,
	\[
	\bar{S}(x, \df x') = \sum_{i=1}^{k} p_i \, \varphi_{i, x_{-i}}(\df x'_i) \, \delta_{x_{-i}}(\df x'_{-i}), \quad S(x, \df x') = \sum_{i=1}^{k} p_i \, H_{i, x_{-i}}(x_i, \df x'_i) \, \delta_{x_{-i}}(\df x'_{-i}),
	\]
	where $(p_1, \dots, p_k)$ is a probability vector.
	
	Our goal is to demonstrate that Proposition \ref{pro:qin} is implied by Theorem \ref{thm:uniform}.
	To this end, we derive Proposition \ref{pro:qin} using Proposition \ref{pro:andrieu}, which is in turn a special case of Theorem \ref{thm:uniform}.
	To be specific, we will show that $\bar{S}$ and $S$ can be viewed as special cases of the ideal and hybrid data augmentation Mtks studied in Sections \ref{ssec:andrieu} and \ref{ssec:da}.
	
	We identify the elements in Sections \ref{ssec:andrieu} and \ref{ssec:da} as follows:
	\begin{itemize}
		\item $\Z = \bigcup_{i = 1}^k \{i\} \times \X_{-i}$, while  $\C$ is the sigma algebra generated from sets of the form $\{i\} \times \mathsf{C}_i$, where $\mathsf{C}_i \in \A_{-i} := \A_1 \times \cdots \times \A_{i-1} \times \A_{i+1} \times \cdots \times \A_k$.
		
		\item $\tilde{\pi}: \A \times \C \to [0,1]$ is the joint distribution of $(X, (\eta, X_{-\eta}))$, where $X \sim \pi$, and independently, $\eta \in [k]$ is distributed according to the probability vector $(p_1, \dots, p_k)$.
		
		\item $\varpi: \C \to [0,1]$ is the marginal distribution of $(\eta, X_{-\eta})$.
		
		\item For $x \in \X$,
		$\pi_x: \C \to [0,1]$ is the conditional distribution of $(\eta, X_{-\eta}) \mid X = x$, i.e.,
		$\pi_x(\{i\} \times \mathsf{C}_i) = p_i \, \delta_{x_{-i}}(\mathsf{C}_i)$ for $i \in [k]$ and $\mathsf{C}_i \in \A_{-i}$.
		
		\item For $i \in [k]$ and $u \in \X_{-i}$,
		$\varpi_{i,u}: \A \to [0,1]$ is the conditional distribution of $X \mid (\eta, X_{-\eta}) = (i, u)$, i.e., $\varpi_{i,u}(\df x) = \varphi_{i,u}(\df x_i) \, \delta_u(\df x_{-i})$.

		\item For $i \in [k]$ and $u \in \X_{-i}$, $Q_{i,u}: \X \times \A \to [0,1]$ is such that
		\[
		Q_{i,u}(x, \df x') = H_{i,u}(x_i, \df x'_i) \, \delta_u(\df x'_{-i}).
		\]
		It is easy to see that $Q_{i,u}$ is reversible with respect to $\varpi_{i,u}$.
		
	\end{itemize}

	One may then check the following:
	\[
	\bar{S}(x, \df x') = \sum_{i=1}^k p_i \int_{\X_{-i}} \delta_{x_{-i}}(\df u) \, \varphi_{i,u}(\df x'_i) \, \delta_u(\df x'_{-i}) = \sum_{i=1}^k \int_{\X_{-i}} \pi_x(\{i\}, \df u) \, \varpi_{i,u}(\df x'),
	\]
	\[
	S(x, \df x') = \sum_{i=1}^k p_i \int_{\X_{-i}} \delta_{x_{-i}}(\df u) \, H_{i,u}(x_i, \df x'_i) \, \delta_u(\df x'_{-i}) = \sum_{i=1}^k \int_{\X_{-i}} \pi_x(\{i\}, \df u) \, Q_{i,u}(x, \df x').
	\]
	Hence, $\bar{S}$ and $S$ can be viewed as special cases of the ideal and hybrid data augmentation Mtks in Sections \ref{ssec:andrieu} and \ref{ssec:da}.
	
	By Proposition \ref{pro:andrieu}, for $f \in L_0^2(\pi)$,
	\[
	\mathcal{E}_S(f) \geq \left[\min_{i} \essinf_u \mathrm{Gap}(Q_{i,u}) \right] \, \mathcal{E}_{\bar{S}}(f), \quad \mathrm{Gap}(S) \geq \left[\min_{i} \essinf_u \mathrm{Gap}(Q_{i,u}) \right] \, \mathrm{Gap}(\bar{S}).
	\]
	It is easy to see $\mathrm{Gap}(Q_{i,u}) = \mathrm{Gap}(H_{i,u})$.
	Then Proposition \ref{pro:qin} follows.
	
	\subsection{Approximate variance conservation in localization schemes} \label{ssec:local}
	
	Recall the setup of Section \ref{ssec:chen}:
	For $s \in \mathbb{N}$, $\mathsf{W}_s$ is a Polish space and $\mathcal{D}_s$ is its Borel algebra.
	$(W_t)_{t=0}^{\infty}$ is a sequence of random elements such that $W_s$ takes values in $\mathsf{W}_s$.
	The distribution of $(W_s)_{s=0}^t$ is denoted by $\mu_t$, and $\Z_t = \mathsf{W}_0 \times \cdots \times \mathsf{W}_t$.
	The function $v_t: \Z_t \times \X \to [0,\infty]$ satisfies \ref{A1} and \ref{A2}, i.e., $v_t(z,x) \, \pi(\df x)$ defines a probability measure on $(\X,\A)$ for $\mu_t$-a.e. $z \in \Z_t$, and $v_t((W_s)_{s=0}^t,x)$, $t \in \mathbb{N}$, is a martingale initialized at~1 for $x \in \X$. 
	The following Mtk is reversible with respect to $\pi$:
	\[
	K_t(x, \df x') = \int_{\Z_t} v_t(z, x) \, v_t(z, x') \, \pi(\df x') \, \mu_t(\df z).
	\]
	Our goal is to show that Proposition \ref{pro:chen} is implied by Corollary \ref{cor:uniform}.
	To this end, we demonstrate that Proposition \ref{pro:chen} can be derived using Proposition \ref{pro:andrieu}, which is encompassed by Corollary \ref{cor:uniform}.

	Fix $t \in \mathbb{N}$ and $s \in \{0, \dots, t-1\}$.
	To utilize Proposition \ref{pro:andrieu}, we shall demonstrate that $K_s$ corresponds to a data augmentation algorithm, while $K_{s+1}$ is associated with a particular hybrid data augmentation algorithm, as described in Sections \ref{ssec:andrieu} and \ref{ssec:da}.
	The elements in Sections \ref{ssec:andrieu} and \ref{ssec:da} are identified below.
	
	\begin{itemize}
		\item $\Z = \Z_s$ and $\C$ is the corresponding Borel algebra.
		
		\item $\tilde{\pi}(\df x, \df z) = v_s(z, x) \, \pi(\df x) \, \mu_s(\df z)$.
		
		\item $\varpi = \mu_s$.
		
		\item For $x \in \X$, $\pi_x(\df z) = v_s(z,x) \, \mu_s(\df z)$.
		
		\item For $z \in \Z$, $\varpi_z(\df x) = v_s(z,x) \, \pi(\df x)$.
		Note that the random probability measure $\nu_s$ defined in Section \ref{ssec:chen} can be seen as $\varpi_z$ with $z = (W_i)_{i=0}^s$.
		
		\item For $z \in \Z$, let $\mu_{s,z}$ be the conditional distribution of $W_{s+1}$ given $(W_i)_{i=0}^s = z$, and define the Mtk $Q_z: \X \times \A \to [0,1]$ as follows:
		\[
		Q_z(x, \df x') = \int_{\mathsf{W}_{s+1}} \mu_{s,z}(\df w) \, \frac{v_{s+1}((z,w), x)}{v_s(z, x)} \, v_{s+1}((z,w), x') \, \pi(\df x').
		\]
		In other words,
		\[
		Q_z(x, \mathsf{A}) = \mathrm{E} \left[ \frac{\df \nu_{s+1} }{\df \nu_s}(x) \, \nu_{s+1}(\mathsf{A}) \mid (W_i)_{i=0}^s = z\right].
		\]
		Observe that $Q_z$ corresponds to a data augmentation algorithm targeting $\varpi_z$.
		Indeed, due to \ref{A1} and \ref{A2}, we may define the probability measure $\tilde{\pi}_z: \A \times \mathcal{D}_{t+1} \to [0,1]$ as follows:
		\[
		\tilde{\pi}_z(\df (x, w)) = \varpi_z(\df x) \, \pi_{x,z}'(\df w) = \varpi_{w,z}'(\df x) \, \mu_{s,z}(\df w) = v_{s+1}((z,w), x) \, \pi(\df x) \, \mu_{s,z}(\df w),
		\]
		where $\pi_{x,z}'(\df w) = [v_{s+1}((z,w),x)/v_s(z,x)] \, \mu_{s,z}(\df w)$ and $\varpi_{w,z}'(\df x) = v_{s+1}((z,w),x) \, \pi(\df x)$ are conditional distributions.
		Then $Q_z(x, \df x') = \int_{\mathsf{W}_{s+1}} \pi_{x,z}'(\df w) \, \varpi_{w,z}'(\df x')$.
	\end{itemize}
	
	One may then check the following:
	\[
	\begin{aligned}
		\int_{\Z} \pi_x(\df z) \, \varpi_z(\df x') = K_s(x, \df x'), \quad
		\int_{\Z} \pi_x(\df z)  \, Q_z(x, \df x') = K_{s+1}(x, \df x').
	\end{aligned}
	\]
	Then $K_s$ and $K_{s+1}$ can be viewed as (hybrid) data augmentation algorithms.
	
	By Proposition \ref{pro:andrieu},
	\begin{equation} \label{ine:andrieu-chen}
		\mathrm{Gap}(K_{s+1}) \geq \left[ \essinf_z \mathrm{Gap}(Q_{s,z}) \right] \, \mathrm{Gap}(K_s),
	\end{equation}
	where we denote $Q_z$ by $Q_{s,z}$ to emphasize its dependence on~$s$, and the essential infimum is defined with respect to $\varpi = \mu_s$.
	
	We now establish Part (a) of Proposition \ref{pro:chen}.
	Next, we note that \eqref{eq:GapKt} in Proposition \ref{pro:chen} is a manifestation of the following result for generic data augmentation algorithms, which may be traced back to at least \cite{liu1995covariance}:
	\begin{lemma} \citep{liu1995covariance} \label{lem:da-gap-variance}
		For a generic data augmentation algorithm described in Sections \ref{ssec:andrieu} and \ref{ssec:da},
		\begin{equation} \nonumber
			\begin{aligned}
				\mathrm{Gap}(\bar{S}) &= \inf_{f \in L_0^2(\pi) \setminus \{0\}} \frac{\|f\|_{\pi}^2 - \int_{\Z} (\varpi_z f)^2 \, \varpi(\df z) }{\|f\|_{\pi}^2} \\
				&= \inf_{f \in L_0^2(\pi) \setminus \{0\}} \frac{\int_{\Z} \mathrm{var}_{\varpi_z}(f) \, \varpi(\df z)}{\mathrm{var}_{\pi}(f)} .
			\end{aligned}
		\end{equation}
	\end{lemma}
	
	{ 
	\begin{proof}
		For every $f\in L^2(\pi)$, the conditional variance decomposition gives
		\[
		\|f - \pi f\|_{\pi}^2 = \mathrm{var}_{\pi}(f)
		=
		\int_{\Z}\mathrm{var}_{\varpi_z}(f)\,\varpi(\df z)
		+
		\mathrm{var}_{\varpi}(\varpi_z f).
		\]
		If $f\in L_0^2(\pi)$, then $\int_{\Z}\varpi_z f\,\varpi(\df z)=\pi f=0$, and it can be checked that
		\[
		\langle f,\bar S f\rangle_{\pi}
		=
		\int_{\Z}(\varpi_z f)^2\,\varpi(\df z)
		=
		\mathrm{var}_{\varpi}(\varpi_z f).
		\]
		Consequently, for $f \in L_0^2(\pi)$,
		\[
		\mathcal{E}_{\bar S}(f)
		= \|f\|_{\pi}^2-
		\langle f,\bar S f\rangle_{\pi} =
		\|f\|_{\pi}^2-
		\int_{\Z}(\varpi_z f)^2\,\varpi(\df z)
		=
		\int_{\Z}\mathrm{var}_{\varpi_z}(f)\,\varpi(\df z).
		\]
		Dividing this identity by $\|f\|_{\pi}^2=\mathrm{var}_{\pi}(f)$ and taking the infimum over $f\in L_0^2(\pi)\setminus\{0\}$ proves both equalities.
	\end{proof}
	}

	By Lemma \ref{lem:da-gap-variance},
	\begin{equation} \label{eq:da-gap-variance-1}
		\begin{aligned}
			\mathrm{Gap}(Q_{s,z}) &=  \inf_{f \in L_0^2(\varpi_z) \setminus \{0\} } \frac{\int_{\mathsf{W}_{s+1}} \mathrm{var}_{\varpi_{w,z}'}(f) \, \mu_{s,z}(\df w)  }{\mathrm{var}_{\varpi_z}(f) } \\
			&= \inf_{f \in L_0^2(\varpi_z) \setminus \{0\} } \frac{ \mathrm{E} [\mathrm{var}_{\nu_{s+1}}(f) \mid (W_i)_{i=0}^s = z] }{\mathrm{var}_{\nu_s}(f) |_{(W_i)_{i=0}^s = z} }.
		\end{aligned}
	\end{equation}
	
	Continuing to Part (b) of Proposition \ref{pro:chen}. 
	Fix $t$ but let $s$ vary.
	Assume that there are positive constants $\kappa_1, \dots, \kappa_t$ such that, given $s \in \{0,\dots,t-1\}$, almost surely, $\mathrm{E}[\mathrm{var}_{\nu_{s+1}}(f) \mid (W_i)_{i=0}^s] \geq \kappa_{s+1} \mathrm{var}_{\nu_s}(f)$ for $f \in L_0^2(\nu_s)$.
	Then, by \eqref{eq:da-gap-variance-1}, for $\mu_s$-a.e. $z \in \Z_s$, $\mathrm{Gap}(Q_{s,z}) \geq \kappa_{s+1}$.
	It then follows from \eqref{ine:andrieu-chen} that 
	\[
	\mathrm{Gap}(K_t) \geq \left( \prod_{s=1}^t \kappa_s \right) \mathrm{Gap}(K_0).
	\]
	But $K_0(x, \df x') = \pi(\df x')$, so $\mathrm{Gap}(K_0) = 1$ by definition, and the above formula is precisely \eqref{ine:chen-decomposition} in Proposition \ref{pro:chen}.

	\subsection{Hybrid hit-and-run samplers} \label{ssec:hit-1}
	
	Recall the setup of Section \ref{ssec:hit-0}:
	$\X = \mathbb{R}^k$, and $\pi$ admits a density function $\dot{\pi}$.
	$\mathsf{W}$ is the Stiefel manifold consisting of ordered orthonormal bases of $\ell$-dimensional linear subspaces of~$\X$.
	$\nu$ is a distribution on $\mathsf{W}$.
	For $x \in \X$ and $w \in \mathsf{W}$, $\theta(x,w) = x - \sum_{i=1}^{\ell} (w_i^{\top} x) w_i$.
	Moreover, $\varphi_{w, \theta(x,w)}$ is a probability measure on $\mathbb{R}^{\ell}$ with density function proportional to $u \mapsto \dot{\pi}(\theta(x,w) + \sum_{i=1}^{\ell} u_i w_i )$.
	The Mtks $\bar{S}: \X \times \A \to [0,1]$ and $S: \X \times \A \to [0,1]$ are of the form
	\[
	\begin{aligned}
		\bar{S}(x, \mathsf{A}) &= \int_{\mathsf{W}} \nu(\df w) \int_{\mathbb{R}^{\ell}} \varphi_{w, \theta(x,w)}(\df u) \, \ind_{\mathsf{A}}\left(\theta(x,w) + \sum_{i=1}^{\ell} u_i w_i \right), \\
		S(x, \mathsf{A}) &= \int_{\mathsf{W}} \nu(\df w) \int_{\mathbb{R}^{\ell}} H_{w, \theta(x,w)}\left((w_i^{\top} x)_{i=1}^{\ell}, \df u \right) \, \ind_{\mathsf{A}}\left(\theta(x,w) + \sum_{i=1}^{\ell} u_i w_i \right) .
	\end{aligned}
	\]
	Here, $H_{w, \theta(x,w)}$ is an Mtk that is reversible with respect to $\varphi_{w, \theta(x,w)}$.
	The function $(x,w,u) \mapsto H_{w, \theta(x,w)}(u, \mathsf{B})$ is measurable for every measurable $\mathsf{B} \subset \mathbb{R}^{\ell}$.
	
	Our goal is to prove Proposition \ref{pro:hit}.
	To this end, we show that $\bar{S}$ defines a data augmentation algorithm, and $S$ defines a hybrid version of that algorithm, as described in Sections~\ref{ssec:andrieu} and~\ref{ssec:da}.
	
	The elements in Sections \ref{ssec:andrieu} and \ref{ssec:da} are identified below.
	\begin{itemize}
		\item $\Z = \mathsf{W} \times \X$, and $\C$ is its Borel algebra.
		
		\item $\tilde{\pi}(\df (x, w, x')) = \pi(\df x) \, \pi_x(\df (w, x'))$, where, for $x \in \X$, $\pi_x(\df (w, x')) = \nu(\df w) \, \delta_{\theta(x,w)}(\df x')$.
		In other words, $\tilde{\pi}$ is the joint distribution of $X$ and $(W, \theta(X,W))$, where $X \sim \pi$ and, independently, $W \sim \nu$.
		
		\item 
		
		Identify
		$\varpi: \C \to [0,1]$ and $\varpi_{w,x'}: \A \to [0,1]$, where $(w,x') \in \Z$, as follows:
		\[
		\begin{aligned}
			&\varpi(\df (w, x')) =  \nu(\df w) \int_{\X} \pi(\df x) \, \delta_{\theta(x,w)}(\df x'), \\
			&\varpi_{w,x'}(\df x) = \int_{\mathbb{R}^{\ell}} \varphi_{w, \theta(x',w)}(\df u) \, \delta_{\theta(x',w) + \sum_{i=1}^{\ell} u_i w_i }(\df x).
		\end{aligned}
		\]
		With a bit of work, one can check that these are indeed the correct marginal and conditional distributions, i.e., $\tilde{\pi}(\df (x,w,x')) =  \varpi_{w,x'}(\df x) \, \varpi(\df (w, x'))$.
		
		\item For $(w,x') \in \Z$, the Mtk $Q_{w,x'}: \X \times \A \to [0,1]$ is such that
		\[
		Q_{w,x'}(x, \mathsf{A}) = \int_{\mathbb{R}^{\ell}} H_{w, \theta(x',w)}\left( (w_i^{\top} x)_{i=1}^{\ell}, \df u \right) \ind_{\mathsf{A}}\left(\theta(x',w) + \sum_{i=1}^{\ell} u_i w_i \right).
		\]
		It is not difficult to show that $Q_{w,x'}$ is reversible with respect to $\varpi_{w,x'}$.
	\end{itemize}
	
	One can then check that, for $x \in \X$ and $\mathsf{A} \in \A$,
	\[
	\int_{\mathsf{W} \times \X} \pi_x(\df (w, x')) \, \varpi_{w,x'}(\mathsf{A}) = \bar{S}(x, \mathsf{A}), \quad \int_{\mathsf{W} \times \X} \pi_x(\df (w, x')) \, Q_{w,x'}(x, \mathsf{A}) = S(x, \mathsf{A}).
	\]
	When verifying this, it is useful to note that, for $x \in \X$ and $w \in \mathsf{W}$, if $x' = \theta(x,w)$, then $\theta(x',w) = x' = \theta(x,w)$.

	We can now establish Proposition \ref{pro:hit} using Proposition \ref{pro:andrieu}.
	Reversibility is implied by the fact that $\bar{S}$ and $S$ are special cases of the (hybrid) data augmentation algorithms defined in Section \ref{ssec:andrieu}.
	By Proposition~\ref{pro:andrieu}, for $f \in L_0^2(\pi)$,
	\[
	\mathcal{E}_S(f) \geq \left[ \essinf_{w,x} \mathrm{Gap}(Q_{w,x}) \right] \mathcal{E}_{\bar{S}}(f), \quad \mathrm{Gap}(S) \geq \left[ \essinf_{w,x} \mathrm{Gap}(Q_{w,x}) \right] \mathrm{Gap}(\bar{S}).
	\]
	By Lemma \ref{lem:iso-2} in Appendix~\ref{app:lemma}, $\mathrm{Gap}(Q_{w,x}) = \mathrm{Gap}(H_{w, \theta(x,w)})$.
	This establishes the spectral decomposition formula in Proposition \ref{pro:hit}.

	\subsection{Hybrid data augmentation algorithms with two intractable conditionals} \label{ssec:doubly-1}
	
	Recall the setup in Section \ref{ssec:doubly-0}.
	$(\X_1, \A_1, \varphi_1)$ and $(\X_2, \A_2, \varphi_2)$ are two probability spaces.
	The probability measure $\varphi: \A_1 \times \A_2 \to [0,1]$ has the disintegration
	\[
	\varphi(\df (x_1, x_2)) = \varphi_1(\df x_1) \, \varphi_{2, x_1}(\df x_2) = \varphi_2(\df x_2) \, \varphi_{1, x_2}(\df x_1).
	\]
	For $i \in \{1,2\}$ and $x_{-i} \in \X_{-i}$, $H_{i,x_{-i}}: \X_i \times \A_i \to [0,1]$ is an Mtk that is reversible with respect to $\varphi_{i,x_{-i}}$, and  $(x_1, x_2) \mapsto H_{i,x_{-i}}(x_i, \mathsf{A})$ is measurable if $\mathsf{A} \in \A_i$.
	
	The following Mtks are defined:
	\[
	\begin{aligned}
		&\hat{S}_1(x_2, \df x_2') = \int_{\X_1}  \varphi_{1,x_2}(\df x_1) \, H_{2,x_1}(x_2, \df x_2'), \quad \hat{S}_2(x_2, \df x_2') = \int_{\X_1}  \varphi_{1,x_2}(\df x_1) \, H_{2,x_1}^2(x_2, \df x_2'), \\
		&T((x_1, x_2), \df (x_1', x_2'))= H_{1,x_2}(x_1, \df x_1') \, H_{2,x_1'}(x_2, \df x_2').
	\end{aligned}
	\]
	The goal of this subsection is to prove Proposition \ref{pro:doubly}, which states that
	\begin{equation} \label{ine:doubly}
		1 - \|T\|_{\varphi}^2 \geq \left( 1 - \esssup_{x_2} \|H_{1,x_2}\|_{\varphi_{1,x_2}}^2 \right) \mathrm{Gap}(\hat{S}_2).
	\end{equation}
	The proposition also asserts that, if $H_{2,x_1}$ is almost always positive semi-definite, then $\mathrm{Gap}(\hat{S}_2) \geq \mathrm{Gap}(\hat{S}_1)$.

	We begin by establishing \eqref{ine:doubly}.
	We shall apply Corollary \ref{cor:nonreversible}.
	The elements in Sections \ref{ssec:sandwich} and \ref{ssec:nonreveresible} are specified below.
	
	\begin{itemize}
		\item $(\X, \A) = (\X_1 \times \X_2, \A_1 \times \A_2)$, $(\Y, \B) = (\X_1, \A_1)$, $(\Z, \C) = (\X_2, \A_2)$.
		
		\item $\pi = \tilde{\pi} = \varphi$.
		
		\item $\varpi = \varphi_2$, and $\varpi_{x_2} = \varphi_{1,x_2}$ for $x_2 \in \X_2$.
		
		\item $P: L_0^2(\pi) \to L_0^2(\tilde{\pi})$ is such that, for $f \in L_0^2(\pi)$ and $(x_1, x_2) \in \X$,
		\[
		P f (x_1, x_2) = \int_{\X_2} f(x_1, x_2') \, H_{2,x_1}(x_2, \df x_2').
		\]
		Note that $P = P^*$.

		\item $Q_{x_2} = H_{1,x_2}^2$ for $x_2 \in \X_2$.
		
		\item $Q: \X \times \A \to [0,1]$ is given by
		\[
		Q((x_1, x_2), \df (x_1',x_2')) = Q_{x_2}(x_1, \df x_1') \, \delta_{x_2}(\df x_2').
		\]
		It is clearly reversible with respect to $\tilde{\pi}$.
		
		\item $R: L_0^2(\tilde{\pi}) \to L_0^2(\pi)$ is such that, for $h \in L_0^2(\tilde{\pi})$ and $(x_1, x_2) \in \X$,
		\[
		R h (x_1, x_2) = \int_{\X_1} h(x_1', x_2) \, H_{1,x_2}(x_1, \df x_1').
		\]
		Note that $R = R^*$.
		
		\item $\Proj: L_0^2(\tilde{\pi}) \to L_0^2(\varphi_2)$ is such that, for $h \in L_0^2(\tilde{\pi})$ and $x_2 \in \X_2$,
		\[
		\Proj h (x_2) = \int_{\X_1} h(x_1, x_2) \, \varphi_{1,x_2}(\df x_1).
		\]
		
		\item For $g \in L_0^2(\varphi_2)$ and $(x_1, x_2) \in \X$, $\Proj^* g(x_1, x_2) = g(x_2)$.
		
	\end{itemize}
	
	One may then verify the following:
	
	\begin{itemize}
		
		\item $P$ is an Mtk that is reversible with respect to $\tilde{\pi}$, so $\|P^* P\|_{\pi} = \|P\|_{\tilde{\pi}}^2 \leq 1$.
		Thus, $P$ satisfies \ref{H2}.
		
		\item $Q_{x_2}$ is reversible with respect to $\varphi_{1,x_2}$ for $x_2 \in \X_2$, and \ref{H3} is satisfied by $\{Q_{x_2}\}_{x_2 \in \X_2}$.
		
		\item For $(x_1, x_2) \in \X_1 \times \X_2$ and $\mathsf{B} \in \B$,
		\[
		Q((x_1, x_2), (\mathsf{B} \setminus \{x_1\}) \times \{x_2\}) = Q_{x_2}(x_1, \mathsf{B} \setminus \{x_1\}).
		\]
		Thus, \ref{H4} is satisfied by~$Q$ with $c_0 = 1$.
		
		\item $Q = R^*R$, $T = RP$.

	\end{itemize}
	
	The above implies that \ref{C1} in Section \ref{ssec:nonreveresible} is satisfied.
	Since $H_{1,x_2}$ is reversible for $x_2 \in \X_2$, it holds that $\mathrm{Gap}(Q_{x_2}) = 1 - \|H_{1,x_2}\|_{\varphi_{1,x_2}}^2$.
	By Corollary \ref{cor:nonreversible},
	\begin{equation} \label{ine:doubly-1}
		\begin{aligned}
			1 - \|T\|_{\pi}^2 &\geq \left[ \essinf_{x_2} \mathrm{Gap}(Q_{x_2}) \right]  \mathrm{Gap}(P^* \Proj^* \Proj P) 
			= \left( 1 - \esssup_{x_2} \|H_{1,{x_2}}\|_{\varphi_{1,x_2}}^2 \right) \mathrm{Gap}( \Proj P P^* \Proj^*),
		\end{aligned}
	\end{equation}
	where the essential supremum is taken with respect to $\varphi_2$.
	It is easy to check that $\Proj P P^* \Proj^* = \hat{S}_2$.
	Then \eqref{ine:doubly} follows from \eqref{ine:doubly-1}.
	
	We now establish the second assertion in Proposition \ref{pro:doubly}.
	Assume that $H_{2,x_1}$ is positive semi-definite.
	Then $H_{2,x_1} - H_{2,x_1}^2$ is positive semi-definite \citep[][\S 31, Theorem 2]{helmberg2014introduction}.
	Then, applying Proposition~\ref{pro:loewner} in Section \ref{sec:simple} to the hybrid data augmentation sampler with one intractable conditional, one obtains $\mathrm{Gap}(\hat{S}_2) \geq \mathrm{Gap}(\hat{S}_1)$.
	
%
%
	
	\section{Decomposition Without the Sandwich Structure: Ge, Lee, and Risteski's Soft Density Decomposition} \label{app:beyond}
	
	{ 
		As argued in Remark~\ref{rem:H4}, one may replace \ref{H1} and \ref{H4} with \ref{H5}, and still arrive at Theorem~\ref{thm:uniform}.
		As an example, we reproduce a ``soft density decomposition" result from \cite{ge2018simulated} using Remark~\ref{rem:H4} and Theorem~\ref{thm:uniform}.
		
		Let $S: \X \times \A \to [0,1]$ be reversible with respect to~$\pi$.
		As in Section \ref{ssec:andrieu}, assume that there is a joint distribution of the form $\tilde{\pi}(\df (x,z)) = \pi(\df x) \, \pi_x(\df z) = \varpi(\df z) \, \varpi_z(\df x)$.
		
		The following is a discrete-time analogue and extension of Theorem~6.1 of \cite{ge2018simulated}.
		The theorem is originally stated for a finite $\Z$ and continuous-time Markov chains.
		Here, we state it for an arbitrary~$\Z$ and discrete-time Markov chains.
		
		\begin{proposition} \label{pro:ge} \citep{ge2018simulated}
			Suppose that, for $\varpi$-a.e. $z \in \Z$, there is an Mtk $Q_z: \X \times \A \to [0,1]$ that is reversible with respect to $\varpi_z$, and assume that $(x,z) \mapsto Q_z(x, \mathsf{A})$ is measurable for $\mathsf{A} \in \A$.
			Suppose further that, for $f \in L_0^2(\pi)$, 
			\begin{equation} \label{eq:soft-decomposition}
				\mathcal{E}_S(f) \geq \frac{1}{2} \int_{\Z} \varpi(\df z)   \int_{\X^2} \varpi_z(\df x) \, Q_z(x, \df x') \, [f(x) - f(x')]^2.
			\end{equation}
			(In \cite{ge2018simulated}, the authors use a slightly stronger condition with the inequality replaced by equality.)
			Furthermore, suppose that there is a positive number $\bar{\kappa}$ such that, for $g \in L_0^2(\varpi)$,
			\begin{equation} \label{ine:soft-projection}
				\frac{1}{2} \int_{\Z^2} [g(z) - g(z')]^2 \, \frac{\varpi(\df z) \, \varpi(\df z')}{ \chi_{\max}^2(\varpi_z || \varpi_{z'}) } \geq \bar{\kappa} \|g\|_{\varpi}^2,
			\end{equation}
			where
			\[
			\chi_{\max}^2(\varpi_z || \varpi_{z'}) = \max\left\{ \left\| \frac{\df \varpi_{z'}}{\df \varpi_z} - 1 \right\|_{\varpi_z}^2, \left\| \frac{\df \varpi_z}{\df \varpi_{z'}} - 1 \right\|_{\varpi_{z'}}^2 \right\}.
			\]
			(We assume that $\chi_{\max}^2(\varpi_z || \varpi_{z'})$ is a measurable function of $(z,z')$, that it is infinity if $\varpi_z$ and $\varpi_{z'}$ are not equivalent, and we use the convention $0/0 = 0$, $1/0 = \infty$, $1/\infty = 0$.)
			Finally, assume that there is a constant $\kappa_{\dagger} \in [0,1]$ such that $\mathrm{Gap}(Q_z) \geq \kappa_{\dagger}$ for $\varpi$-a.e. $z \in \Z$.
			Then
			\[
			\mathrm{Gap}(S) \geq \frac{\kappa_{\dagger}}{1 + 1/(2\bar{\kappa})} .
			\]
		\end{proposition}
		
		\begin{proof}
			In terms of the objects in Section~\ref{sec:unified}, let $\Y = \X$, so $\tilde{\pi}$ is also defined on the space $\Y \times \Z$.
			For $f \in L_0^2(\pi)$ and $(y,z) \in \Y \times \Z$, let $Pf(y,z) = f(y)$, so that $P$ maps a function in $L_0^2(\pi)$ to one in $L_0^2(\tilde{\pi})$.
			Then, for $h \in L_0^2(\tilde{\pi})$ and $x \in \X$, $P^*h(x) = \int_{\Z} h(x,z) \, \pi_x(\df z)$.
			Condition \ref{H2} clearly holds.
			By the assumptions on $Q_z$, \ref{H3} holds.
			By \eqref{eq:soft-decomposition} and the two disintegrations of $\tilde{\pi}$, 
			\[
			\mathcal{E}_S(f) \geq \frac{1}{2} \int_{\X} \pi(\df x) \int_{\Z \times \X} \pi_x(\df z) \, Q_z(x, \df x') \, [f(x) - f(x')]^2 = \mathcal{E}_{P^*\hat{Q}P}(f),
			\]
			so \ref{H5} holds with $c_0 = 1$.
			In light of Remark~\ref{rem:H4}, we can claim that
			\begin{equation} \label{ine:soft-pre}
				\mathrm{Gap}(S) \geq \kappa_{\dagger} \mathrm{Gap}(P^*E^*EP),
			\end{equation}
			where $E: L_0^2(\tilde{\pi}) \to L_0^2(\varpi)$ is defined in \eqref{eq:Proj}.
			For $f \in L_0^2(\pi)$ and $x \in \X$, it is easy to verify that
			\[
			P^*E^*EP f(x) = \int_{\Z} \pi_x(\df z) \, \varpi_z(\df x') f(x').
			\]
			That is, $\bar{S} = P^*E^*EP $ is the data augmentation algorithm described in Section~\ref{ssec:andrieu}.
			In fact, the bound \eqref{ine:soft-pre} is just the spectral gap bound from Proposition~\ref{pro:andrieu}, with the assumption $$S(x, \df x') = \int_{\Z} \pi_x(\df z) \, Q_z(x, \df x') = P^*\hat{Q}P(x, \df x')$$ replaced by the softer assignment $\mathcal{E}_S(f) \geq \mathcal{E}_{P^*\hat{Q}P}(f)$.
			
			For $f \in L_0^2(\pi)$,
			\[
			\begin{aligned}
				\mathcal{E}_{P^*E^*EP}(f) &= \|f\|_{\pi}^2 - \langle EPf, EPf \rangle_{\varpi} \\
				&= \|f\|_{\pi}^2 - \int_{\Z} \varpi(\df z) \left[\int_{\X} \varpi_z(\df x) \, f(x) \right]^2.
			\end{aligned}
			\]
			Apply \eqref{ine:soft-projection} with $g(z) = \varpi_z f$ together with the most recent display to obtain
			\[
			\mathcal{E}_{P^*E^*EP}(f)  \geq \|f\|_{\pi}^2 - \frac{1}{2 \bar{\kappa}} \int_{\Z^2}  (\varpi_z f - \varpi_{z'} f)^2 \, \frac{\varpi(\df z) \, \varpi(\df z')}{\chi_{\max}^2(\varpi_z || \varpi_{z'})}.
			\]
			We then use an elementary inequality from \cite{ge2018simulated}, which follows from Cauchy-Schwarz:
			\[
			\begin{aligned}
				(\varpi_z f - \varpi_{z'} f)^2 &= \left\{ \int_{\X} [f(x) - \varpi_z f] \left[ 1- \frac{\df \varpi_{z'}}{\df \varpi_z}(x) \right] \varpi_z(\df x) \right\}^2 \\
				&\leq \int_{\X} [f(x) - \varpi_z f]^2 \, \varpi_z(\df x) \, \left\| \frac{\df \varpi_{z'}}{\df \varpi_z} - 1 \right\|_{\varpi_z}^2 ,
			\end{aligned}
			\]
			$\varpi \times \varpi$-a.e. over the set of $(z,z')$ such that $\chi_{\max}^2(\varpi_z || \varpi_{z'}) < \infty$.
			Combining the two most recent displays yields
			\[
			\begin{aligned}
				\mathcal{E}_{P^*E^*EP}(f) &\geq \|f\|_{\pi}^2 - \frac{1}{2 \bar{\kappa}} \int_{\Z} \int_{\X} [f(x) - \varpi_z f]^2 \, \varpi_z(\df x) \, \varpi(\df z) \\
				&= \|f\|_{\pi}^2 - \frac{1}{2 \bar{\kappa}} (\|f\|_{\pi}^2 - \langle EPf, EPf \rangle_{\varpi}).
			\end{aligned}
			\]
			In other words,
			\[
			\mathcal{E}_{P^*E^*EP}(f) \geq \frac{\|f\|_{\pi}^2}{1 + 1/(2\bar{\kappa})}.
			\]
			Combining this with \eqref{ine:soft-pre} yields the desired result.
		\end{proof}
		
		The preceding proof clarifies the relationship between our framework and the fairly general
		Theorem~6.1 of \cite{ge2018simulated}. After replacing conditions
		\ref{H1} and \ref{H4} by the weaker condition \ref{H5},
		Proposition~\ref{pro:ge}, which is our discrete-time analogue
		of that theorem, follows in two steps. First, Corollary~\ref{cor:uniform},
		applied with $\Y=\X$ and
		\[
		Pf(x,z)=f(x),
		\]
		gives the comparison in \eqref{ine:soft-pre}. Second,
		\eqref{ine:soft-projection} supplies a lower bound on
		$\mathrm{Gap}(P^*E^*EP)$. The same canonical choice of $P$ underlies the
		derivations of Proposition~\ref{pro:madras-1}, concerning decomposition over
		an overlapping cover; Proposition~\ref{pro:andrieu}, concerning data
		augmentation with one intractable conditional; Proposition~\ref{pro:qin},
		concerning random-scan Gibbs sampling; Proposition~\ref{pro:chen},
		concerning stochastic localization; and Proposition~\ref{pro:hit},
		concerning hit-and-run.
		
		This common representation does not mean that Theorem~6.1 of
		\cite{ge2018simulated} subsumes all these results. More precisely,
		Proposition~\ref{pro:ge} combines the extension of
		Proposition~\ref{pro:andrieu} under condition \ref{H5} with the particular
		$\chi_{\max}^2$-based lower bound on
		$\mathrm{Gap}(P^*E^*EP)$ supplied by
		\eqref{ine:soft-projection}. By contrast,
		Propositions~\ref{pro:madras-1}, \ref{pro:andrieu}, \ref{pro:qin},
		\ref{pro:chen}, and \ref{pro:hit} establish analogues of
		\eqref{ine:soft-pre}, up to constant factors, without imposing
		\eqref{ine:soft-projection}. They therefore leave
		$\mathrm{Gap}(P^*E^*EP)$ to be controlled by other methods. This is
		illustrated in Section~\ref{sec:example-hybrid}, where a recent
		entropy-contraction result is used to bound this quantity for the
		hit-and-run sampler.
		For the natural line augmentation used in
		Section~\ref{ssec:hit-1} with $\ell = 1$ and $k\geq2$, the
		component measures $\varpi_z$ are supported on lines in $\mathbb{R}^k$ determined by $z = (w,x)$, where $w$ is a the direction of a line, and $x$ is a point on the line. For
		$\varpi \times \varpi$-almost every pair of $z$ and $z'$, the corresponding line
		measures are mutually singular, and hence
		$\chi_{\max}^2(\varpi_z || \varpi_{z'})=\infty$. Consequently, under the
		convention $0/0 = 0$, $1/0 = \infty$, $1/\infty=0$, the projected Dirichlet form in
		\eqref{ine:soft-projection} vanishes whenever $z \mapsto g(z)$ is a function of the line determined by~$z$. 
		Then \eqref{ine:soft-projection} cannot hold in general with $\bar{\kappa} > 0$.
		Thus,
		Proposition~\ref{pro:ge} yields no nontrivial bound through this
		natural decomposition.
		
		Proposition~\ref{pro:madras}, concerning decomposition over a partition
		$\{\X_z\}$ of the state space, has a different operator representation:
		its derivation uses
		\[
		Pf(x,z)=M^{1/2}f(x)
		\]
		rather than the canonical embedding above. Moreover, the natural choice
		of $\varpi_z$ is $\pi$ conditioned on $\X_z$. Since the sets $\X_z$ are
		disjoint, $\chi_{\max}^2(\varpi_z || \varpi_{z'})=\infty$ for every distinct
		pair $z,z'$. Thus, for this natural decomposition, all off-diagonal weights in
		\eqref{ine:soft-projection} vanish, so the projected Dirichlet form is
		zero and the inequality cannot hold with $\bar{\kappa}>0$.
		
		Finally, Theorem~6.1 does not provide the
		non-reversible and weak-Poincar\'e extensions developed in
		Sections~\ref{ssec:nonreveresible} and \ref{ssec:weak}, respectively.
		To our knowledge, neither the non-reversible comparison result in
		Section~\ref{ssec:nonreveresible} nor its application to data
		augmentation with two intractable conditionals in
		Section~\ref{ssec:doubly-0} is covered by previous spectral-gap
		decomposition results.
	}

	\bibliographystyle{plain}
	\bibliography{qinbib}

\begin{thebibliography}{10}

\bibitem{anari2021spectral}
Nima Anari, Kuikui Liu, and Shayan~Oveis Gharan.
\newblock Spectral independence in high-dimensional expanders and applications
  to the hardcore model.
\newblock {\em SIAM Journal on Computing}, pages FOCS20--1, 2021.

\bibitem{andersen2007hit}
Hans~C Andersen and Persi Diaconis.
\newblock Hit and run as a unifying device.
\newblock {\em Journal de la Societe Francaise de Statistique \& Revue de
  Statistique Appliquee}, 148(4):5--28, 2007.

\bibitem{andrieu2022comparison}
Christophe Andrieu, Anthony Lee, Sam Power, and Andi~Q Wang.
\newblock {Comparison of Markov chains via weak Poincar{\'e} inequalities with
  application to pseudo-marginal MCMC}.
\newblock {\em Annals of Statistics}, 50(6):3592--3618, 2022.

\bibitem{andrieu2024explicit}
Christophe Andrieu, Anthony Lee, Sam Power, and Andi~Q Wang.
\newblock {Explicit convergence bounds for Metropolis Markov chains:
  Isoperimetry, spectral gaps and profiles}.
\newblock {\em Annals of Applied Probability}, 34(4):4022--4071, 2024.

\bibitem{andrieu18supplement}
Christophe Andrieu, Anthony Lee, and Matti Vihola.
\newblock {Supplement to ``Uniform ergodicity of the iterated conditional SMC
  and geometric ergodicity of particle Gibbs samplers”}.
\newblock DOI: 10.3150/15-BEJ785SUPP, 2018.

\bibitem{andrieu2018uniform}
Christophe Andrieu, Anthony Lee, and Matti Vihola.
\newblock {Uniform ergodicity of the iterated conditional SMC and geometric
  ergodicity of particle Gibbs samplers}.
\newblock {\em Bernoulli}, 24:842--872, 2018.

\bibitem{ascolani2024entropy}
Filippo Ascolani, Hugo Lavenant, and Giacomo Zanella.
\newblock Entropy contraction of the {G}ibbs sampler under log-concavity.
\newblock arXiv preprint, 2026.

\bibitem{ascolani2024scalability}
Filippo Ascolani, Gareth~O Roberts, and Giacomo Zanella.
\newblock {Scalability of Metropolis-within-Gibbs schemes for high-dimensional
  Bayesian models}.
\newblock {\em Journal of the Royal Statistical Society Series B: Statistical
  Methodology}, page qkaf084, 2026.

\bibitem{atchade2021approximate}
Yves~F Atchad{\'e}.
\newblock {Approximate spectral gaps for Markov chain mixing times in high
  dimensions}.
\newblock {\em SIAM Journal on Mathematics of Data Science}, 3(3):854--872,
  2021.

\bibitem{baxendale2005renewal}
Peter~H. Baxendale.
\newblock Renewal theory and computable convergence rates for geometrically
  ergodic {M}arkov chains.
\newblock {\em Annals of Applied Probability}, 15(1B):700--738, 2005.

\bibitem{bou2025accelerated}
Nawaf Bou-Rabee, Andreas Eberle, and Stefan Oberd{\"o}rster.
\newblock {Accelerated convergence in Hit-and-Run Monte Carlo and a
  coordinate-free randomized Kaczmarz algorithm}.
\newblock {\em Electronic Journal of Probability}, 30:1--28, 2025.

\bibitem{caprio2023calculus}
Rocco Caprio and Adam~M Johansen.
\newblock {A calculus for Markov chain Monte Carlo: studying approximations in
  algorithms}.
\newblock arXiv preprint, 2023.

\bibitem{caracciolo1990nonlocal}
Sergio Caracciolo, Andrea Pelissetto, and Alan~D Sokal.
\newblock {Nonlocal Monte Carlo algorithm for self-avoiding walks with fixed
  endpoints}.
\newblock {\em Journal of Statistical Physics}, 60:1--53, 1990.

\bibitem{caracciolo1992two}
Sergio Caracciolo, Andrea Pelissetto, and Alan~D Sokal.
\newblock Two remarks on simulated tempering.
\newblock Unpublished manuscript, 1992.

\bibitem{carlen2003determination}
Eric~A Carlen, Maria~C Carvalho, and Michael Loss.
\newblock Determination of the spectral gap for {K}ac's master equation and
  related stochastic evolution.
\newblock {\em Acta Mathematica}, 191(1):1--54, 2003.

\bibitem{chen2024rapid}
Xiaoyu Chen, Weiming Feng, Yitong Yin, and Xinyuan Zhang.
\newblock {Rapid mixing of Glauber dynamics via spectral independence for all
  degrees}.
\newblock {\em SIAM Journal on Computing}, pages FOCS21--224, 2024.

\bibitem{chen2022localization}
Yuansi Chen and Ronen Eldan.
\newblock {Localization schemes: A framework for proving mixing bounds for
  Markov chains}.
\newblock In {\em 2022 IEEE 63rd Annual Symposium on Foundations of Computer
  Science}, pages 110--122. IEEE, 2022.

\bibitem{chen2021rapid}
Zongchen Chen, Andreas Galanis, Daniel {\v{S}}tefankovi{\v{c}}, and Eric
  Vigoda.
\newblock Rapid mixing for colorings via spectral independence.
\newblock In {\em Proceedings of the 2021 ACM-SIAM Symposium on Discrete
  Algorithms}, pages 1548--1557. SIAM, 2021.

\bibitem{chen2023rapid}
Zongchen Chen, Kuikui Liu, and Eric Vigoda.
\newblock {Rapid mixing of Glauber dynamics up to uniqueness via contraction}.
\newblock {\em SIAM Journal on Computing}, 52(1):196--237, 2023.

\bibitem{chewi2022query}
Sinho Chewi, Patrik~R Gerber, Chen Lu, Thibaut Le~Gouic, and Philippe Rigollet.
\newblock The query complexity of sampling from strongly log-concave
  distributions in one dimension.
\newblock In {\em Conference on Learning Theory}, pages 2041--2059. PMLR, 2022.

\bibitem{chewi2023entropic}
Sinho Chewi and Aram-Alexandre Pooladian.
\newblock {An entropic generalization of Caffarelli’s contraction theorem via
  covariance inequalities}.
\newblock {\em Comptes Rendus. Math{\'e}matique}, 361(G9):1471--1482, 2023.

\bibitem{chlebicka2023solidarity}
Iwona Chlebicka, Krzysztof {\L}atuszy{\'n}ski, and B{\l}a{\.z}ej Miasojedow.
\newblock Solidarity of {G}ibbs samplers: the spectral gap.
\newblock Technical report, 2023.

\bibitem{diaconis1993comparison}
Persi Diaconis and Laurent Saloff-Coste.
\newblock {Comparison theorems for reversible Markov chains}.
\newblock {\em Annals of Applied Probability}, 3(3):696--730, 1993.

\bibitem{diaconis1996logarithmic}
Persi Diaconis and Laurent Saloff-Coste.
\newblock {Logarithmic Sobolev inequalities for finite Markov chains}.
\newblock {\em Annals of Applied Probability}, 6(3):695--750, 1996.

\bibitem{douc2018markov}
Randal Douc, Eric Moulines, Pierre Priouret, and Philippe Soulier.
\newblock {\em Markov Chains}.
\newblock Springer, 2018.

\bibitem{eldan2013thin}
Ronen Eldan.
\newblock Thin shell implies spectral gap up to polylog via a stochastic
  localization scheme.
\newblock {\em Geometric and Functional Analysis}, 23(2):532--569, 2013.

\bibitem{feng2022rapid}
Weiming Feng, Heng Guo, Yitong Yin, and Chihao Zhang.
\newblock {Rapid mixing from spectral independence beyond the Boolean domain}.
\newblock {\em ACM Transactions on Algorithms}, 18(3):1--32, 2022.

\bibitem{gaitonde2024comparison}
Jason Gaitonde and Elchanan Mossel.
\newblock Comparison theorems for the mixing times of systematic and random
  scan dynamics.
\newblock In {\em Proceedings of the 2026 Annual ACM-SIAM Symposium on Discrete
  Algorithms (SODA)}, pages 2575--2587. SIAM, 2026.

\bibitem{gao2026weak}
Mengxi Gao, Gareth~O Roberts, and Andi~Q Wang.
\newblock {Weak Poincar\'{e} inequalities for deterministic-scan
  Metropolis-within-Gibbs samplers}.
\newblock arXiv preprint, 2026.

\bibitem{ge2018simulated}
Rong Ge, Holden Lee, and Andrej Risteski.
\newblock {Simulated tempering Langevin Monte Carlo II: An improved proof using
  soft Markov chain decomposition}.
\newblock arXiv preprint, 2018.

\bibitem{guan2007small}
Yongtao Guan and Stephen~M Krone.
\newblock {Small-world MCMC and convergence to multi-modal distributions: From
  slow mixing to fast mixing}.
\newblock {\em Annals of Applied Probability}, 17:284--304, 2007.

\bibitem{he2016scan}
Bryan~D He, Christopher~M De~Sa, Ioannis Mitliagkas, and Christopher R{\'e}.
\newblock {Scan order in Gibbs sampling: Models in which it matters and bounds
  on how much}.
\newblock In {\em Advances in neural information processing systems}, pages
  1--9, 2016.

\bibitem{helmberg2014introduction}
Gilbert Helmberg.
\newblock {\em {Introduction to Spectral Theory in Hilbert Space}}.
\newblock Elsevier, 2014.

\bibitem{hobert2008theoretical}
James~P Hobert and Dobrin Marchev.
\newblock A theoretical comparison of the data augmentation, marginal
  augmentation and {PX}-{DA} algorithms.
\newblock {\em Annals of Statistics}, 36(2):532--554, 2008.

\bibitem{jerrum2004elementary}
Mark Jerrum, Jung-Bae Son, Prasad Tetali, and Eric Vigoda.
\newblock {Elementary bounds on Poincar{\'e} and log-Sobolev constants for
  decomposable Markov chains}.
\newblock {\em Annals of Applied Probability}, 14:1741--1765, 2004.

\bibitem{jones2001honest}
Galin~L Jones and James~P Hobert.
\newblock Honest exploration of intractable probability distributions via
  {M}arkov chain {M}onte {C}arlo.
\newblock {\em Statistical Science}, 16(4):312--334, 2001.

\bibitem{jones2014convergence}
Galin~L Jones, Gareth~O Roberts, and Jeffrey~S Rosenthal.
\newblock Convergence of conditional {M}etropolis-{H}astings samplers.
\newblock {\em Advances in Applied Probability}, 46(2):422--445, 2014.

\bibitem{latuszynski2014convergence}
Krzysztof {\L}atuszy{\'n}ski and Daniel Rudolf.
\newblock Convergence of hybrid slice sampling via spectral gap.
\newblock {\em Advances in Applied Probability}, 56:1440--1466, 2024.

\bibitem{lawler1988bounds}
Gregory~F Lawler and Alan~D Sokal.
\newblock Bounds on the $l^2$ spectrum for {M}arkov chains and {M}arkov
  processes: {A} generalization of {C}heeger’s inequality.
\newblock {\em Transactions of the American Mathematical Society},
  309(2):557--580, 1988.

\bibitem{lee2024fast}
Holden Lee and Kexin Zhang.
\newblock {Fast mixing of data augmentation algorithms: Bayesian probit, logit
  and lasso regression}.
\newblock {\em The Annals of Statistics}, 54(4):2083--2109, 2026.

\bibitem{liu1995covariance}
Jun~S Liu, Wing~H Wong, and Augustine Kong.
\newblock Covariance structure and convergence rate of the {G}ibbs sampler with
  various scans.
\newblock {\em Journal of the Royal Statistical Society, Series B},
  57(1):157--169, 1995.

\bibitem{liu2023spectral}
Kuikui Liu.
\newblock {\em Spectral Independence: A New Tool to Analyze Markov Chains}.
\newblock University of Washington, 2023.

\bibitem{liu2026spectral}
Shuigen Liu and Xin~T Tong.
\newblock {Spectral gap of Metropolis algorithms for non-smooth distributions
  under isoperimetry}.
\newblock {\em Japanese Journal of Statistics and Data Science}, pages 1--27,
  2026.

\bibitem{lovasz1999hit}
L{\'a}szl{\'o} Lov{\'a}sz.
\newblock Hit-and-run mixes fast.
\newblock {\em Mathematical Programming}, 86:443--461, 1999.

\bibitem{lovasz2004hit}
L{\'a}szl{\'o} Lov{\'a}sz and Santosh Vempala.
\newblock Hit-and-run from a corner.
\newblock In {\em Proceedings of the Thirty-sixth Annual ACM Symposium on
  Theory of Computing}, pages 310--314, 2004.

\bibitem{madras1999importance}
Neal Madras and Mauro Piccioni.
\newblock Importance sampling for families of distributions.
\newblock {\em Annals of Applied Probability}, 9(4):1202--1225, 1999.

\bibitem{madras1996factoring}
Neal Madras and Dana Randall.
\newblock Factoring graphs to bound mixing rates.
\newblock In {\em Proceedings of 37th Conference on Foundations of Computer
  Science}, pages 194--203. IEEE, 1996.

\bibitem{madras2002markov}
Neal Madras and Dana Randall.
\newblock Markov chain decomposition for convergence rate analysis.
\newblock {\em Annals of Applied Probability}, pages 581--606, 2002.

\bibitem{peskun1973optimum}
Peter~H Peskun.
\newblock {Optimum Monte-Carlo sampling using Markov chains}.
\newblock {\em Biometrika}, 60(3):607--612, 1973.

\bibitem{pillai2014ergodicity}
Natesh~S Pillai and Aaron Smith.
\newblock {Ergodicity of approximate MCMC chains with applications to large
  data sets}.
\newblock arXiv preprint, 2014.

\bibitem{popoviciu1935equations}
Tiberiu Popoviciu.
\newblock Sur les {\'e}quations alg{\'e}briques ayant toutes leurs racines
  r{\'e}elles.
\newblock {\em Mathematica}, 9(129-145):20, 1935.

\bibitem{power2024weak}
Sam Power, Daniel Rudolf, Bj{\"o}rn Sprungk, and Andi~Q Wang.
\newblock {Weak Poincar{\'e} inequality comparisons for ideal and hybrid slice
  sampling}.
\newblock arXiv preprint, 2024.

\bibitem{pozza2024fundamental}
Francesco Pozza and Giacomo Zanella.
\newblock {On the fundamental limitations of multi-proposal Markov chain Monte
  Carlo algorithms}.
\newblock {\em Biometrika}, 112(2):asaf019, 2025.

\bibitem{qin2023geometric}
Qian Qin.
\newblock {Geometric ergodicity of trans-dimensional Markov chain Monte Carlo
  algorithms}.
\newblock {\em Journal of the American Statistical Association}, 2025.

\bibitem{qin2026simulation}
Qian Qin.
\newblock {Simulation code for Section 7.2 of ``On spectral gap decomposition
  for Markov chains''}.
\newblock GitHub repository,
  \url{https://github.com/qianqin11/on-spectral-gap-decomposition-for-Markov-chains},
  2026.

\bibitem{qin2020convergence}
Qian Qin and Galin~L Jones.
\newblock {Convergence Rates of Two-Component MCMC Samplers}.
\newblock {\em Bernoulli}, 28:859--885, 2022.

\bibitem{qin2024spectral}
Qian Qin, Nianqiao Ju, and Guanyang Wang.
\newblock Spectral gap bounds for reversible hybrid {G}ibbs chains.
\newblock {\em Annals of Statistics}, 53:1613--1638, 2025.

\bibitem{qin2024spectrala}
Qian Qin and Guanyang Wang.
\newblock {Spectral telescope: Convergence rate bounds for random-scan Gibbs
  samplers based on a hierarchical structure}.
\newblock {\em Annals of Applied Probability}, 34(1B):1319--1349, 2024.

\bibitem{roberts1997geometric}
Gareth~O Roberts and Jeffrey~S Rosenthal.
\newblock Geometric ergodicity and hybrid {M}arkov chains.
\newblock {\em Electronic Communications in Probability}, 2(2):13--25, 1997.

\bibitem{rudolf2018perturbation}
Daniel Rudolf and Nikolaus Schweizer.
\newblock {Perturbation theory for Markov chains via Wasserstein distance}.
\newblock {\em Bernoulli}, 24(4A):2610--2639, 2018.

\bibitem{rudolf2018comparison}
Daniel Rudolf and Mario Ullrich.
\newblock Comparison of hit-and-run, slice sampler and random walk
  {M}etropolis.
\newblock {\em Journal of Applied Probability}, 55(4):1186--1202, 2018.

\bibitem{secchi2025spectral}
Cecilia Secchi and Giacomo Zanella.
\newblock {Spectral gap of Metropolis-within-Gibbs under log-concavity}.
\newblock arXiv preprint, 2025.

\bibitem{tanner1987calculation}
Martin~A Tanner and Wing~Hung Wong.
\newblock The calculation of posterior distributions by data augmentation (with
  discussion).
\newblock {\em Journal of the American Statistical Association},
  82(398):528--540, 1987.

\bibitem{tierney1998note}
Luke Tierney.
\newblock {A note on Metropolis-Hastings kernels for general state spaces}.
\newblock {\em Annals of Applied Probability}, pages 1--9, 1998.

\bibitem{wang2022theoretical}
Guanyang Wang.
\newblock On the theoretical properties of the exchange algorithm.
\newblock {\em Bernoulli}, 28(3):1935--1960, 2022.

\bibitem{woodard2009conditions}
Dawn~B Woodard, Scott~C Schmidler, and Mark Huber.
\newblock Conditions for rapid mixing of parallel and simulated tempering on
  multimodal distributions.
\newblock {\em Annals of Applied Probability}, 19(1):617--640, 2009.

\end{thebibliography}
	
\end{document}